\documentclass[11pt,a4paper]{article}

\usepackage[margin=2.5cm]{geometry}
\usepackage{graphicx}
\usepackage{float}
\usepackage{amsmath}
\usepackage{amssymb}
\usepackage{authblk}
\usepackage{xcolor}
\usepackage{textcomp}
\usepackage{booktabs}
\usepackage{tabularx}
\usepackage{makecell}
\usepackage{xltabular}
\usepackage{caption}
\usepackage{array}
\usepackage{hyperref}
\usepackage[numbers,sort&compress]{natbib}
\usepackage[version=4,arrows=pgf-filled,textfontname=sffamily,mathfontname=mathsf]{mhchem}
\usepackage{setspace}
\usepackage{lineno}
\renewcommand{\vec}[1]{\boldsymbol{#1}}

\DeclareUnicodeCharacter{03BC}{\ensuremath{\mu}}

\usepackage[T1]{fontenc}
\usepackage{lmodern}
\usepackage[varqu,scaled=0.95]{zi4}
\usepackage{upgreek}
\usepackage{listings}
\definecolor{juliabackground}{RGB}{246,246,246}
\definecolor{juliaframe}{RGB}{200,200,200}
\definecolor{juliakeyword}{RGB}{0,72,150}
\definecolor{juliamacro}{RGB}{140,0,140}
\definecolor{juliastring}{RGB}{160,40,40}
\definecolor{juliacomment}{RGB}{100,100,100}
\lstdefinelanguage{Julia}{
    morekeywords={begin,end,function,return,for,in,if,else,elseif,while,using,import,let,do,true,false,nothing,const,local,global,struct,mutable,module,export,try,catch,finally,break,continue,quote,macro},
    morekeywords=[2]{@reaction_network,@species,@parameters,@variables,@SymbolicNeuralNetwork},
    alsoletter={@},
    sensitive=true,
    morecomment=[l]{\#},
    morecomment=[n]{\#=}{=\#},
    morestring=[b]",
}
\lstdefinestyle{julia}{
    language=Julia,
    basicstyle=\ttfamily\small,
    keywordstyle=\color{juliakeyword}\bfseries,
    keywordstyle=[2]\color{juliamacro}\bfseries,
    stringstyle=\color{juliastring},
    commentstyle=\color{juliacomment},
    backgroundcolor=\color{juliabackground},
    frame=single,
    rulecolor=\color{juliaframe},
    framesep=4pt,
    xleftmargin=5pt,
    xrightmargin=5pt,
    columns=fullflexible,
    keepspaces=true,
    showstringspaces=false,
    upquote=true,
    breaklines=true,
    aboveskip=8pt,
    belowskip=8pt,
    extendedchars=true,
    literate={μ}{{$\upmu$}}1 {α}{{$\upalpha$}}1 {β}{{$\upbeta$}}1 {γ}{{$\upgamma$}}1 {δ}{{$\updelta$}}1 {θ}{{$\uptheta$}}1 {₀}{{$_0$}}1 {₁}{{$_1$}}1 {₂}{{$_2$}}1 {₃}{{$_3$}}1,
}
\lstdefinestyle{juliaoutput}{
    style=julia,
    language={},
    basicstyle=\ttfamily\small\color{black!75},
    backgroundcolor=\color{white},
    frame=leftline,
    framerule=1.5pt,
}
\lstnewenvironment{juliacode}[1][]{\lstset{style=julia,#1}}{}
\lstnewenvironment{juliaoutput}[1][]{\lstset{style=juliaoutput,#1}}{}

\title{\bf Functional and parametric identifiability for \\ universal differential equations}
\date{}
\author[1]{Torkel E Loman}
\author[1]{Ruth E Baker}
\affil[1]{Mathematical Institute, University of Oxford, Oxford, OX2 6GG, United Kingdom}

\begin{document}

\maketitle


\begin{abstract}
Mechanistic knowledge of a studied system can be encoded as predictive mathematical models, but this raises the question of what to do when that knowledge is uncertain or incomplete. Hybrid mechanistic/data-driven modelling, where unknown system dynamics are learnt from data using machine learning components, offers a potential solution. This work focuses on one such approach, \emph{universal differential equations} (UDEs), and asks when these can, and cannot, learn unknown dynamics from data. We investigate this question through the lens of identifiability. Specifically, we show that UDE identifiability can be decomposed into \emph{parametric} and \emph{functional} identifiability, which concern the model's mechanistic and data-driven parts, respectively. We consider these two types of identifiability along a second, orthogonal, division: \emph{structural} and \emph{practical} identifiability, which distinguish ambiguity inherent to the model formulation from ambiguity arising under finite, noisy observations. Specifically, we consider the problem of how to assess practical identifiability for UDEs. As a modelling setting, we use ordinary differential equations derived from chemical reaction network models. Such models arise naturally across systems biology, chemistry, pharmacology, epidemiology, and population dynamics, making them a broadly relevant test case. Through a wide range of examples, we show that chemical reaction networks are a promising setting for UDEs, yielding learnt dynamics that are highly interpretable. We then examine how model properties---including the number of neural networks, the constraints imposed on them, and data availability---affect identifiability. A key finding is that generalising a fully mechanistic model to a UDE can sometimes have remarkably little impact on its identifiability. We demonstrate how this makes UDEs a principled way to counteract the erroneous predictions that model misspecification can cause. Taken together, our results provide strong evidence that UDEs offer an intuitive means of integrating data-driven elements into mechanistic models while preserving predictive capability.
\end{abstract}


\section{Introduction}\label{section:introduction}

Mechanistic modelling seeks to describe system dynamics by translating knowledge of the underlying processes into mathematical models, often composed of systems of differential equations. Such models are widely used to predict future system behaviour (e.g.~epidemiological forecasting) and investigate the consequences of perturbations (e.g.~gene knockdowns). Their construction, however, depends on sufficient knowledge of the underlying mechanisms. When this knowledge is incomplete, purely data-driven approaches can instead be used to infer relationships directly from observations. Advances in both large-scale data collection and machine learning have made such approaches increasingly attractive. However, purely data-driven models typically neglect established mechanistic knowledge, despite the fact that incorporating such information can substantially improve predictive performance~\cite{baker_mechanistic_2018,noordijk_rise_2024}. This has motivated growing interest in hybrid modelling approaches that combine mechanistic descriptions with data-driven components. One such framework is provided by \emph{universal differential equations} (UDEs), which replace selected functional components within a differential equation model with flexible representations inferred from data~\cite{rackauckas_universal_2021}.

For most purposes UDEs function like standard differential equations, however, they contain at least one universal function approximator---a generic and parameterised function $U(\vec{X}; \vec{\theta})$ where $U:\mathbb{R}^n\to\mathbb{R}$ is such that, if the parameterisation is made arbitrarily large, it can approximate any function $f:\mathbb{R}^n\to\mathbb{R}$ to any degree of accuracy. During model fitting, the neural network parameters $\vec{\theta}$ are estimated simultaneously with the mechanistic model parameters. Through this process, the UDE learns unknown functional forms. While many neural network applications typically involve large networks trained on huge datasets, UDE-based neural networks are typically small in size. Furthermore, while neural networks are a common function approximator choice, alternatives, such as Gaussian processes or splines are also possible. UDEs have now been applied across diverse scientific domains~\cite{bolibar_universal_2023,buckner_recovering_2024,martinez_comparative_2024,elgazzar2025}, differential equation classes~\cite{jia_neural_2020,oh_stable_2025,jo_neural_2025,zhu_neural_2021,loman_learning_2026}, and specialised modelling settings~\cite{dandekar_bayesian_2022,philipps_non-negative_2024,de_rooij_conditional_2025,schmid_assessment_2025}. We note that UDEs differ conceptually from physics-informed neural networks (PINNs)~\cite{raissi_physics-informed_2019}. While PINNs typically use neural networks to approximate latent solutions while enforcing known governing equations as a soft constraint, UDEs embed neural networks within the governing equations themselves to represent unknown model components, with the resulting dynamics obtained through numerical integration (i.e.~the governing equations are enforced as a hard constraint). 

The use of UDEs to recover unknown functional relationships raises a fundamental question: when can the recovered functions be interpreted as providing insight into previously unknown processes, rather than merely representing one of many possible explanations of the data? This depends on whether the functions are uniquely determined by the available observations, a question of \emph{identifiability}. Classical identifiability theory addresses whether model parameters can be uniquely recovered from observations. In UDEs, however, the unknown quantities include not only scalar mechanistic parameters but also functional components. Applying conventional parameter identifiability methods directly to the weights of a neural network is uninformative: different network parameterisations can represent exactly the same function, so the network weights may be non-identifiable even when the function itself is uniquely determined. To assess whether recovered functions can be meaningfully interpreted, identifiability must therefore be considered at the level of the function itself.

Classical parameter identifiability distinguishes between two complementary concepts. Structural identifiability asks whether parameters can, in principle, be uniquely recovered from ideal, noise-free observations, whereas practical identifiability asks whether they can be reliably estimated from the finite, noisy data available in practice. The corresponding question for state variables, namely whether unmeasured species can be inferred from measured ones, is known as \emph{observability}~\cite{villaverde_observability_2019}. With the introduction of functional identifiability, identifiability can be divided along two complementary axes: parametric versus functional, and structural versus practical. Structural parametric identifiability can be assessed using differential algebra approaches~\cite{LJUNG1994265}, and there now also exists software that automates this process~\cite{structidjl}. In this work, we briefly introduce structural functional identifiability (see our recent results in~\cite{loman2026structuralfunctionalidentifiabilitymodel} for more details). Next, we consider practical identifiability in the context of UDEs. We show how to assess practical parameter identifiability using established approaches such as profile likelihood analysis~\cite{kreutz_profile_2013,simpson2026parameter}. For practical functional identifiability, we propose an ensemble-based diagnostic in which variability among independently fitted functions is used to assess how strongly the available data constrain the recovered functional forms.

Throughout this work, we develop and illustrate our UDE identifiability framework using ordinary differential equation (ODE) models derived from chemical reaction networks~\cite{gunawardena_chemical_nodate,hahl_comparison_2016,feinberg_foundations_2019,loman_2022}. Chemical reaction networks arise naturally across systems biology~\cite{goodwin_oscillatory_1965}, epidemiology~\cite{avram_advancing_2024,brauer_compartmental_2008}, synthetic biology~\cite{elowitz_synthetic_2000}, pharmacology~\cite{zou_application_2020}, population dynamics~\cite{boros_center_2018}, and chemistry~\cite{field_oscillations_1974,warnatz_combustion_1999}. To simplify complex and/or unknown mechanisms in such models, many reactions are associated with a non-constant reaction rate i.e.~the reaction rate is a function that depends on some of the reactants. Hill functions, originally developed to capture cooperative ligand-receptor binding, are a widely used example of this. Such widespread use of non-constant reaction rates, often associated with unknown mechanisms, makes chemical reaction network models a natural use case for UDEs. In this work, we specifically demonstrate how phenomena such as gene transcription, disease infection events, and competition-based mortality in population dynamics can be straightforwardly modelled using UDEs.

Across a range of example chemical reaction network models, we analyse how both parametric and functional identifiability depend on model structure and data. Although the relationship is nuanced, several results emerge. First, incorporating prior knowledge, such as monotonicity constraints, directly into neural network architectures can substantially improve identifiability. Second, extending ODE-based models to UDE models often, but not always, has limited impact on parametric identifiability, and in many cases does not degrade the functional identifiability of other model components. Finally, we show how representing unknown mechanisms as learnt functions---rather than forcing simplified parametric forms---can offer robustness to model misspecification errors. Taken together, these results indicate that UDEs can be used to infer unknown mechanisms without necessarily compromising identifiability, provided that appropriate structure and constraints are imposed.


\section{Methodologies}\label{section:methods}

In this section we provide an overview of the main methodologies used in this work, including chemical reaction network modelling, synthetic data generation, UDE fitting, and practical identifiability analysis.


\subsection{Chemical reaction network modelling}\label{section:methods_chemical_reaction_networks}

Throughout this work, we represent dynamical systems as chemical reaction networks and from these generate ODE models using standard formalisms. A chemical reaction network consists of $M$ species, denoted $S_1,S_2,\dots,S_M$, and $K$ reactions, denoted $R_1,R_2,\dots,R_K$. The system state is given by $\vec{X}(t)=[X_1(t),X_2(t),\dots,X_M(t)]$, where $X_i(t)$ denotes the concentration of species $S_i$ at time $t$.

Each reaction $R_k$ can be written as
\begin{equation}
\ce{$\alpha_1^k S_1 + \alpha_2^k S_2 + \cdots + \alpha_M^k S_M$ ->C[$c_k(\vec{X};\vec{p}_k)$] $\beta_1^k S_1 + \beta_2^k S_2 + \cdots + \beta_M^k S_M$,}
\end{equation}
where $\vec{\alpha}^{\,k}=[\alpha_1^k,\ldots,\alpha_M^k]$ and $\vec{\beta}^{\,k}=[\beta_1^k,\ldots,\beta_M^k]$ are the substrate and product stoichiometric vectors, respectively, and $c_k(\vec{X};\vec{p}_k)\geq0$ is the corresponding reaction-rate function\footnote{For simplicity, we assume throughout that reaction rates are independent of time.}. Defining the net stoichiometric vector as $\vec{\nu}^{\,k}=\vec{\beta}^{\,k}-\vec{\alpha}^{\,k}$, the deterministic ODE model is generated according to the reaction-rate equations:
\begin{equation}
\label{eq:RRE}
\frac{\mathrm{d}X_m}{\mathrm{d}t}
=
\sum_{k=1}^{K}\nu_m^k\,a_k(\vec{X};\vec{p}_k),
\qquad m=1,\ldots,M,
\end{equation}
where the reaction propensity $a_k$ is defined as
\begin{equation*}
a_k(\vec{X};\vec{p}_k) = \prod_{m=1}^M\frac{X_m^{\alpha_m^k}}{\alpha_m^k!}c_k(\vec{X};\vec{p}_k).
\end{equation*}

As a simple example, consider a reversible dimerisation system
\begin{equation}
\ce{$2X_1$ ->C[$k_b$] $X_2$},
\qquad
\ce{$X_2$ ->C[$k_d$] $2X_1$,}
\end{equation}
which gives the following ODE model:
\begin{align}
\frac{\mathrm{d}X_1}{\mathrm{d}t}
&=
2k_dX_2-k_bX_1^2,
\\
\frac{\mathrm{d}X_2}{\mathrm{d}t}
&=
\frac{k_b}{2}X_1^2-k_dX_2.
\end{align}

In many applications, individual reaction mechanisms are too complex or insufficiently understood to be modelled directly. Instead, they are represented using effective state-dependent reaction rates that capture their aggregate behaviour. A common example arises in gene regulation, where production rates are frequently represented using Hill functions. A simple self-activation loop (where a single species activates its own production), for example, might be described using
\begin{equation}
\ce{\text{\o} ->C[$v\dfrac{X_1^n}{X_1^n+K^n}$] $X_1$,}
\qquad
\ce{$X_1$ ->C[$d$] \text{\o},}
\end{equation}
for which the corresponding ODE model is
\begin{equation}
\frac{\mathrm{d}X_1}{\mathrm{d}t}
=
v\frac{X_1^n}{X_1^n+K^n}
-dX_1.
\end{equation}
While prescribed functional forms such as Hill functions often provide convenient and biologically plausible modelling assumptions, they are ultimately heuristic choices whose validity depends on the application. UDEs provide a principled means of relaxing these assumptions by replacing prescribed reaction-rate functions with universal function approximators whose functional forms are recovered from data.


\subsection{Synthetic data generation}\label{section:detailed_methodology_synthetic_data}

All analyses in this paper were carried out using data synthetically generated from ground-truth models. Each model originates from a chemical reaction network modelled using the form in Equation~\eqref{eq:RRE}. The model is simulated numerically, using a specified initial condition, over the time interval $t\in[0,t_\text{end}]$. We assume that observations of the system are collected at $N$ time points $t_1,\ldots,t_N$ equally distributed in the interval $[0,t_\text{end}]$. For each measured species, $X_j$, $j=1,\ldots,M_\text{obs}$, we assume a measurement noise model such that the observed quantity, $\hat{X}_j$, is given by
\begin{equation} \label{eq:measurement_formula}
\hat{X}_j(t_i) \sim \mathcal{N}(X_j(t_i),(\sigma X_j(t_i))^2),
\end{equation}
for $i=1,\ldots,N$ and $j=1,\ldots,M_\text{obs}$, where $\mathcal{N}(\mu,\nu^2)$ is the normal distribution with mean $\mu$ and variance $\nu^2$ and $\sigma>0$ is a constant parameter that controls the noise level. For simplicity, throughout this work we assume that $\sigma$ is known \textit{a priori}. All relevant details for the synthetic data generation settings of each example model, such as parameters and initial conditions, can be found in Supplementary Section~\ref{ssection:sup_tables}.


\subsection{Universal differential equation fitting}\label{section:detailed_methodology_UDE_fitting}

Throughout this work, we use small feedforward neural networks (at most 12 hidden units) as universal function approximators. Small networks reduce computational cost, and previous studies suggest that network architecture has only a limited influence on UDE performance~\cite{philipps_current_2025}. Unless otherwise stated, all activation functions are softplus~\cite{dugas_incorporating_2000} (the networks of the incoherent feedforward loop use a fast sigmoid in their first layer), whose smoothness is advantageous for ODE solvers. Specifically, by using non-negative activation functions (like softplus) in all neural network output layers, we ensure that $U(\vec{X};\vec{\theta})\geq 0$, consistent with the non-negative reaction rates required by chemical reaction network models. In Section~\ref{section:results_linear_pathway}, we consider modified architectures that additionally enforce structural constraints such as monotonicity and boundedness. We apply no explicit regularisation to the fitted functions, although using relatively small neural network architectures, which limits the complexity of the learnt functions, has a similar effect. Remaining permissive in this way lets us explore, as extensively as possible, the range of functions that can explain the data, with as few prior assumption about how complex these functions can, and cannot, be.

We fit the mechanistic and neural network parameters of each UDE simultaneously. We do so by maximising the log-likelihood of the data, under the assumption of the noise model described in Section~\ref{section:detailed_methodology_synthetic_data} with $\sigma$ known. For mechanistic parameters $\vec{p}$ and neural network parameters $\vec{\theta}$, this generates the following log-likelihood:
\begin{equation} \label{eq:likelihood_formula}
    \ell(\vec{p},\vec{\theta}|\hat{\vec{X}}) = -\sum_{j = 1}^{M_\text{obs}} \sum_{i = 1}^N \left[ \log(\sigma \, X_j(t_i)) + \frac{1}{2} \left( \frac{X_j(t_i) - \hat{X}_j(t_i)}{\sigma \, X_j(t_i)} \right)^2 \right],
\end{equation}
where $t_1,\ldots,t_N$ are the sample time points, $X_1(t_i),\ldots,X_{M_\text{obs}}(t_i)$ are the measured species' simulated values at these time points (which depend on $\vec{p}$ and $\vec{\theta}$), $\hat{X}_1(t_i),\ldots,\hat{X}_{M_\text{obs}}(t_i)$ are the corresponding measured values, summarised as $\hat{\vec{X}}$, and $\sigma$ the noise parameter. In Supplementary Section~\ref{ssection:noise_parameter_learning}, we also demonstrate how an unknown $\sigma$ can be learnt together with the unknown mechanistic parameters and functions.

In each case, we use multistart optimisation with between $100$ and $400$ starts (exact numbers vary from case to case and are detailed in the code). Initial values for mechanistic parameters are generated using Latin hypercube sampling, while neural network parameters are initialised using the default initialisation provided by the Lux.jl package. We start with between $10{,}000$ and $100{,}000$ iterations using Adam (exact numbers vary from case to case and are detailed in the code), followed by L-BFGS until convergence is achieved (using the standard tolerance settings)~\cite{kingma_adam_2017,liu_limited_1989}. Alternative optimisation approaches can be used, and are considered in Supplementary Section~\ref{ssection:optimisation_approaches}. Mechanistic parameters are estimated on the logarithmic scale, improving numerical stability and optimisation performance for ODE models~\cite{schmiester_petabinteroperable_2021}.

All models are implemented in Julia. Neural networks are constructed using Lux.jl, chemical reaction networks are specified using Catalyst.jl, and UDEs are generated via ModelingToolkitNeuralNets.jl~\cite{pal_lux_2022,loman_catalyst_2023,noauthor_scimlmodelingtoolkitneuralnetsjl_2025}. This combination of packages enables a notably straightforward UDE implementation interface, which is exemplified in Supplementary Section~\ref{ssection:code_example}. Optimisation is performed using Optimization.jl (for Adam) and Optim.jl (for L-BFGS)~\cite{mogensen_optim_2018,dixit2023optimization}, and systems of ODEs are solved using OrdinaryDiffEq.jl~\cite{rackauckas_differentialequationsjl_2017}.


\subsection{Practical identifiability analysis}\label{section:detailed_methodology_identifiability_analysis}

\begin{figure}[htbp]
     \centering
     \includegraphics[width=0.8\linewidth]{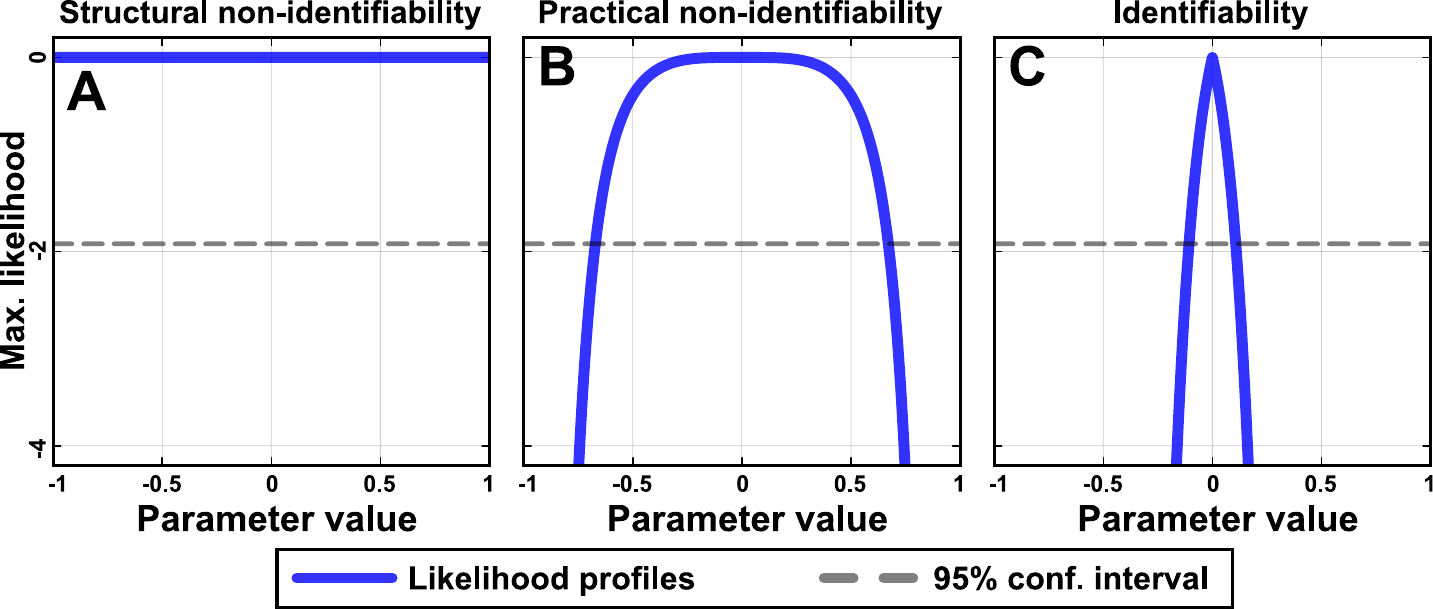}
     \caption{\textbf{Profile likelihood analysis for parameter identifiability.} For each parameter, the profile likelihood shows how the optimal likelihood varies as that parameter is fixed across a range of values, while all other parameters are re-optimised~\cite{kreutz_profile_2013}. The $x$-axis denotes the parameter value and the $y$-axis the normalised likelihood. (A) A flat profile indicates structural non-identifiability, as the parameter can vary without affecting the fit. (B) A broad plateau indicates practical non-identifiability, where a wide range of values yields similarly good fits. (C) A sharp maximum indicates identifiability, as only values near the optimum are consistent with the data. For maximum likelihood estimation, normalising the profile by subtracting its maximum value allows a horizontal threshold at $y \approx -1.92$ to define a $95\%$ confidence interval for a parameter's value. Smaller intervals correspond to higher practical identifiability.} 
     \label{figure:profile_likelihood_diagram_intro}
 \end{figure}

We use profile likelihood analysis to determine practical parametric identifiability across all models. In particular, for a UDE-based model depending on the unknown mechanistic parameters $\vec{p} = (p_1, p_2, \ldots, p_n)$ and neural network parameters $\vec{\theta}$, together with data $\hat{\vec{X}}$, then for each mechanistic parameter, $p_i$, we define a corresponding profile likelihood, $l_i$, as
\begin{equation}
    l_i(p_i|\hat{\vec{X}}) = \underset{\vec{p}_{-i},\vec{\theta}}{\sup} \left[\ell(p_i,\vec{p}_{-i},\vec{\theta} | \hat{\vec{X}})\right],
\end{equation}
where $\vec{p}_{-i}$ is $\vec{p}$ with the $i^\text{th}$ element removed. This augments the definition of a profile likelihood that is used for ODE-based models in the natural way by additionally optimising over neural network parameters $\vec{\theta}$. The shape of the resulting profile can be used as a measure of identifiability (Figure \ref{figure:profile_likelihood_diagram_intro}). In particular, if we define the normalised likelihood profile as
\begin{equation*}
    \hat{l}_i(p_i|\hat{\vec{X}}) = l_i(p_i|\hat{\vec{X}}) - \underset{\vec{p},\vec{\theta}}{\sup} \left[\ell(\vec{p},\vec{\theta} | \hat{\vec{X}})\right],
\end{equation*}
then an approximate $95\%$ likelihood-based confidence interval can be obtained from the intersections of the normalised profile with a horizontal line at $y\approx -1.92$~\cite{simpson2026parameter}. We use the length of this interval as a measure of parameter identifiability. Accurately computing these profiles can be substantially more computationally expensive for UDEs compared to ODE-based models.

For practical functional identifiability in UDE-based models, we generated ensembles of independently fitted models using different random initialisations. Typically, between $100$ and $200$ successful fits were retained for each analysis, although the exact number varied between cases (see code for details). Because UDE training constitutes a non-convex optimisation problem, some runs can converge to poor local optima. To exclude such optimisation failures from the functional ensemble, we required each retained fit to achieve a log-likelihood at least as high as that obtained using the ground-truth parameter set,
\begin{equation}
\label{eq:likelihood_threshold}
    \ell(\vec{p}_{\mathrm{fit}},\vec{\theta}_{\mathrm{fit}}|\hat{\vec{X}}) > \ell(\vec{p}_{\mathrm{true}}|\hat{\vec{X}}),
\end{equation}
where $\vec{p}_{\mathrm{fit}}$ and $\vec{\theta}_{\mathrm{fit}}$ denote the fitted parameters and $\vec{p}_{\mathrm{true}}$ the ground-truth parameters (see Equation~\eqref{eq:likelihood_formula}). Runs that did not satisfy this criterion were repeated from new initialisations until the required number of successful fits was obtained. To facilitate comparison with the ODE-based models, we perform the same procedures to assess practical parametric and functional identifiability in the case where only mechanistic parameters are estimated.

The purpose of the likelihood thresholding in Equation~\eqref{eq:likelihood_threshold} is to distinguish poor optimisation outcomes from alternative well-fitting solutions. In this work, all non-misspecified models should achieve this threshold when fitted to the synthetic datasets using any adequate optimisation procedure. Each retained fit therefore provides a well-fitting candidate function that is consistent with the observed data under the chosen model. We use variability across this ensemble as an empirical diagnostic of practical functional identifiability. Small variability indicates that the optimisation procedure repeatedly identifies similar functions among well-fitting solutions, whereas large variability provides direct evidence that the available data do not uniquely constrain the function. Importantly, small ensemble variability does not by itself establish functional identifiability, since alternative well-fitting functions may remain inaccessible to the chosen architecture or optimisation procedure. Throughout most of the work variability is assessed visually from ensemble plots. In Section~\ref{section:results_mosquitoes}, we additionally quantify it using the mean $L^2$ distance between functions in the ensemble (computed numerically using a grid of 200 equally-spaced grid points). Supplementary Table~\ref{stable:figure_ensemble_mean_l2_distances} also reports the mean $L^2$ distance for every functional ensemble plot considered in the main text.

The likelihood of the ground-truth parameter set is used here only as a convenient optimisation-quality threshold and is available because all analyses use synthetic data. However, it is not a formal confidence threshold for functional uncertainty. For experimental data, where the ground-truth parameters are unknown, an alternative criterion would be required, for example one based on the maximum likelihood achieved across optimisation runs. Likewise, rerunning fits that fail to reach the threshold is not an intrinsic part of the ensemble-based identifiability procedure; it is used here simply to ensure that the reported functional variability reflects differences among well-fitting solutions rather than failures of numerical optimisation. In Supplementary Section~\ref{ssection:ensemble_plot_thresholds}, we demonstrate the ensemble analysis using an alternative threshold, which does not require knowledge of the ground-truth model.


\section{Results}\label{section:results}

In the following sections, we first introduce UDEs in the context of chemical reaction networks, and then apply them to a series of models to investigate how generalising mechanistic models to UDEs affects parametric and functional identifiability. We examine structural and practical non-identifiability, the relationship between identifiability and predictive performance, the effects of learning multiple functions and imposing functional constraints, and the trade-off between identifiability and model misspecification. Supplementary analyses using additional models assess the generality of these findings.


\subsection{Universal differential equations for chemical reaction network models}\label{section:results_crn_udes}

Section~\ref{section:methods_chemical_reaction_networks} illustrated two types of reaction rates: constant rates, corresponding to well-characterised mechanisms, and state-dependent rates, which are often introduced to represent complex or incompletely understood processes. In practice, such state-dependent reaction rates are frequently chosen from convenient parametric families, for example Hill functions, despite limited mechanistic justification. UDEs provide a principled alternative by replacing selected reaction-rate functions with universal function approximators, whose functional forms are subsequently recovered from data. This yields a natural class of hybrid models in which established mechanisms are retained while unknown mechanisms are inferred.

Formally, consider a chemical reaction with reaction rate $c(\vec{X};\vec{p})$, depending on the system state $\vec{X}$ and parameters $\vec{p}$. In a UDE formulation, the prescribed reaction-rate function is replaced by a universal function approximator $U(\vec{X};\vec{\theta})$~\cite{loman_catalyst_2023}. Since $U$ behaves as an ordinary (albeit highly parameterised) function, the resulting model remains a system of ODEs derived from the underlying chemical reaction network. UDEs can therefore be viewed as a direct generalisation of chemical reaction network models in which selected reaction-rate functions are promoted from prescribed functional forms to data-driven representations.

For example, replacing the production rate in the self-activation loop from Section~\ref{section:methods_chemical_reaction_networks} with a universal function approximator yields
\begin{equation}
\ce{\text{\o} ->C[$U(X;\vec{\theta})$] $X$,} 
\qquad 
\ce{$X$ ->C[$d$] \text{\o},}
\end{equation}
which yields the following UDE model
\begin{equation}
\frac{\mathrm{d}X}{\text{d}t} = U(X;\vec{\theta}) - dX.
\end{equation}


\subsection{A simple self-activation loop exhibits structural functional non-identifiability}\label{section:results_simple_sa}

We first investigate a simple setting in which replacing a mechanistic reaction-rate function with a universal function approximator introduces structural functional non-identifiability~\cite{loman2026structuralfunctionalidentifiabilitymodel}. We consider a minimal self-activation loop in which the species $X$ promotes its own production through a state-dependent reaction rate and undergoes linear decay. We model this using the following chemical reaction network
\begin{equation}
\ce{ \text{\o} ->C[$A(X;\vec{p})$] $X$, } 
\qquad
\ce{ $X$ ->C[$d$] \text{\o},}
\end{equation}
which yields the following ODE model
\begin{equation}
\frac{\mathrm{d}X}{\mathrm{d}t} = A(X;\vec{p}) - dX.
\end{equation}
For the ground-truth model, we take the activation function $A(X;\vec{p})$ to be a Hill function,
\begin{equation}
A(X;\vec{p})=\frac{vX^n}{X^n+K^n},
\end{equation}
with parameters $\vec{p}=(v,K,n)$. In the corresponding UDE model, this prescribed activation function is replaced by a universal function approximator $U(X;\vec{\theta})$.

We compare the identifiability properties of three models that differ only in the assumed prior knowledge of the activation function $A(X;\vec{p})$: (1) a fully specified model, in which both the functional form $A(X;\vec{p})=vX^n/(X^n+K^n)$ and the parameters $\vec{p}=(v,K,n)$ are assumed known, so that only the decay parameter $d$ is estimated; (2) a parameterised model, in which the Hill-function form is assumed known but the parameters $\vec{p}=(v,K,n)$ and $d$ are estimated; and (3) a UDE model, in which $A(X;\vec{p})$ is replaced by a universal function approximator $U(X;\vec{\theta})$, with both $\vec{\theta}$ and $d$ estimated. 

\begin{figure}[htbp]
    \centering
    \includegraphics[width=1.0\textwidth]{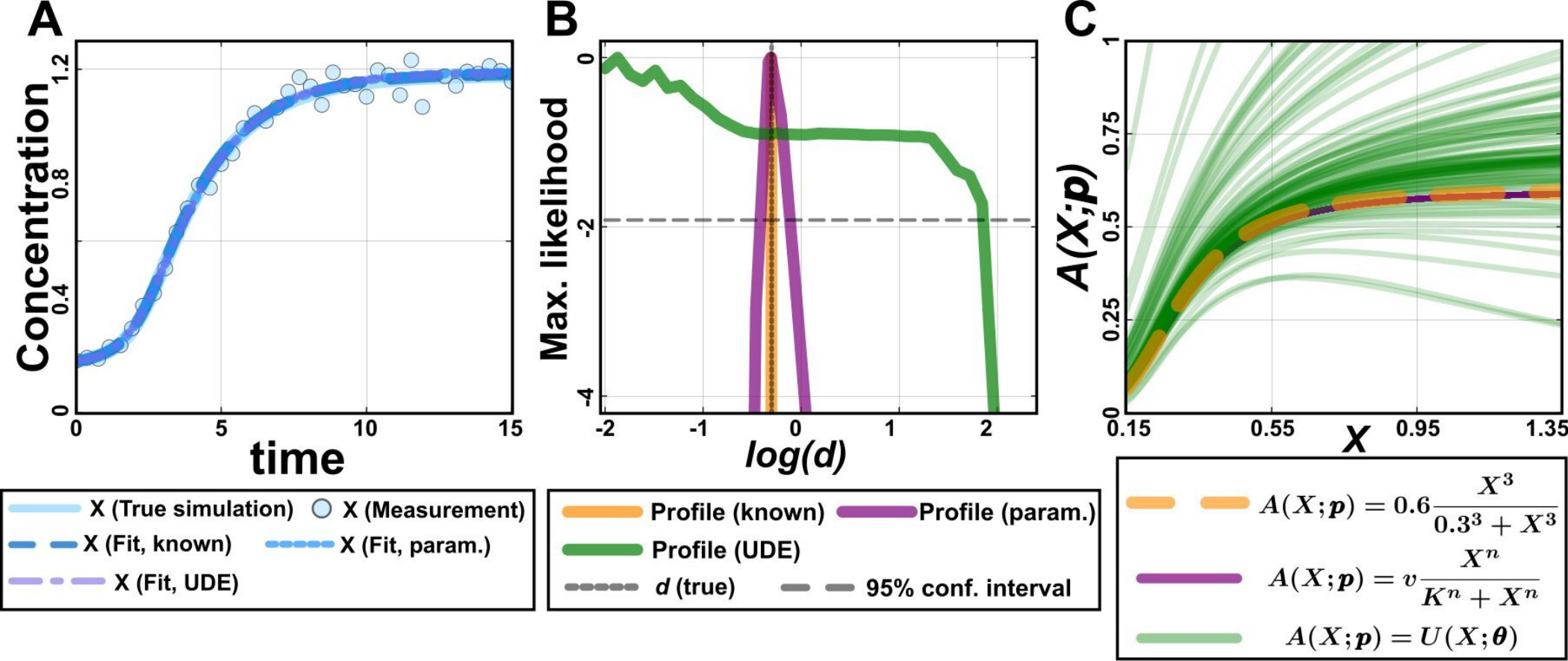}
    \caption{\textbf{Replacing a mechanistic reaction-rate function with a universal function approximator induces structural functional non-identifiability.} (A) Ground-truth simulation (solid line), synthetic data (points), and fitted trajectories for a fully specified model (dashed), a parameterised Hill-function model (dotted), and a UDE model (dash-dotted). All three models fit the observed data equally well. (B) Profile likelihoods for the decay parameter $d$. The fully specified and parameterised models remain identifiable, whereas the UDE model exhibits non-identifiability. (C) Activation functions recovered from repeated fits. The parameterised model consistently recovers the ground-truth Hill function, whereas the UDE admits multiple distinct functions that are all consistent with the observed data, demonstrating functional non-identifiability.}
    \label{figure:simple_sal_struc_nonident}
 \end{figure}

To illustrate the resulting identifiability properties numerically, we fit each of the three model formulations to the same synthetic dataset, where the data generation and fitting procedures are described in Section~\ref{section:detailed_methodology_synthetic_data} and Section~\ref{section:detailed_methodology_UDE_fitting}, respectively. We assess identifiability of the decay parameter $d$ using profile likelihood analysis, and identifiability of the activation function using ensemble plots of functional forms obtained from repeated independent fits, as described in Section~\ref{section:detailed_methodology_identifiability_analysis}. Figure~\ref{figure:simple_sal_struc_nonident}A shows the corresponding fitted trajectories, Figure~\ref{figure:simple_sal_struc_nonident}B the profile likelihoods for $d$, and Figure~\ref{figure:simple_sal_struc_nonident}C the recovered activation functions.


All three models fit the observed data equally well (Figure~\ref{figure:simple_sal_struc_nonident}A), but their identifiability properties differ markedly. The decay parameter $d$ is well constrained in the fully specified and parameterised models, whereas the UDE model admits a broad range of values of $d$ with essentially unchanged likelihood (Figure~\ref{figure:simple_sal_struc_nonident}B). Similarly, repeated fits of the parameterised model recover essentially the same activation function, while the UDE model admits a wide range of distinct functions that are all consistent with the observed data (Figure~\ref{figure:simple_sal_struc_nonident}C). Thus, although the three model formulations are indistinguishable in terms of their fitted trajectories, replacing the prescribed activation function with a universal function approximator introduces substantial ambiguity in both the decay parameter and the recovered activation function.

The recovered functions in Figure~\ref{figure:simple_sal_struc_nonident}C suggest a specific source of this ambiguity: changes in the decay rate can be compensated by adding a linear component to the activation function. This is not merely a consequence of finite or noisy data, but reflects a structural non-identifiability that can be demonstrated directly. Ignoring for the moment the non-negativity constraint on the reaction rate, consider any $d^*>0$, and define
\begin{equation}
\label{eq:ustar}
A^*(X)=A(X)+(d^*-d)X.
\end{equation}
Then
\begin{equation}
\frac{\mathrm{d}X}{\mathrm{d}t}
=A^*(X)-d^*X=A(X)-dX.
\end{equation}
Thus, for every admissible $d^*$, there exists a corresponding function $A^*$ that generates exactly the same vector field, and hence exactly the same trajectory as $(d,A)$ for any initial condition. Neither $d$ nor $A$ is therefore uniquely determined by ideal, noise-free observations of the system. Only their combination $A(X)-dX$, i.e.~the right-hand side of the differential equation, is identifiable. The activation function and decay parameter are therefore jointly structurally non-identifiable.

Physical constraints on the learnt function can reduce, or even eliminate, this structural non-identifiability. For example, suppose that $A^*$ is required to saturate, such that $A^{*\prime}(X)\to 0$ as $X\to\infty$. Since
\begin{equation}
A^{*\prime}(X)=A'(X)+(d^*-d),
\end{equation}
and $A'(X)\to0$ as $X\to\infty$, this constraint requires $d^*=d$, thereby uniquely recovering both $d$ and $A$. A weaker constraint is that $A^*$ is non-decreasing. Requiring
\begin{equation}
A^{*\prime}(X)=A'(X)+(d^*-d)\geq0,
\end{equation}
for all $X\geq0$, together with $A'(X)\to0$ as $X\to\infty$, implies $d^*\geq d$. Thus, monotonicity does not eliminate the structural non-identifiability, but restricts the family of indistinguishable solutions to $d^*\geq d$. Non-negativity of the learnt reaction rate provides a further restriction. From Equation~\eqref{eq:ustar},
\begin{equation}
A^*(X)\geq0
\quad\Longleftrightarrow\quad
A(X)+(d^*-d)X\geq0,
\end{equation}
and hence, for \(X>0\),
\begin{equation}
d^*\geq d-\frac{A(X)}{X}.
\end{equation}
If non-negativity is imposed globally, this condition must hold for all $X>0$. In practice, however, the fitted neural network is primarily constrained over the range of $X$ explored by the available data. The effective restriction on $d^*$ can therefore be substantially weaker when the observations cover only a limited state-space range. For the data used in Figure~\ref{figure:simple_sal_struc_nonident}B, for example, the non-negativity constraint does not impose the upper bound visible within the plotted profile range. The apparent upper bound on $d^*$ in Figure~\ref{figure:simple_sal_struc_nonident}B has a different origin: it is not implied by these structural constraints, but likely reflects the finite neural network architecture becoming unable to represent the increasingly steep, approximately linear functions required as $d^*$ increases.

More generally, this example illustrates that structural functional non-identifiability can arise even in very simple UDE settings, where the learnt function can absorb part of the mechanistic component. A detailed description of how to assess structural functional identifiability in UDEs, including for partially observed multi-variable models, is given in~\cite{loman2026structuralfunctionalidentifiabilitymodel}. In the remainder of this work, we focus on practical identifiability in UDEs, considering both mechanistic parameters and learnt functions.


\subsection{Practical identifiability depends on observation design and data quality in an extended self-activation loop}\label{section:results_extended_sa}

We next investigate practical functional identifiability in a setting where the structural non-identifiability encountered in Section~\ref{section:results_simple_sa} is removed. Specifically, we ask how recovery of an unknown function depends on the available observations and the quality of the data. We consider a minimal extension of the previous self-activation loop that is structurally identifiable. Here, we introduce a second species $Y$, such that $X$ promotes production of $Y$, which in turn promotes production of $X$. The production rate of $X$ is governed by the activation function $A$, which now depends on $Y$, while the production rate of $Y$ scales linearly with $X$. Both species decay at the same constant rate $d$. We model this using the following chemical reaction network
\begin{equation}
    \ce{ \text{\o} ->C[$A(Y;\vec{p})$] $X$},
    \quad
    \ce{ $X$ ->C[$d$] \text{\o}},
    \quad
    \ce{ \text{\o} ->C[$X$] $Y$},
    \quad
    \ce{ $Y$ ->C[$d$] \text{\o}},
\end{equation}
which yields the following ODE model
\begin{align}
\frac{\mathrm{d}X}{\mathrm{d}t} &= A(Y;\vec{p}) - dX, \\
\frac{\mathrm{d}Y}{\mathrm{d}t} &= X - dY.
\end{align}
For the ground-truth model, we take the activation function $A(Y;\vec{p})$ to be a Hill function
\begin{equation}
    A(Y;\vec{p})=vY^n/(Y^n + K^n),
\end{equation}
with parameters $\vec{p}=(v,K,n)$.

As in the previous section, we compare three levels of prior knowledge regarding the activation function $A(Y;\vec{p})$: (1) a fully specified model, in which both the functional form $A(Y;\vec{p})=vY^n/(Y^n+K^n)$ and the parameters $\vec{p}=(v,K,n)$ are assumed known and only the decay parameter $d$ is estimated; (2) a parameterised model, in which the Hill-function form is assumed known but the parameters $\vec{p}=(v,K,n)$ and $d$ are estimated; and (3) a UDE model, in which $A(Y;\vec{p})$ is replaced by a universal function approximator $U(Y;\vec{\theta})$, with both $\vec{\theta}$ and $d$ estimated. Synthetic data are generated as described in Section~\ref{section:detailed_methodology_synthetic_data} under three observation scenarios: (1) measurements of $X$ only with low observation noise; (2) measurements of both $X$ and $Y$ with low observation noise; and (3) measurements of both $X$ and $Y$ with high observation noise. For each scenario, we use the methods described in Section~\ref{section:detailed_methodology_UDE_fitting} and Section~\ref{section:detailed_methodology_identifiability_analysis} to fit all three model formulations, and to assess practical identifiability of $d$ using profile likelihood analysis and of $A(Y;\vec{p})$ using ensembles of the fitted function $U(Y;\vec{\theta})$.


\begin{figure}[htbp]
    \centering
    \includegraphics[width=\textwidth]{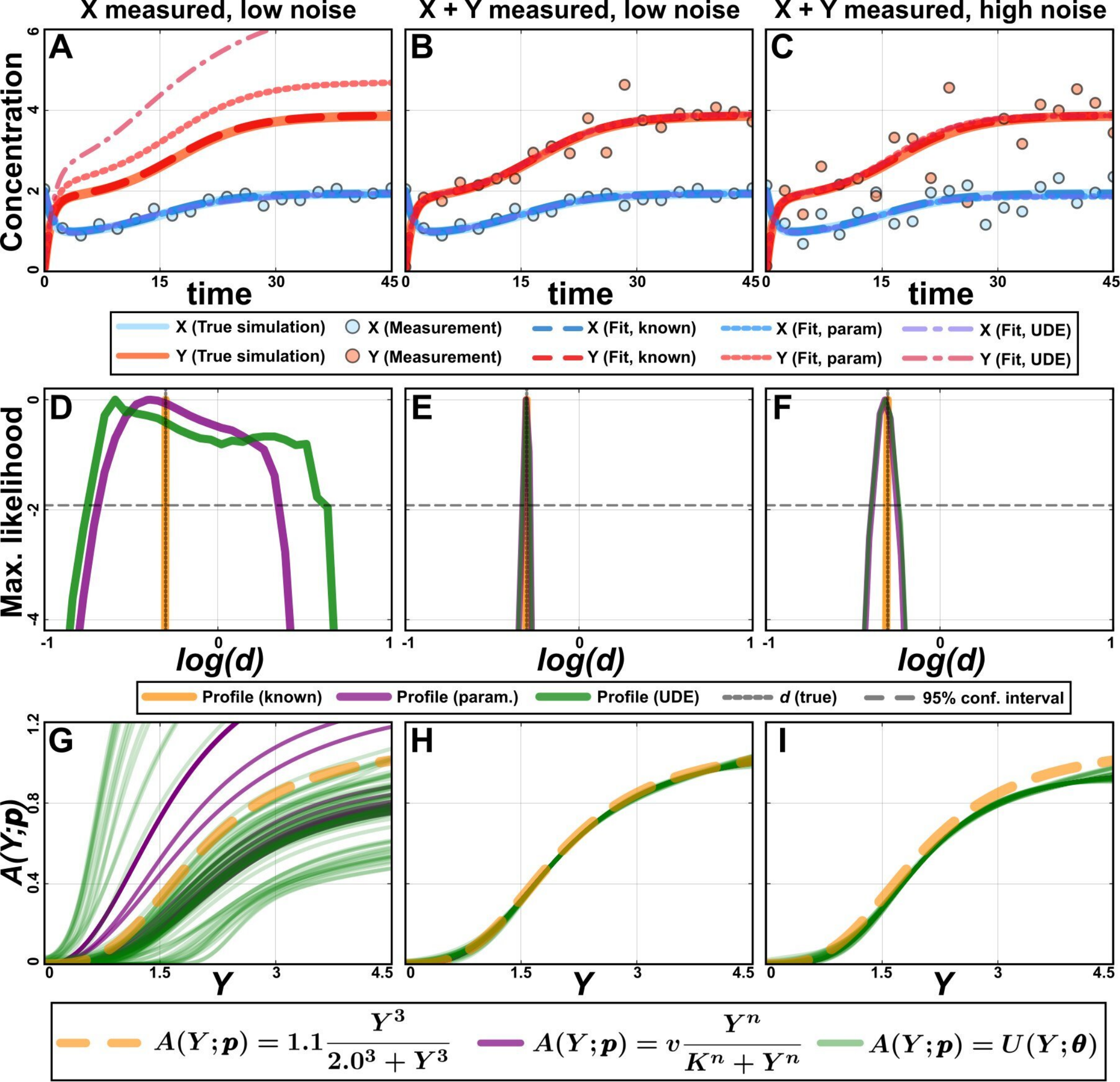}
    \caption{\textbf{Practical functional identifiability depends strongly on the available observations.} Three model formulations are fitted under three observation scenarios: (1) measurements of $X$ only at low noise (A,D,G); (2) measurements of both $X$ and $Y$ at low noise (B,E,H); and (3) measurements of both species at high noise (C,F,I). (A--C) All models reproduce the observed data, but when only $X$ is observed, the cases where the activation function is not fully known produce erroneous predictions for $Y$. (D--F) Corresponding profile likelihoods for the decay parameter $d$. Practical parameter non-identifiability occurs in the same cases as when predictions of $Y$ are incorrect. (G--I) Ensembles of recovered activation functions. Functional non-identifiability occurs under the same conditions as parameter non-identifiability.}
    \label{figure:extended_sal_pract_ident}
\end{figure}


All three models reproduce the observed data under each observation scenario (Figure~\ref{figure:extended_sal_pract_ident}A--C), but their ability to recover the unobserved dynamics (in this case $Y$) differs substantially. When only $X$ is measured, both the parameterised and UDE models produce inaccurate trajectories for the unobserved species $Y$ (Figure~\ref{figure:extended_sal_pract_ident}A), with substantial variability across independently fitted models (Supplementary Figure~\ref{sfigure:extended_sal_pred_uncertain}). Next, when both species are observed, all three formulations recover the underlying trajectories accurately, even at the higher observation-noise level considered here (Figure~\ref{figure:extended_sal_pract_ident}B,C). The same pattern is apparent in both parameter and functional identifiability. In the $X$-only setting, the parameterised and UDE models exhibit broad profile likelihoods for the decay parameter $d$ and substantial variability among the recovered activation functions (Figure~\ref{figure:extended_sal_pract_ident}D,G). When both species are observed, the profiles for $d$ become substantially narrower and the recovered activation functions converge to a common form close to the ground truth (Figure~\ref{figure:extended_sal_pract_ident}E,F,H,I). Thus, for this model and these observation regimes, poor recovery of the unobserved dynamics coincides with poor practical parameter and functional identifiability. We highlight that the non-identifiability observed when only $X$ is measured is practical rather than structural. With more densely sampled and substantially lower-noise observations of $X$, both the mechanistic parameters and the activation function can be recovered accurately (Supplementary Figure~\ref{sfigure:extended_sal_pract_ident_maxdata}). 

Across the observation scenarios considered above, replacing the parameterised activation function $A(Y;\vec{p})$ with a universal function approximator produces only a small additional loss of identifiability. In this model, therefore, knowledge of the precise Hill-function form contributes relatively little additional information over the tested observation regimes. This behaviour is reproduced across independent datasets generated from the same ground-truth model (Supplementary Figures~\ref{sfigure:extended_sal_pract_ident_rep1} and \ref{sfigure:extended_sal_pract_ident_rep2}), an alternative parameter set (Supplementary Figure~\ref{sfigure:extended_sal_pract_ident_alt_ps}), and a different model (the Goodwin oscillator; Supplementary Figure~\ref{sfigure:gwo_pract_ident}). It is not universal, however: under an alternative parameter regime and initial condition, the parameterised Hill-function model achieves a measurable identifiability advantage over the UDE formulation (Supplementary Figure~\ref{sfigure:extended_sal_pract_ident_lowconc_ps}).

To examine the dependence on data quality more systematically, we next consider a composite data-quality axis in which sampling density is increased while observation noise is simultaneously reduced. Using the extended self-activation loop, we generate datasets ranging from sparsely sampled, high-noise observations to densely sampled, low-noise observations. Across this range, we quantify parametric identifiability using the width of confidence intervals derived from profile likelihood analysis, functional identifiability using the mean $L^2$ distance between fitted functions, and predictive performance using the prediction error for the unobserved species. Across the regimes considered, the parameterised and UDE models exhibit similar behaviour (Supplementary Figures~\ref{sfigure:extended_sal_datavar},\ref{sfigure:extended_sal_sample_fits}), with parametric identifiability, functional identifiability, and prediction of the unobserved dynamics all improving as the data become more informative.


\subsection{An insect population model demonstrates a complex relationship between parametric identifiability, functional identifiability, and predictive power}\label{section:results_mosquitoes}

The extended self-activation loop in Section~\ref{section:results_extended_sa} exhibits a close correspondence between (1) practical parameter identifiability, (2) practical functional identifiability, and (3) ability to predict unobserved dynamics (see Supplementary Figure~\ref{sfigure:extended_sal_datavar}). We next ask whether the correlation between these properties can be generalised to other models. To investigate this, we consider a multi-stage population dynamics model of an insect population (e.g.~mosquitoes)~\cite{kaye_impact_2024}. The population progresses through three life stages, egg ($E$), larvae ($L$), and mature individuals ($M$), with linear transition rates between stages and linear mortality rates for eggs and mature individuals. To represent density-dependent competition, larval mortality depends on population size through the function $d_L(L;\vec{p})$. We model this using the following chemical reaction network
\begin{equation}
\ce{ $E$ ->C[$k_{E,L}$] $L$},
\quad
\ce{ $E$ ->C[$d_E$] \text{\o}},
\quad
\ce{ $L$ ->C[$k_{L,M}$] $M$},
\quad
\ce{ $L$ ->C[$d_L(L;\vec{p})$] \text{\o}},
\quad
\ce{ $M$ ->C[$r$] $M + E$},
\quad
\ce{ $M$ ->C[$d_M$] \text{\o}},
\end{equation}
which yields the following ODE model
\begin{align}
\frac{\mathrm{d}E}{\mathrm{d}t} &= r M - (k_{E,L} + d_E) E,\\
\frac{\mathrm{d}L}{\mathrm{d}t} &= k_{E,L} E - (k_{L,M} + d_L(L;\vec{p}))L,\\
\frac{\mathrm{d}M}{\mathrm{d}t} &= k_{L,M} L - d_M M.
\end{align}
We model competition using the ground-truth quadratic mortality function
\begin{equation}
d_L(L;\vec{p})=d_{L,1}+d_{L,2}L+d_{L,3}L^2,    
\end{equation}
with parameters $\vec{p}=(d_{L,1},d_{L,2},d_{L,3})$ where $d_{L,i}>0$ for $i=1,2,3$.

We note that the larval equation alone contains a potential confounding interaction between the transition rate $k_{L,M}$ and the mortality function $d_L(L;\vec{p})$, since they enter the equation through the sum $k_{L,M}+d_L(L;\vec{p})$. Unlike in Section~\ref{section:results_simple_sa}, however, $k_{L,M}$ also appears independently in the equation for $M$, providing additional dynamical information with which to distinguish these contributions. Despite this, the full model with all mechanistic parameters unknown is structurally non-identifiable, and we therefore restrict the subsequent analysis to a structurally identifiable subset of unknown quantities. In particular, $k_{E,L}$, $d_E$, $k_{L,M}$, and $d_M$ are assumed known throughout this section. 

As in the previous examples, we compare three levels of prior knowledge: (1) a fully specified model, in which the functional form and parameters of $d_L(L;\vec{p})$ are assumed known and only the mechanistic parameter $r$ is estimated; (2) a parameterised model, in which the functional form is assumed known but $\vec{p}$ and $r$ are estimated; and (3) a UDE model, in which $d_L(L;\vec{p})$ is replaced by a universal function approximator $U(L;\vec{\theta})$, with both $\vec{\theta}$ and $r$ estimated. Using this model, we repeat the data-quality analysis performed for the extended self-activation loop, generating independent synthetic datasets while simultaneously increasing sampling density and decreasing observation noise. We fit each of the three model formulations to measurements of $M$ and $L$, and assess predictive performance using the resulting trajectories of $E$. Parametric identifiability is quantified using the width of the profile-likelihood confidence interval for $r$ and functional identifiability using the mean $L^2$ distance between fitted mortality functions (Supplementary Section~\ref{ssection:detailed_methodology_datavar_scans}).

\begin{figure}[htbp]
    \centering
    \includegraphics[width=\textwidth]{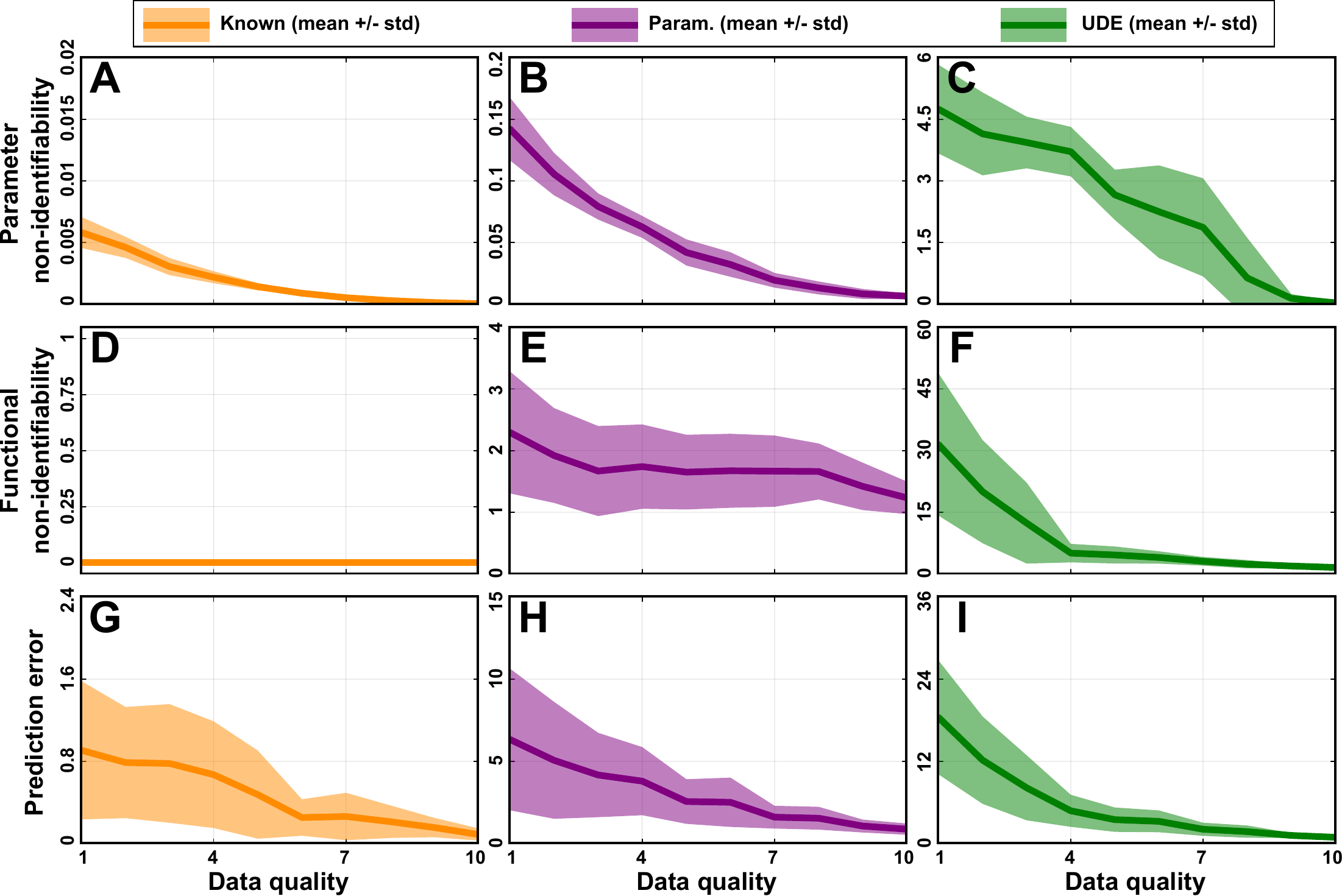}
    \caption{\textbf{Parametric and functional identifiability do not need to be correlated.} Three versions of the insect population model are compared: (A,D,G) the mortality function is fully specified; (B,E,H) the functional form is known but its parameters are estimated; and (C,F,I) the mortality function is represented by a universal function approximator. Parametric identifiability (A--C), functional identifiability (D--F), and predictive performance (G--I) are shown as functions of data quality, which increases from sparsely sampled, high-noise observations to densely sampled, low-noise observations (Supplementary Section~\ref{ssection:detailed_methodology_datavar_scans}). Relative to the mechanistic models, the UDE formulation exhibits large reductions in both predictive performance and parametric and functional identifiability. Furthermore, unlike for the extended self-activation loop (Supplementary Figure~\ref{sfigure:extended_sal_datavar}), the different measures are not necessarily correlated e.g. for the UDE, parameter identifiability improves more slowly with improvements to data, as compared to functional identifiability and predictive performance. Example fits, profile likelihoods, and recovered functions are shown in Supplementary Figure~\ref{sfigure:insects_example_fits}.}
    \label{figure:mosquitoes_data_quality_scan}
\end{figure}

The resulting behaviour differs substantially from that of the extended self-activation loop. For both the parameterised and UDE models, the correlation between improvements in parametric and functional identifiability, as data quality improves, is lost (Figure~\ref{figure:mosquitoes_data_quality_scan}). Furthermore, all three measures (predictive performance, parametric identifiability, and functional identifiability) deteriorate drastically as we move from the specified, to the parameterised, to the UDE model. This is in contrast to the extended self-activation loop, where the parameterised and UDE models exhibited similar predictive performance and degree of identifiability across all datasets. These results show that practical parameter identifiability, practical functional identifiability, and predictive performance do not need to be correlated. They also demonstrate that the identifiability cost of replacing a prescribed functional form with a universal function approximator is model-dependent.


\subsection{An incoherent feedforward loop can exhibit functional identifiability for multiple learnt functions simultaneously} \label{section:results_ic_feedforward}

Thus far, we have considered models containing a single unknown function. We now ask whether representing multiple functions simultaneously using universal function approximators necessarily leads to a larger loss of identifiability. To investigate this, we consider an incoherent feedforward loop, a well-studied systems biology motif that generates behaviours such as pulse responses~\cite{goentoro_incoherent_2009}. In this system, species $X$ represses $Y$ and activates $Z$, while $Y$ activates $Z$. We model this using the following chemical reaction network
\begin{equation}
\ce{ \text{\o} ->C[$1$] $X$,}
\quad
\ce{ \text{\o} ->C[$R(X;\vec{p}_R)$] $Y$,}
\quad
\ce{ $Y$ ->C[$A(X;\vec{p}_A)$] $Z$,}
\quad
\ce{ $X$ ->C[$d$] \text{\o},}
\quad
\ce{ $Y$ ->C[$d$] \text{\o},}
\quad
\ce{ $Z$ ->C[$d$] \text{\o},}
\end{equation}
which yields the following ODE model
\begin{align}
\frac{\mathrm{d}X}{\mathrm{d}t} &= 1 - d X, \\
\frac{\mathrm{d}Y}{\mathrm{d}t} &= R(X;\vec{p}_R) - A(X;\vec{p}_A) Y - d Y,\\
\frac{\mathrm{d}Z}{\mathrm{d}t} &= A(X;\vec{p}_A) Y - d Z.
\end{align}
We use Hill-type functional forms for the ground-truth activation and repression functions:
\begin{equation}
A(X;\vec{p}_A)=v_A\frac{X^{n}}{K_A^{n} + X^{n}},
\qquad
R(X;\vec{p}_R)=v_R\frac{K_R^{n}}{K_R^{n} + X^{n}}.
\end{equation}

We compare five levels of prior knowledge regarding the activation and repression functions: (1) both functions and their parameters are known; (2) both functional forms are known but their parameters are estimated; (3) the repression function $R$ is parameterised while the activation function $A$ is represented by a universal function approximator; (4) the activation function $A$ is parameterised while the repression function $R$ is represented by a universal function approximator; and (5) both functions are represented by universal function approximators. In all cases, the mechanistic decay parameter $d$ is estimated. Synthetic data are generated with all three species observed\footnote{Note that successful fits to data can be achieved while observing a single species only; however, this requires measurements from multiple experimental conditions---see Supplementary Figure~\ref{sfigure:ic_feedforward_measurement_type_identifiability}.} (Section~\ref{section:detailed_methodology_synthetic_data}). For each model formulation, we assess practical identifiability of $d$ using profile likelihood analysis and practical identifiability of $A$ and $R$ using ensemble plots (Sections~\ref{section:detailed_methodology_UDE_fitting} and \ref{section:detailed_methodology_identifiability_analysis}).
 

\begin{figure}[htbp]
    \centering
    \includegraphics[width=\textwidth]{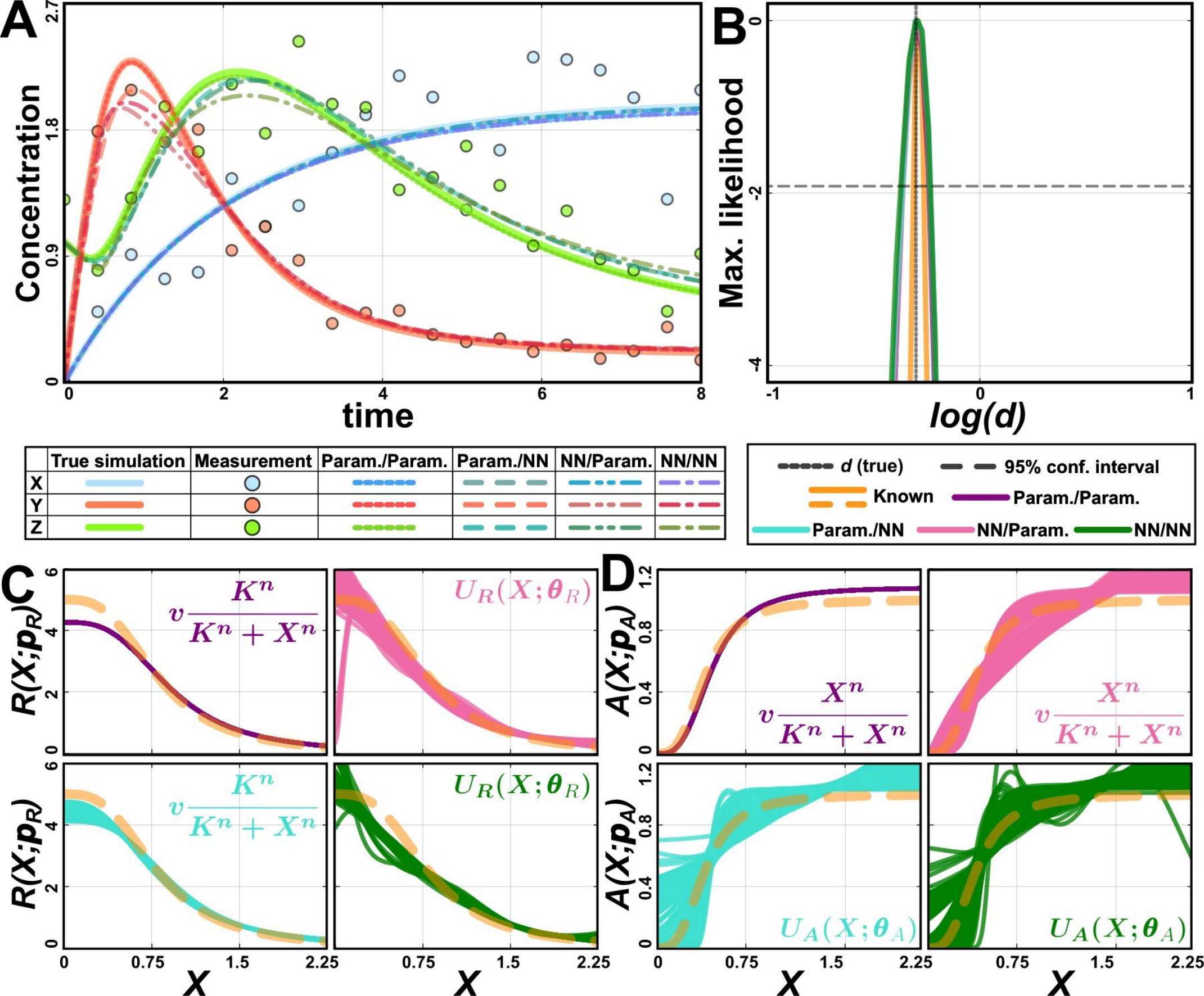}
    \caption{\textbf{Multiple unknown functions can be recovered simultaneously with only a modest additional loss of identifiability in the incoherent feedforward loop.} The incoherent feedforward loop model is fitted under five levels of prior knowledge, ranging from fully specified activation and repression functions to both functions being represented by universal function approximators. (A) All model variants accurately reproduce the observed dynamics. (B) Profile likelihoods for the mechanistic decay parameter $d$, which remains well constrained across all configurations, including when both functions are represented by universal function approximators. (C,D) Ensembles of recovered repression and activation functions, respectively. Replacing a prescribed functional form with a universal function approximator does reduce that function's identifiability. However, doing so has a relatively small impact on the other function's identifiability (whenever it is parameterised or fully unknown). Colours denote the different model formulations described in the text.}
    \label{figure:ic_feedforward_2nn}
\end{figure}


All five model variants accurately reproduce the observed dynamics (Figure~\ref{figure:ic_feedforward_2nn}A), despite substantial differences in how the activation and repression functions are represented. The practical identifiability of the decay parameter $d$ is also largely unaffected by replacing either or both Hill functions with universal function approximators: the likelihood profile remains well constrained across all five formulations (Figure~\ref{figure:ic_feedforward_2nn}B). Next, we consider how generalising a function from a Hill function to a universal function approximator affects functional identifiability---both for the function in question, but also for the other system function. While, across all cases, the ground-truth function is recovered, there is a noticeable reduction in identifiability for functions that are generalised to universal function approximators. However, in all cases, this generalisation has relatively little impact on the other function's identifiability. These results are repeated (using the same model and simulation conditions) in Supplementary Figure~\ref{sfigure:ic_feedforward_2nn_repeat}, demonstrating that they are robust to the realisation of the noise.

As in the previous example, these results are unlikely to generalise across all model structures and observation schemes. Nevertheless, they demonstrate that successively replacing known mechanistic terms with universal function approximators does not necessarily result in a cumulative loss of identifiability. In particular, introducing additional unknown functions need not reduce the identifiability of the remaining model components, whether parameters or other functions.


\subsection{Architectural constraints can substantially improve practical identifiability in a linear pathway model}\label{section:results_linear_pathway}

Thus far, the universal function approximators considered in this work have been constrained only to produce non-negative outputs, reflecting the requirement that reaction rates are non-negative. In many applications, however, additional qualitative or quantitative information about an unknown function may be available. We therefore ask whether incorporating such prior knowledge directly into the neural network architecture can improve practical identifiability. Specifically, we consider two classes of constraints: monotonicity, which restricts the function to be strictly increasing or decreasing, and boundedness, which restricts the function to a prescribed range $U_{\min}<U(\vec{X};\vec{\theta})<U_{\max}$. Both constraints are enforced directly through the neural network architecture, as described in Supplementary Section~\ref{ssection:detailed_methodology_nn_constraints}.

To investigate their effect, we consider a linear pathway model consisting of four species, in which each species represses production of the next species in the pathway. We model this using the following chemical reaction network
\begin{equation}
\ce{$X_1$ ->C[$d$] \text{\o}},
\,\,
\ce{\text{\o} ->C[$R(X_1;\vec{p})$] $X_2$},
\,\,
\ce{$X_2$ ->C[$d$] \text{\o}},
\,\,
\ce{\text{\o} ->C[$R(X_2;\vec{p})$] $X_3$},
\,\,
\ce{$X_3$ ->C[$d$] \text{\o}},
\,\,
\ce{\text{\o} ->C[$R(X_3;\vec{p})$] $X_4$},
\,\,
\ce{$X_4$ ->C[$d$] \text{\o},}
\end{equation}
which yields the following ODE model
\begin{align}
\frac{\mathrm{d}X_1}{\mathrm{d}t} &= - d X_1, \\
\frac{\mathrm{d}X_2}{\mathrm{d}t} &= R(X_1;\vec{p}) - d X_2, \\
\frac{\mathrm{d}X_3}{\mathrm{d}t} &= R(X_2;\vec{p}) - d X_3, \\
\frac{\mathrm{d}X_4}{\mathrm{d}t} &= R(X_3;\vec{p}) - d X_4.
\end{align}
We use the following common ground-truth repression function $R$ across the pathway, 
\begin{equation}
    R(X;\vec{p}) = v_0 + vK^n/(X^n + K^n),
\end{equation} 
with $\vec{p}=(v_0,v,K,n)$. 

We compare four levels of prior knowledge regarding the repression function: (1) a parameterised model in which the functional form $R(X;\vec{p})$ is known but its parameters are estimated; (2) a UDE model in which $R$ is constrained to be monotonically decreasing and bounded (assuming the true bounds $(v_0, v_0 + v)$ are known and used); (3) a UDE model constrained only to be monotonically decreasing; and (4) an unconstrained UDE model (aside from the non-negativity constraint imposed throughout this work). In all cases, the mechanistic decay parameter $d>0$ is estimated. Synthetic data are generated with observations of $X_4$ only (Section~\ref{section:detailed_methodology_synthetic_data}). For each formulation, we assess practical identifiability of $d$ using profile likelihood analysis and functional identifiability using ensemble plots.

\begin{figure}[htbp]
    \centering
    \includegraphics[width=\textwidth]{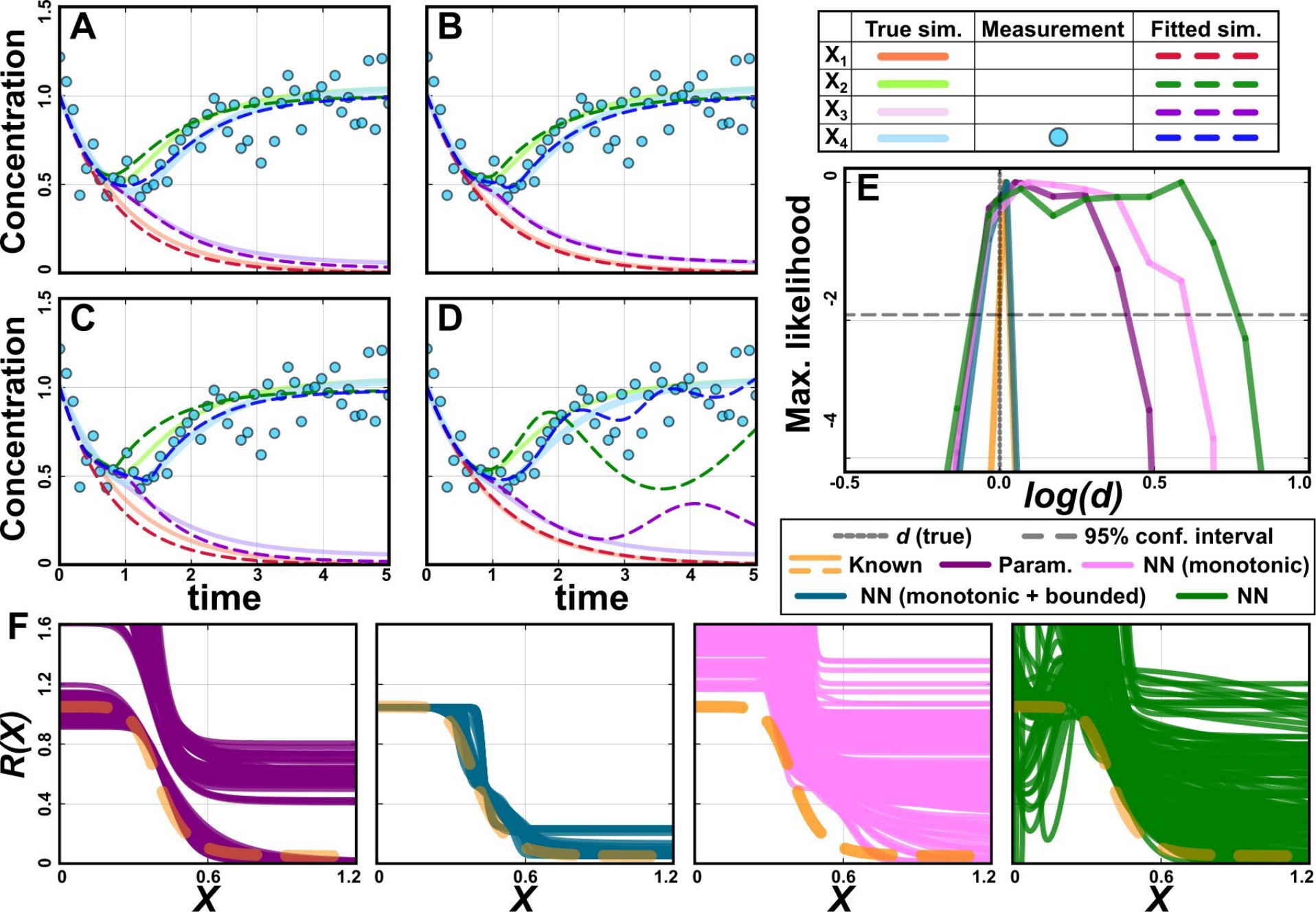}
    \caption{\textbf{Architectural constraints substantially improve practical identifiability in the linear pathway model.} Four models are compared: (A) a parameterised repression function; (B) a UDE constrained to be monotonically decreasing and bounded; (C) a UDE constrained only by monotonicity; and (D) an unconstrained UDE. (A--D) Fitted trajectories. The parameterised model and the monotonic and bounded UDE accurately recover the unobserved dynamics, whereas progressively weaker constraints lead to poorer recovery. (E) Profile likelihoods for the decay parameter $d$. Practical identifiability of $d$ deteriorates as constraints on the unknown function are removed. (F) Ensembles of recovered repression functions. Functional variability follows the same ordering, with the monotonic and bounded UDE recovering much of the functional constraint obtained by prescribing the fully known Hill-function form.}
    \label{figure:ln_pathway_nn_constraints}
\end{figure}

The effect of incorporating constraints is substantial (Figure~\ref{figure:ln_pathway_nn_constraints}). The parameterised model and the UDE constrained by both monotonicity and boundedness recover the unobserved dynamics accurately, whereas the unconstrained UDE produces substantial errors (Figure~\ref{figure:ln_pathway_nn_constraints}A--D). Imposing monotonicity alone yields intermediate performance with respect to predicting unobserved dynamics. Similar ordering is apparent in both parameter and functional identifiability, however, here the monotonic and bounded UDE also outperforms the parameterised Hill model. The profile likelihood for $d$ remains relatively narrow for the monotonic and bounded UDE, approaching that where the function is known exactly. The parameterised Hill, the monotonic, and the unconstrained UDE models exhibit progressively lower identifiability (Figure~\ref{figure:ln_pathway_nn_constraints}E). Similarly, the monotonic and bounded UDE yields a tightly clustered ensemble of recovered repression functions, while the other models exhibit substantial variability in recovered functions (Figure~\ref{figure:ln_pathway_nn_constraints}F). This analysis is repeated across repeated synthetic datasets for the same linear pathway model (Supplementary Figures~\ref{sfigure:ln_pathway_nn_constraints_repeat2},\ref{sfigure:ln_pathway_nn_constraints_repeat3}) and using the Goodwin oscillator model (Supplementary Figures~\ref{sfigure:gwo_nn_constraints_repeat1},\ref{sfigure:gwo_nn_constraints_repeat2}) with similar results.


Again, the impact of including constraints on the learnt function likely varies with model and data. Furthermore, bounds as tight as those used here may not always be available in practice, and similar bounds could also be imposed on the parameterised Hill function model. However, the dramatic improvement in identifiability for the monotonic and bounded UDE is still remarkable. Indeed, this combination of constraints carries almost identical information to knowing the exact functional form (i.e.~both the Hill function form, and its parameter values, are known). These results show that relatively basic quantitative knowledge about the system can translate into dramatic gains in model performance. More generally, they support treating UDEs not as an all-or-nothing choice between a known and an unknown component, but as a granular framework for encoding exactly what is, and is not, known about a system.


\subsection{A UDE can mitigate functional-form misspecification while preserving practical parameter identifiability}\label{section:results_modified_sir}

An advantage of UDEs is that they can model complex functional forms that do not follow exact algebraic expressions. While such algebraic expressions are frequently used in contemporary literature, they may generate errors due to model misspecification. If the chosen expression is a sufficiently good approximation, the errors should be small, however, UDEs offer a more robust approach~\cite{browning_cautionary_2025}. Here, we investigate the trade-off between model misspecification error and the loss of identifiability due to a UDE generalisation.

For this example, we consider a modified susceptible--infected--recovered (SIR) infectious disease model~\cite{brauer_compartmental_2008,avram_advancing_2024}. The classical SIR model is typically expressed through the reactions
\begin{equation}
\ce{ $S + I$ ->C[$\beta$] $2I$},
\qquad
\ce{ $I$ ->C[$\gamma$] $R$},
\end{equation}
which can be equivalently written as
\begin{equation}
\ce{ $S$ ->C[$\beta I$] $I$},
\qquad
\ce{ $I$ ->C[$\gamma$] $R$},
\end{equation}
where the infection rate depends linearly on $I$. This linear dependence is a modelling assumption and may be inappropriate when the effective infection process depends nonlinearly on the population state. In this case, a more appropriate description could be written as
\begin{equation}
\ce{ $S$ ->C[$B(S,I;\vec{p})$] $I$},
\qquad
\ce{ $I$ ->C[$\gamma$] $R$}.
\end{equation}
This corresponds to the following ODE model
\begin{align}
\frac{\mathrm{d}S}{\mathrm{d}t} &= - B(S,I;\vec{p})S,\\
\frac{\mathrm{d}I}{\mathrm{d}t} &= B(S,I;\vec{p})S - \gamma I,\\
\frac{\mathrm{d}R}{\mathrm{d}t} &= \gamma I,
\end{align}
for which we will use the following ground-truth infection function
\begin{equation*}
    B(S,I;\vec{p}) = \beta S^{\alpha_1} I^{\alpha_2},
\end{equation*}
with $\vec{p} = (\beta, \alpha_1,\alpha_2)$~\cite{taghvaei_fractional_2020}.

We compare four models with different assumptions about the infection function: (1) a fully specified model, in which both $B(S,I;\vec{p})$ and its parameters are known; (2) a parameterised model, in which the functional form of $B(S,I;\vec{p})$ is known but its parameters are estimated; (3) a UDE model, in which $B(S,I;\vec{p})$ is replaced by a universal function approximator; and (4) a misspecified model that assumes the standard mass-action infection rate $B(S,I)=\beta I$. In all cases, the recovery rate $\gamma$ is estimated. Synthetic data are generated for $I$ only (Section~\ref{section:detailed_methodology_synthetic_data}).


\begin{figure}[htbp]
    \centering
    \includegraphics[width=\textwidth]{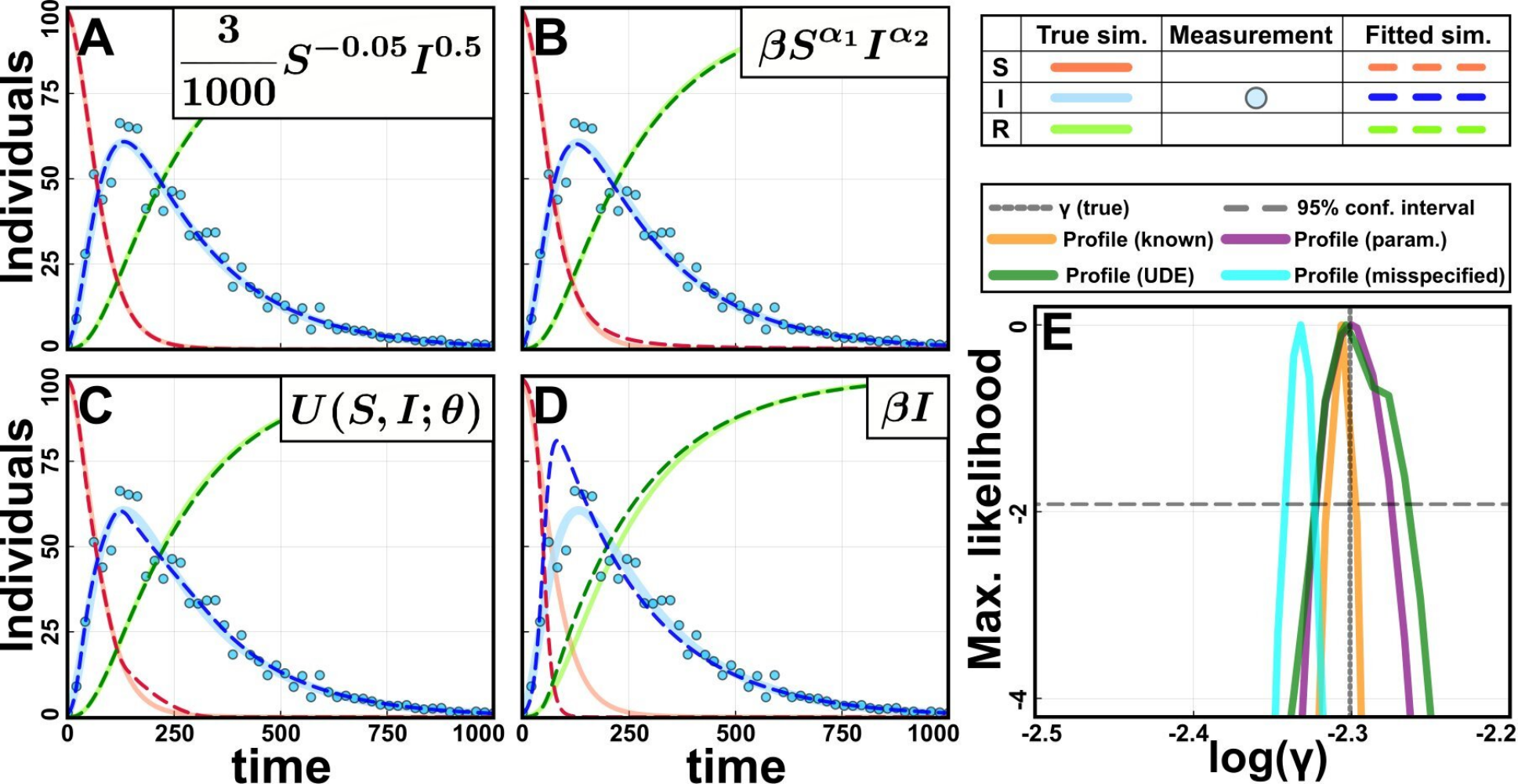}
    \caption{\textbf{A flexible infection-rate representation mitigates functional-form misspecification with little loss of practical parameter identifiability.} Four model formulations are fitted to synthetic data generated from the modified SIR model: (A) fully specified; (B) parameterised; (C) UDE; and (D) misspecified. (A--D) All four models reproduce the observed $I$ data, but the misspecified model produces substantially different predictions for the unobserved $S$ and $R$ dynamics with a larger outbreak than the ground truth. (E) Profile likelihoods for the recovery parameter $\gamma$. The misspecified model yields a narrow profile centred away from the true value, whereas the fully specified, parameterised, and UDE models yield profiles consistent with the ground-truth value. Note the smaller than usual $x$-axis, indicating that the misspecification error is relatively small. The UDE incurs only a modest loss of practical identifiability of $\gamma$ relative to the parameterised model.}
    \label{figure:modified_sir_misspecification}
\end{figure} 


All four models reproduce the observed $I$ data closely (Figure~\ref{figure:modified_sir_misspecification}A--D), indicating that, in this example, agreement with the observed data alone is insufficient to reveal the incorrect functional assumption. Notably, however, the trajectories of the misspecified model do not align with the ground-truth simulations, and it predicts a larger overall outbreak (Figure~\ref{figure:modified_sir_misspecification}D). The misspecification also affects inference of the mechanistic recovery parameter $\gamma$. The misspecified model yields a narrow profile likelihood, but its optimum is displaced from the ground-truth value, demonstrating that a parameter can appear practically well identified while nevertheless being biased by an incorrect functional assumption (Figure~\ref{figure:modified_sir_misspecification}E). In contrast, the fully specified, parameterised, and UDE models yield profile likelihoods consistent with the ground-truth value of $\gamma$. The UDE profile likelihood is negligibly broader than that of the parameterised model, despite imposing substantially weaker assumptions on the functional form of the infection process. These results are repeated (using the same model and simulation conditions) in Supplementary Figure~\ref{sfigure:modified_sir_misspecification_repeat}, demonstrating that they are robust to the realisation of the noise. Functional identifiability analysis is omitted in this analysis as the discovered function takes two inputs, which prevents us from displaying easily legible ensemble plots.

The consequences of misspecification are relatively modest in this example. Stronger misspecification errors occur in the incoherent feedforward loop model (Supplementary Figure~\ref{sfigure:ic_feedforward_misspecification}), although the misspecification is also more apparent there. Together, these examples illustrate an important benefit of functional flexibility: when the prescribed functional form is uncertain, a UDE can reduce sensitivity to this source of model misspecification while retaining useful predictive properties. This does not imply that UDEs eliminate model misspecification more generally, since the surrounding mechanistic structure remains assumed. Rather, the result illustrates why some loss of identifiability may be an acceptable cost of relaxing uncertain functional assumptions: a tightly constrained but misspecified model can yield precise yet biased mechanistic inference, something a UDE formulation might avoid.


\section{Discussion}\label{section:discussion}

Since their introduction, UDEs have been applied across a wide range of fields, leveraging their ability to account for unknown system dynamics. This flexibility, however, brings an increased propensity for overfitting and spurious discovery. Until now, we have lacked systematic tools for detecting this. In this work, we have shown that identifiability analysis can be extended to UDEs by decomposing the problem into parametric and functional components. Just as parametric identifiability admits a well-established distinction between structural and practical variants, we have shown that the same distinction applies to functional identifiability. We also showed that practical identifiability of mechanistic parameters can be assessed using established methods such as profile likelihood analysis and functional identifiability can be assessed using an ensemble plot-based approach.
In this work, we apply this framework to models spanning systems biology, epidemiology, and population dynamics---including settings with multiple simultaneously fitted functions, shared functions across equations, and multivariate inputs. First, we showed how structural functional non-identifiability can be reduced, or even removed, by a combination of physical constraints and data. Next, through multiple examples, we analysed how replacing a mechanistically specified function with a universal function approximator affects identifiability. We noted how, for some, but not all, cases, there exists a strong correlation between parametric identifiability, functional identifiability, and model predictive power. In another example, we showed a simple model for which multiple unknown functions could be identified simultaneously, and where the identifiability of one function was mostly unaffected by our level of knowledge of the other function. We also showed how functional constraints, such as monotonicity and bounds, could dramatically affect a model's identifiability properties. In our example, relatively simple constraints negated almost the full identifiability penalty for introducing the UDE. Finally, we considered the trade-off where the introduction of a UDE does reduce identifiability, but also reduces the effect of model misspecification. Specifically, in situations where the UDE-associated identifiability penalty is low, UDEs offer a potentially more robust modelling approach.

Throughout many of these examples, the identifiability cost of generalising a mechanistically-defined function to a universal function approximator is strikingly low. In some cases, identifiability is maintained even under substantial noise. While this might not always bring a UDE's identifiability on par with a model where the function is fully known (as we observed in this paper), it underscores how important the incorporation of such constraints can be. Finally, we also identified cases where UDEs introduce a noticeable identifiability penalty. This highlights how the generalisation to a UDE can have negligible to severe impact on a model. Our goal here is to present tools for assessing when this is the case, allowing practitioners to determine when the UDE approach is tractable for their applications. 
 
We specifically apply UDEs to chemical reaction network models. This approach carries three main advantages. First, just as rules exist for deriving ODEs from chemical reaction networks, analogous rules generate other model types, providing a direct route to non-ordinary UDEs (such as universal stochastic differential equations). Second, the non-constant rate functions that arise in these models---used to represent complex or poorly characterised mechanisms---are precisely the kind of unknown functional components that UDEs are designed to learn. Identifying reactions with such rates therefore directly identifies natural use cases for UDEs. Third, because these rate functions are associated with specific physical processes (e.g.~a protein's production dynamics or inter-species competition), the learnt functions are immediately interpretable: they can be plotted and related to the underlying biology without further post-processing.

This mode of interpretation distinguishes our function-centric approach from an alternative strategy where symbolic regression is applied to the fitted neural network to recover explicit algebraic expressions~\cite{rackauckas_universal_2021,cranmer_interpretable_2023}. The symbolic regression approach has clear strengths: it can reveal new mechanistic rules and, by projecting onto a restricted function space, may promote identifiability. However, it also introduces additional implementation choices that make the methodology more involved and may affect results: selecting the space of candidate expressions from which to construct equations, and determining the threshold for including terms. Applying UDEs directly requires only specifying the network architecture, a choice to which results are generally robust. More fundamentally, the assumption underlying symbolic regression---that the true dynamics conform to neat algebraic expressions---cannot always be justified. In chemical reaction network modelling, even widely used functional forms such as Hill functions are adopted for analytical convenience rather than because they are expected to reflect the true underlying mechanisms. For such systems, omitting symbolic regression yields a simpler workflow without sacrificing interpretability.

Several directions for future work emerge naturally from this study. From a computational perspective, profile likelihood analysis remains substantially more expensive for UDEs than for purely mechanistic models because of the large number of neural network parameters. Developing algorithms specifically designed for this setting will therefore be important. Likewise, our ensemble-based approach to practical functional identifiability is simple, and primarily applies to single-input functions. Development of more principled approaches, with proven validity, for characterising functional identifiability is required. From a theoretical perspective, structural and practical functional identifiability remain relatively new concepts, and extending them to broader classes of differential equations and more general forms of learnt model components represents an important direction for future research.

More broadly, hybrid mechanistic/data-driven models offer a powerful alternative to modelling approaches that assume either complete mechanistic knowledge or none at all. UDEs provide a particularly attractive realisation of this philosophy by combining mechanistic structure with the flexibility to recover unknown dynamics from data. Our results demonstrate that this flexibility often comes at a surprisingly modest identifiability cost, especially when prior knowledge is incorporated through appropriate architectural constraints. More fundamentally, they show that the ability to learn unknown mechanisms is only scientifically meaningful when those mechanisms are identifiable. By providing a practical framework for assessing both functional and parametric identifiability, this work establishes a foundation for using UDEs not only as predictive models, but also as reliable tools for mechanistic discovery.


\section*{Acknowledgements}

The authors thank Frank Sch\"afer for helpful discussion, and would like to acknowledge the use of the University of Oxford Advanced Research Computing (ARC) facility in carrying out this work~\cite{richards_university_2015}. REB and TL acknowledge support of a grant from the Simons Foundation (MP-SIP-00001828).


\section*{Data availability}
Only synthetic data are used in the study. Code to reproduce all results in the study, and all simulated data, can be found at \url{https://github.com/TorkelE/UDE-Identifiability-Loman-Baker}.

\section*{Code availability}
The code used for carrying out all analyses, as well as for generating all figures, can be found at \url{https://github.com/TorkelE/UDE-Identifiability-Loman-Baker}.


\begin{singlespace}
\bibliographystyle{abbrvnat}
\bibliography{references} 

@misc{rackauckas_universal_2021,
	title = {Universal differential equations for scientific machine learning},
	doi = {10.48550/arXiv.2001.04385},
	publisher = {arXiv},
	author = {Rackauckas, Christopher and Ma, Yingbo and Martensen, Julius and Warner, Collin and Zubov, Kirill and Supekar, Rohit and Skinner, Dominic and Ramadhan, Ali and Edelman, Alan},
	month = nov,
	year = {2021},
	note = {arXiv:2001.04385 [cs]},
}

@phdthesis{loman_2022,
	title={How bacteria tune mixed positive/negative feedback loops to generate diverse gene expression dynamics},
	url={https://www.repository.cam.ac.uk/handle/1810/345440},
	DOI={10.17863/CAM.92862},
	school={University of Cambridge},
	author={Loman, Torkel},
	year={2022},
}

@article{villaverde_observability_2019,
	author  = {Villaverde, Alejandro F.},
	title   = {Observability and Structural Identifiability
             of Nonlinear Biological Systems},
	journal = {Complexity},
	year    = {2019},
	volume  = {2019},
	pages   = {8497093},
	doi     = {10.1155/2019/8497093},
	publisher = {Hindawi},
}

@article{philipps_non-negative_2024,
	series = {10th {IFAC} {Conference} on {Foundations} of {Systems} {Biology} in {Engineering} {FOSBE} 2024},
	title = {Non-negative universal differential equations with applications in systems biology},
	volume = {58},
	doi = {10.1016/j.ifacol.2024.10.005},
	number = {23},
	journal = {IFAC-PapersOnLine},
	author = {Philipps, Maren and Körner, Antonia and Vanhoefer, Jakob and Pathirana, Dilan and Hasenauer, Jan},
	month = jan,
	year = {2024},
	pages = {25--30},
}

@article{schmid_assessment_2025,
	title = {Assessment of uncertainty quantification in universal differential equations},
	volume = {383},
	doi = {10.1098/rsta.2024.0444},
	number = {2293},
	journal = {Philosophical Transactions of the Royal Society A: Mathematical, Physical and Engineering Sciences},
	author = {Schmid, Nina and Fernandes del Pozo, David and Waegeman, Willem and Hasenauer, Jan},
	month = apr,
	year = {2025},
	pages = {20240444},
}

@article{noordijk_rise_2024,
	title = {The rise of scientific machine learning: a perspective on combining mechanistic modelling with machine learning for systems biology},
	volume = {4},
	doi = {10.3389/fsysb.2024.1407994},
	journal = {Frontiers in Systems Biology},
	author = {Noordijk, Ben and Garcia Gomez, Monica L. and ten Tusscher, Kirsten H. W. J. and de Ridder, Dick and van Dijk, Aalt D. J. and Smith, Robert W.},
	month = aug,
	year = {2024},
	pages = {1407994},
}

@article{raissi_physics-informed_2019,
	title = {Physics-informed neural networks: {A} deep learning framework for solving forward and inverse problems involving nonlinear partial differential equations},
	volume = {378},
	doi = {10.1016/j.jcp.2018.10.045},
	journal = {Journal of Computational Physics},
	author = {Raissi, M. and Perdikaris, P. and Karniadakis, G. E.},
	month = feb,
	year = {2019},
	pages = {686--707},
}

@article{loman_catalyst_2023,
	title = {Catalyst: {Fast} and flexible modeling of reaction networks},
	volume = {19},
	doi = {10.1371/journal.pcbi.1011530},
	number = {10},
	journal = {PLOS Computational Biology},
	author = {Loman, Torkel E. and Ma, Yingbo and Ilin, Vasily and Gowda, Shashi and Korsbo, Niklas and Yewale, Nikhil and Rackauckas, Chris and Isaacson, Samuel A.},
	month = oct,
	year = {2023},
	pages = {e1011530},
}

@article{rackauckas_differentialequationsjl_2017,
	title = {{DifferentialEquations.jl} – {A} performant and feature-rich ecosystem for solving differential equations in {Julia}},
	volume = {5},
	doi = {10.5334/jors.151},
	number = {1},
	journal = {Journal of Open Research Software},
	author = {Rackauckas, Christopher and Nie, Qing},
	month = may,
	year = {2017},
	pages = {15},
}

@misc{noauthor_scimlmodelingtoolkitneuralnetsjl_2025,
	author = {Micluta-Campeanu, Sebastian and Rackauckas, Christopher},
	title = {{SciML}/{ModelingToolkitNeuralNets}.jl},
	url = {https://github.com/SciML/ModelingToolkitNeuralNets.jl},
	publisher = {SciML Open Source Scientific Machine Learning},
	month = sep,
	year = {2025},
}

@misc{dixit2023optimization,
	title={{Optimization.jl}: {A} unified optimization package},
	author={Dixit, Vaibhav Kumar and Rackauckas, Christopher},
	note={Zenodo, doi:10.5281/zenodo.7738525},
	year={2023},
}

@article{mogensen_optim_2018,
	title = {Optim: {A} mathematical optimization package for {Julia}},
	volume = {3},
	doi = {10.21105/joss.00615},
	number = {24},
	journal = {Journal of Open Source Software},
	author = {Mogensen, Patrick K. and Riseth, Asbjørn N.},
	month = apr,
	year = {2018},
	pages = {615},
}

@inproceedings{kingma_adam_2017,
	title = {Adam: {A} method for stochastic optimization},
	booktitle = {3rd International Conference on Learning Representations ({ICLR})},
	year = {2015},
	url = {https://arxiv.org/abs/1412.6980},
	author = {Kingma, Diederik P. and Ba, Jimmy},
}

@article{baker_mechanistic_2018,
	title = {Mechanistic models versus machine learning, a fight worth fighting for the biological community?},
	volume = {14},
	doi = {10.1098/rsbl.2017.0660},
	number = {5},
	journal = {Biology Letters},
	author = {Baker, Ruth E. and Peña, Jose-Maria and Jayamohan, Jayaratnam and Jérusalem, Antoine},
	month = may,
	year = {2018},
	pages = {20170660},
}

@misc{pal_lux_2022,
	title = {Lux: {E}xplicit parameterization of deep neural networks in {Julia}},
	url = {https://github.com/LuxDL/Lux.jl},
	author = {Pal, Avik},
	month = may,
	year = {2022},
	note = {Zenodo, doi:10.5281/zenodo.7808903},
}

@book{feinberg_foundations_2019,
	address = {Cham},
	series = {Applied {Mathematical} {Sciences}},
	title = {Foundations of {Chemical} {Reaction} {Network} {Theory}},
	volume = {202},
	isbn = {978-3-030-03857-1 978-3-030-03858-8},
	publisher = {Springer International Publishing},
	author = {Feinberg, Martin},
	year = {2019},
	doi = {10.1007/978-3-030-03858-8},
}

@article{avram_advancing_2024,
	title = {Advancing mathematical epidemiology and chemical reaction network theory via synergies between them},
	volume = {26},
	doi = {10.3390/e26110936},
	number = {11},
	journal = {Entropy},
	author = {Avram, Florin and Adenane, Rim and Neagu, Mircea},
	month = nov,
	year = {2024},
	pages = {936},
}

@incollection{brauer_compartmental_2008,
	title = {Compartmental models in epidemiology},
	volume = {1945},
	series = {Lecture Notes in Mathematics},
	publisher = {Springer},
	doi = {10.1007/978-3-540-78911-6_2},
	booktitle = {Mathematical Epidemiology},
	author = {Brauer, Fred},
	year = {2008},
	pages = {19--79},
}

@book{warnatz_combustion_1999,
	address = {Berlin, Heidelberg},
	title = {Combustion},
	isbn = {978-3-642-98029-9 978-3-642-98027-5},
	publisher = {Springer},
	author = {Warnatz, Jürgen and Maas, Ulrich and Dibble, Robert W.},
	year = {1999},
	doi = {10.1007/978-3-642-98027-5},
}

@article{elowitz_synthetic_2000,
	title = {A synthetic oscillatory network of transcriptional regulators},
	volume = {403},
	doi = {10.1038/35002125},
	number = {6767},
	journal = {Nature},
	author = {Elowitz, Michael B. and Leibler, Stanislas},
	month = jan,
	year = {2000},
	pages = {335--338},
}

@ARTICLE{elgazzar2025,
	AUTHOR={El-Gazzar, Ahmed  and van Gerven, Marcel },
	TITLE={Universal differential equations as a unifying modeling language for neuroscience},
	JOURNAL={Frontiers in Computational Neuroscience},
	VOLUME={19},
	PAGES={1677930},
	YEAR={2025},
	DOI={10.3389/fncom.2025.1677930},
}

@article{zou_application_2020,
	title = {Application of pharmacokinetic-pharmacodynamic modeling in drug delivery: {D}evelopment and challenges},
	volume = {11},
	doi = {10.3389/fphar.2020.00997},
	journal = {Frontiers in Pharmacology},
	author = {Zou, Huixi and Banerjee, Parikshit and Leung, Sharon Shui Yee and Yan, Xiaoyu},
	month = jul,
	year = {2020},
	pages = {997},
}

@article{hahl_comparison_2016,
	title = {A comparison of deterministic and stochastic modeling approaches for biochemical reaction systems: {O}n fixed points, means, and modes},
	volume = {7},
	doi = {10.3389/fgene.2016.00157},
	journal = {Frontiers in Genetics},
	author = {Hahl, Sayuri K. and Kremling, Andreas},
	month = aug,
	year = {2016},
	pages = {157},
}

@article{jo_neural_2025,
	title = {Neural network-based parameter estimation for nonautonomous differential equations with discontinuous signals},
	journal = {SIAM Journal on Applied Mathematics},
	volume = {86},
	pages = {327--347},
	doi = {10.1137/25M1741340},
	author = {Jo, Hyeontae and Josić, Krešimir and Kim, Jae Kyoung},
	year = {2026},
}

@misc{gunawardena_chemical_nodate,
	title = {Chemical reaction network theory for in-silico biologists},
	note = {Lecture notes, Bauer Center for Genomics Research, Harvard University},
	url = {http://vcp.med.harvard.edu/papers/crnt.pdf},
	author = {Gunawardena, Jeremy},
	year = {2003},
}

@article{goodwin_oscillatory_1965,
	title = {Oscillatory behavior in enzymatic control processes},
	volume = {3},
	doi = {10.1016/0065-2571(65)90067-1},
	journal = {Advances in Enzyme Regulation},
	author = {Goodwin, Brian C.},
	month = jan,
	year = {1965},
	pages = {425--437},
}

@article{goentoro_incoherent_2009,
	title = {The incoherent feedforward loop can provide fold-change detection in gene regulation},
	volume = {36},
	doi = {10.1016/j.molcel.2009.11.018},
	number = {5},
	journal = {Molecular Cell},
	author = {Goentoro, Lea and Shoval, Oren and Kirschner, Marc and Alon, Uri},
	month = dec,
	year = {2009},
	pages = {894--899},
}

@article{liu_limited_1989,
	title = {On the limited memory {BFGS} method for large scale optimization},
	volume = {45},
	doi = {10.1007/BF01589116},
	number = {1},
	journal = {Mathematical Programming},
	author = {Liu, Dong C. and Nocedal, Jorge},
	month = aug,
	year = {1989},
	pages = {503--528},
}

@article{boros_center_2018,
	title = {The center problem for the {Lotka} reactions with generalized mass-action kinetics},
	volume = {17},
	doi = {10.1007/s12346-017-0243-2},
	number = {2},
	journal = {Qualitative Theory of Dynamical Systems},
	author = {Boros, Balázs and Hofbauer, Josef and Müller, Stefan and Regensburger, Georg},
	year = {2018},
	pages = {403--410},
}

@article{field_oscillations_1974,
	title = {Oscillations in chemical systems. {IV}. {Limit} cycle behavior in a model of a real chemical reaction},
	volume = {60},
	doi = {10.1063/1.1681288},
	number = {5},
	journal = {The Journal of Chemical Physics},
	author = {Field, Richard J. and Noyes, Richard M.},
	month = mar,
	year = {1974},
	pages = {1877--1884},
}

@article{philipps_current_2025,
	title = {Current state and open problems in universal differential equations for systems biology},
	volume = {11},
	doi = {10.1038/s41540-025-00550-w},
	number = {1},
	journal = {npj Systems Biology and Applications},
	author = {Philipps, Maren and Schmid, Nina and Hasenauer, Jan},
	month = aug,
	year = {2025},
	pages = {101},
}

@inproceedings{dugas_incorporating_2000,
	title = {Incorporating second-order functional knowledge for better option pricing},
	volume = {13},
	url = {https://proceedings.neurips.cc/paper_files/paper/2000/hash/44968aece94f667e4095002d140b5896-Abstract.html},
	booktitle = {Advances in {Neural} {Information} {Processing} {Systems}},
	publisher = {MIT Press},
	author = {Dugas, Charles and Bengio, Yoshua and Bélisle, François and Nadeau, Claude and Garcia, René},
	year = {2000},
}

@article{schmiester_petabinteroperable_2021,
	title = {{PEtab}—interoperable specification of parameter estimation problems in systems biology},
	volume = {17},
	doi = {10.1371/journal.pcbi.1008646},
	number = {1},
	journal = {PLOS Computational Biology},
	author = {Schmiester, Leonard and Schälte, Yannik and Bergmann, Frank T. and Camba, Tacio and Dudkin, Erika and Egert, Janine and Fröhlich, Fabian and Fuhrmann, Lara and Hauber, Adrian L. and Kemmer, Svenja and Lakrisenko, Polina and Loos, Carolin and Merkt, Simon and Müller, Wolfgang and Pathirana, Dilan and Raimundez, Elba and Refisch, Lukas and Rosenblatt, Marcus and Stapor, Paul L. and Städter, Philipp and Wang, Dantong and Wieland, Franz-Georg and Banga, Julio R. and Timmer, Jens and Villaverde, Alejandro F. and Sahle, Sven and Kreutz, Clemens and Hasenauer, Jan and Weindl, Daniel},
	month = jan,
	year = {2021},
	pages = {e1008646},
}

@article{structidjl,
	author  = {Dong, R. and Goodbrake, C. and Harrington, H. A. and Pogudin, G.},
	doi     = {10.1137/22M1469067},
	title   = {Differential Elimination for Dynamical Models via Projections with Applications to Structural Identifiability},
	journal = {SIAM Journal on Applied Algebra and Geometry},
	year    = {2023},
	volume  = {7},
	number  = {1},
	pages   = {194-235},
}

@article{kreutz_profile_2013,
	title = {Profile likelihood in systems biology},
	volume = {280},
	doi = {10.1111/febs.12276},
	number = {11},
	journal = {The FEBS Journal},
	author = {Kreutz, Clemens and Raue, Andreas and Kaschek, Daniel and Timmer, Jens},
	year = {2013},
	pages = {2564--2571},
}

@article{de_rooij_conditional_2025,
	title = {Conditional universal differential equations capture population dynamics and interindividual variation in {C}-peptide production},
	volume = {11},
	doi = {10.1038/s41540-025-00570-6},
	number = {1},
	journal = {npj Systems Biology and Applications},
	author = {de Rooij, Max and van Riel, Natal A. W. and O’Donovan , Shauna D.},
	month = jul,
	year = {2025},
	pages = {84},
}

@inproceedings{jia_neural_2020,
	title = {Neural jump stochastic differential equations},
	booktitle = {Advances in Neural Information Processing Systems},
	volume = {32},
	publisher = {Curran Associates},
	year = {2019},
	url = {https://papers.nips.cc/paper/9177-neural-jump-stochastic-differential-equations},
	author = {Jia, Junteng and Benson, Austin R.},
}

@inproceedings{oh_stable_2025,
	title = {Stable neural stochastic differential equations in analyzing irregular time series data},
	booktitle = {12th International Conference on Learning Representations ({ICLR})},
	year = {2024},
	url = {https://proceedings.iclr.cc/paper_files/paper/2024/hash/a61023ce36d21010f1423304f8ec49af-Abstract-Conference.html},
	author = {Oh, YongKyung and Lim, Dong-Young and Kim, Sungil},
}

@misc{dandekar_bayesian_2022,
	title = {Bayesian neural ordinary differential equations},
	doi = {10.48550/arXiv.2012.07244},
	publisher = {arXiv},
	author = {Dandekar, Raj and Chung, Karen and Dixit, Vaibhav and Tarek, Mohamed and Garcia-Valadez, Aslan and Vemula, Krishna Vishal and Rackauckas, Chris},
	month = feb,
	year = {2022},
	note = {arXiv:2012.07244 [cs]},
}

@article{bolibar_universal_2023,
	title = {Universal differential equations for glacier ice flow modelling},
	volume = {16},
	doi = {10.5194/gmd-16-6671-2023},
	number = {22},
	journal = {Geoscientific Model Development},
	author = {Bolibar, Jordi and Sapienza, Facundo and Maussion, Fabien and Lguensat, Redouane and Wouters, Bert and Pérez, Fernando},
	month = nov,
	year = {2023},
	pages = {6671--6687},
}

@misc{martinez_comparative_2024,
	title = {A comparative study of {NeuralODE} and universal {ODE} approaches to solving {Chandrasekhar} white dwarf equation},
	doi = {10.48550/arXiv.2410.14998},
	publisher = {arXiv},
	author = {{Vazquez Martinez}, Raymundo and Dandekar, Raj Abhijit and Dandekar, Rajat and Panat, Sreedath},
	month = oct,
	year = {2024},
	note = {arXiv:2410.14998 [cs]},
}

@article{PEtabBioinformatics2025,
	title={{PEtab.jl}: advancing the efficiency and utility of dynamic modelling},
	doi={10.1093/bioinformatics/btaf497},
	author={Persson, Sebastian and Fr{\"o}hlich, Fabian and Grein, Stephan and Loman, Torkel and Ognissanti, Damiano and Hasselgren, Viktor and Hasenauer, Jan and Cvijovic, Marija},
	journal={Bioinformatics},
	volume={41},
	number={9},
	pages={btaf497},
	year={2025},
	publisher={Oxford University Press},
}

@article{buckner_recovering_2024,
	title = {Recovering complex ecological dynamics from time series using state-space universal dynamic equations},
	journal = {Communications Earth \& Environment},
	volume = {7},
	pages = {112},
	doi = {10.1038/s43247-025-03130-2},
	author = {Buckner, Jack H. and Meunier, Zechariah D. and Arroyo-Esquivel, Jorge and Fitzpatrick, Nathan and Greiner, Ariel and McManus, Lisa C. and Watson, James R.},
	year = {2026},
}

@inproceedings{zhu_neural_2021,
	title = {Neural delay differential equations},
	booktitle = {9th International Conference on Learning Representations ({ICLR})},
	year = {2021},
	url = {https://iclr.cc/virtual/2021/poster/3120},
	author = {Zhu, Qunxi and Guo, Yao and Lin, Wei},
}

@article{taghvaei_fractional_2020,
	title = {Fractional {SIR} epidemiological models},
	journal = {Scientific Reports},
	volume = {10},
	pages = {20882},
	year = {2020},
	doi = {10.1038/s41598-020-77849-7},
	author = {Taghvaei, Amirhossein and Georgiou, Tryphon T. and Norton, Larry and Tannenbaum, Allen},
}

@misc{cranmer_interpretable_2023,
	title = {Interpretable machine learning for science with {PySR} and {SymbolicRegression.jl}},
	doi = {10.48550/arXiv.2305.01582},
	publisher = {arXiv},
	author = {Cranmer, Miles},
	month = may,
	year = {2023},
	note = {arXiv:2305.01582 [astro-ph]},
}

@article{browning_cautionary_2025,
	title = {A cautionary tale of model misspecification and identifiability},
	journal = {Bulletin of Mathematical Biology},
	volume = {88},
	number = {1},
	pages = {5},
	doi = {10.1007/s11538-025-01573-4},
	author = {Browning, Alexander P. and Flegg, Jennifer A. and Murphy, Ryan J.},
	year = {2026},
}

@article{kaye_impact_2024,
	title = {The impact of natural climate variability on the global distribution of \emph{{Aedes aegypti}}: a mathematical modelling study},
	volume = {8},
	doi = {10.1016/S2542-5196(24)00238-9},
	number = {12},
	journal = {The Lancet Planetary Health},
	author = {Kaye, Alexander R. and Obolski, Uri and Sun, Lantao and Hart, William S. and Hurrell, James W. and Tildesley, Michael J. and Thompson, Robin N.},
	month = dec,
	year = {2024},
	pages = {e1079--e1087},
}

@misc{richards_university_2015,
	title = {University of {Oxford} {Advanced} {Research} {Computing}},
	doi = {10.5281/zenodo.22558},
	note = {Zenodo, doi:10.5281/zenodo.22558},
	author = {Richards, Andrew},
	month = aug,
	year = {2015},
}

@misc{loman_learning_2026,
	title         = "Learning functional components of {PDEs} from data using
                   neural networks",
	author        = "Loman, Torkel E and Salmaniw, Yurij and {Leon Villares}, Antonio
                   and Carrillo, Jose A and Baker, Ruth E",
	note          = "arXiv:2602.13174",
	month         =  feb,
	year          =  2026,
}

@misc{loman2026structuralfunctionalidentifiabilitymodel,
	title={Structural functional identifiability and model discovery in differential equation models},
	author={Torkel E Loman and Alexander P Browning and Ruth E Baker},
	year={2026},
	note={arXiv:2606.30289},
}

@article{simpson2026parameter,
	author = {Simpson, Matthew J. and Baker, Ruth E.},
	title = {Parameter Identifiability, Parameter Estimation, and Model Prediction for Differential Equation Models},
	journal = {SIAM Review},
	volume = {68},
	number = {1},
	pages = {153-171},
	doi = {10.1137/24M1667968},
	year = {2026},
}

@article{LJUNG1994265,
	title = {On global identifiability for arbitrary model parametrizations},
	journal = {Automatica},
	volume = {30},
	number = {2},
	pages = {265-276},
	year = {1994},
	doi = {10.1016/0005-1098(94)90029-9},
	author = {Lennart Ljung and Torkel Glad},
}
\end{singlespace}


\pagebreak


\appendix
\clearpage
\setcounter{page}{1}
\setcounter{section}{0}
\setcounter{equation}{0}
\setcounter{figure}{0}
\setcounter{table}{0}
\renewcommand{\thepage}{S\arabic{page}}
\renewcommand{\thesection}{S\arabic{section}}
\renewcommand{\theequation}{S\arabic{equation}}
\renewcommand{\thefigure}{S\arabic{figure}}
\renewcommand{\thetable}{S\arabic{section}.\arabic{table}}
\numberwithin{equation}{section}
\renewcommand{\appendixname}{}
\makeatletter
\renewcommand{\@seccntformat}[1]{\csname the#1\endcsname.\quad}
\makeatother
\begin{center}
    {\Large\bfseries Functional and parametric identifiability for \\ universal differential equations --- Supplementary Information}\\
\date{}
\end{center}


\section{Encoding neural network constraints}\label{ssection:detailed_methodology_nn_constraints}

In this section, we describe various ways to encode constraints on fitted functions through model architecture. With the exception of the non-negativity constraint (which is used throughout all UDE models), these are only used in Section~\ref{section:results_linear_pathway}. 

\paragraph{Non-negativity constraints.} By ensuring that all functions in the last layer of the neural network are non-negative, we can ensure that all fitted functions are non-negative (fulfilling a strict requirement for chemical reaction network rate functions). We utilise this for all fitted neural networks, primarily through the softplus function.

\paragraph{Lower and upper bounds constraints.} To encode the constraint that a fitted neural network is bounded by $U_\text{min}<U(\vec{X}; \vec{\theta})<U_\text{max}$, we first ensure that the final output function of the neural network is bounded by $(0,1)$ (here, using the fast sigmoid function). Next, we replace the rate
\begin{equation}
\ce{$\dots$ ->C[$U(\vec{X}; \vec{\theta})$] $\dots$} 
\end{equation}
with
\begin{equation}
\ce{$\dots$ ->C[$U_\text{min} + (U_\text{max} - U_\text{min})U(\vec{X}; \vec{\theta})$] $\dots$} 
\end{equation}
which ensures the prescribed bounds.

\paragraph{Monotonicity constraints.} To ensure that a fitted neural network is a monotonically increasing function of the input, we ensure that all activation functions are monotonically increasing and all neural network weight parameters are nonnegative. This means that, for a neural network of the form
\begin{equation}
    U(x;\vec{w},\vec{b})
    = w_n\,\sigma_{n-1}\!\left(w_{n-1}\,\sigma_{n-2}\!\left(\cdots\,
      \sigma_1\!\left(w_1x+b_1\right)\cdots\right)+b_{n-1}\right)+b_n,
\end{equation}
with activation functions $\sigma_i$, weights $\vec{w}=(w_1,\ldots,w_n)$, and biases $\vec{b}=(b_1,\ldots,b_n)$, we have set all weights $w_i = \text{softplus}(v_i)$ (where $v_i$ are the unmodified weights). Here, $\text{softplus}(x) = \log(1 + e^x)$ is simply used as a smooth function that ensures that all considered weight values are positive. For the monotonically decreasing case, we will instead use non-positive weights, i.e.~$w_i = -\text{softplus}(v_i)$, and also ensure that the neural network has an odd number of layers.

When all activation functions are taken to be equal, this approach for enforcing monotonicity has a tendency to bias the fitted functions towards the shape of the activation function. To prevent this effect when enforcing monotonicity, we used a mixture of different activation functions (softplus and the fast sigmoid).

\pagebreak

\section{Data quality variation scans}\label{ssection:detailed_methodology_datavar_scans}

\begin{figure}
    \centering
    \includegraphics[width=\textwidth]{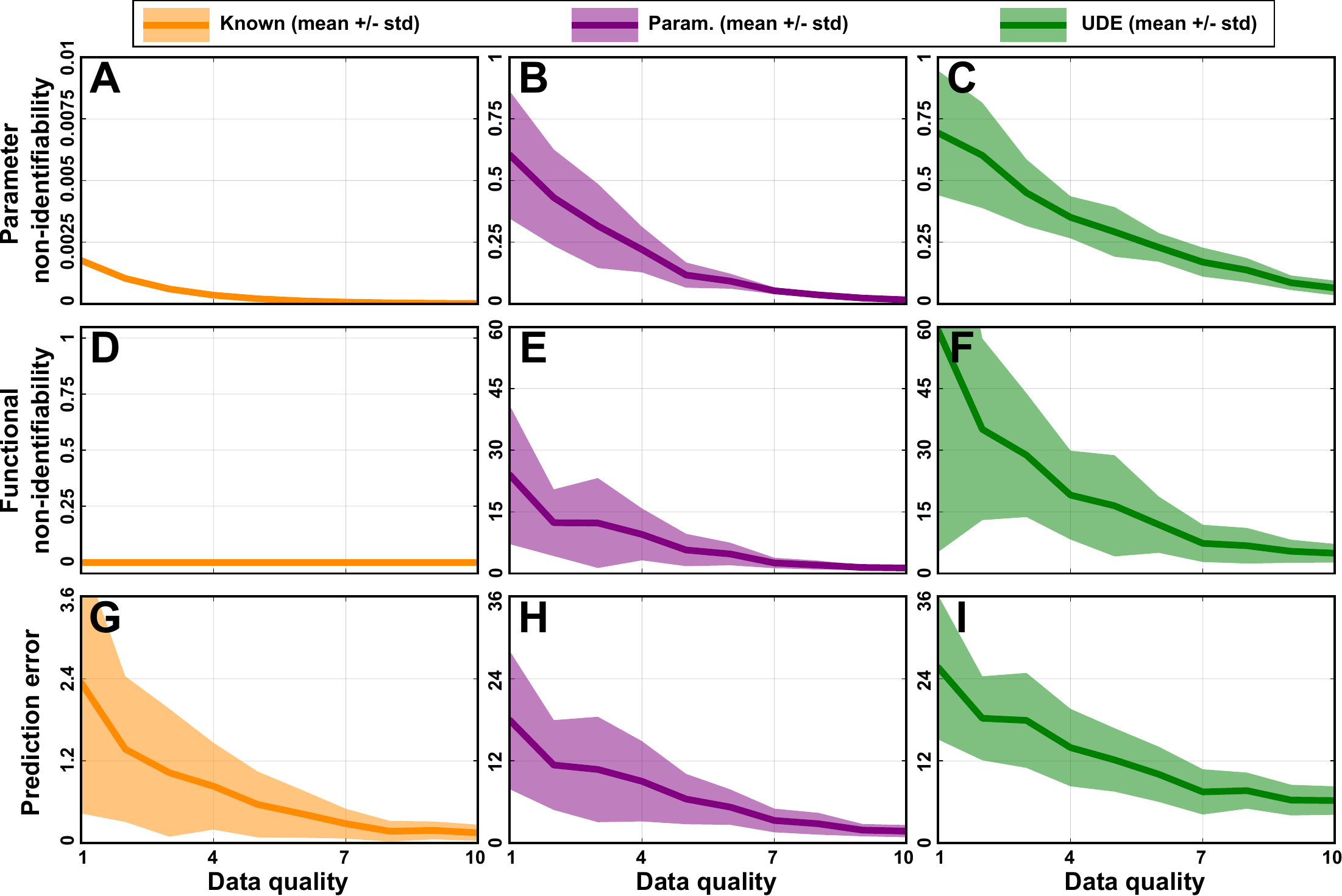}
    \caption{\textbf{Parametric and functional practical identifiability are correlated with predictive power in an extended self-activation loop model.} For the extended self-activation loop model, we consider three cases: (1) $A(Y;\vec{p})=vY^n/(Y^n + K^n)$ with $\vec{p}=(v,K,n)$ known (yellow; A, D, G); (2) $A(Y;\vec{p})=vY^n/(Y^n + K^n)$ with $\vec{p}=(v,K,n)$ unknown (purple; B, E, H); and (3) $A(Y;\vec{p})=U(Y;\vec{\theta})$ (green; C, F, I). In each case, the mechanistic parameter, $d$, is unknown. We examine how parametric identifiability (A--C), functional identifiability (D--F), and predictive performance (G--I) vary with data quality (measures are detailed in Supplementary Section~\ref{ssection:detailed_methodology_datavar_scans}). Data quality increases along the $x$-axis, from sparsely sampled, high-noise observations to densely sampled, low-noise observations (synthetic data generation detailed in Supplementary Section~\ref{ssection:detailed_methodology_datavar_scans}). Across the three model cases, identifiability and predictive performance improve as data quality increases. We note that, across the full data range, all three measures have similar magnitudes for the parameterised and UDE model. Next, we note that the three measures deteriorate in a coordinated manner across all models (functional identifiability for the known model is constant $\equiv0$, as the function is known). Example fitted trajectories, likelihood profiles, and function ensemble plots, for the three models at different data qualities, are provided in Supplementary Figure~\ref{sfigure:extended_sal_sample_fits}.}
    \label{sfigure:extended_sal_datavar}
\end{figure}

In Figure~\ref{figure:mosquitoes_data_quality_scan} and Supplementary Figure~\ref{sfigure:extended_sal_datavar}, we investigate how identifiability and predictive error are affected by data quality. In each case, we generate data across ten different data qualities (denoted $1, 2, \dots, 10$). In each case, we generate ten different sets of synthetic data (as described in Section~\ref{section:detailed_methodology_synthetic_data}) where $n_\text{sample}$ and $\sigma$ are varied between $n_{\text{sample},1}$ and $\sigma_1$ and $n_{\text{sample},10}$ and $\sigma_{10}$, respectively. For each level of data quality, we generate $40$ different datasets. We fit the respective model to each dataset, and for each, compute three different measures as described in the following subsections.


\subsection{Parametric identifiability measure}\label{ssection:detailed_methodology_paramident_measures}

For likelihood-based parameter fitting, the $95\%$ confidence interval for a parameter value can be obtained from the intersections of the profile likelihood curve with a horizontal line at $y\approx -1.92$ (plotted in each likelihood profile figure in the paper). We use the length of this interval as a measure of parameter identifiability.


\subsection{Functional identifiability measure}\label{ssection:detailed_methodology_funcident_measures}

In Section~\ref{section:detailed_methodology_identifiability_analysis}, we describe how we generate fitted functions for the ensemble plots of functional identifiability. For any two functions, we can compute the $L^2$ distance between them over the range of input values explored by the dataset. We do this by computing the $L^2$ distance using a densely sampled ($N = 200$) grid of values across the interval. The functional identifiability measure is then taken as the mean $L^2$ distance across the ensemble.


\subsection{Predictive power measure}\label{ssection:detailed_methodology_predpow_measures}

In this work, we assess a model's predictive power through its ability to accurately predict the values of unobserved species. Specifically, we simulate the model for the fitted parameter set and the ground-truth parameter set. For each unobserved species, we then compute the $L^2$ distance between its fitted and ground-truth trajectories. The sum of these $L^2$ distances across all unobserved species is taken as the prediction error, and is computed using a densely sampled set of time points ($N = 100$) across the simulation time interval.

\pagebreak

\section{Supplementary analyses}\label{ssection:supp_analysis}


\subsection{Alternative optimisation approaches}\label{ssection:optimisation_approaches}

As described in Section~\ref{section:detailed_methodology_UDE_fitting}, we use a combination of Adam and L-BFGS to fit all models to data. However, while this is a commonly used optimisation scheme for UDEs, simpler ones also perform well. In Supplementary Figure~\ref{sfigure:alt_optimisation}, we show three different alternative optimisation schemes that are also able to fit a simple UDE to data. These are L-BFGS, Nelder-Mead, and conjugate gradient. In practice, if high performance is required, it is prudent to analyse the performance of multiple schemes. 


\begin{figure}[htbp]
   \centering
   \includegraphics[width=\textwidth]{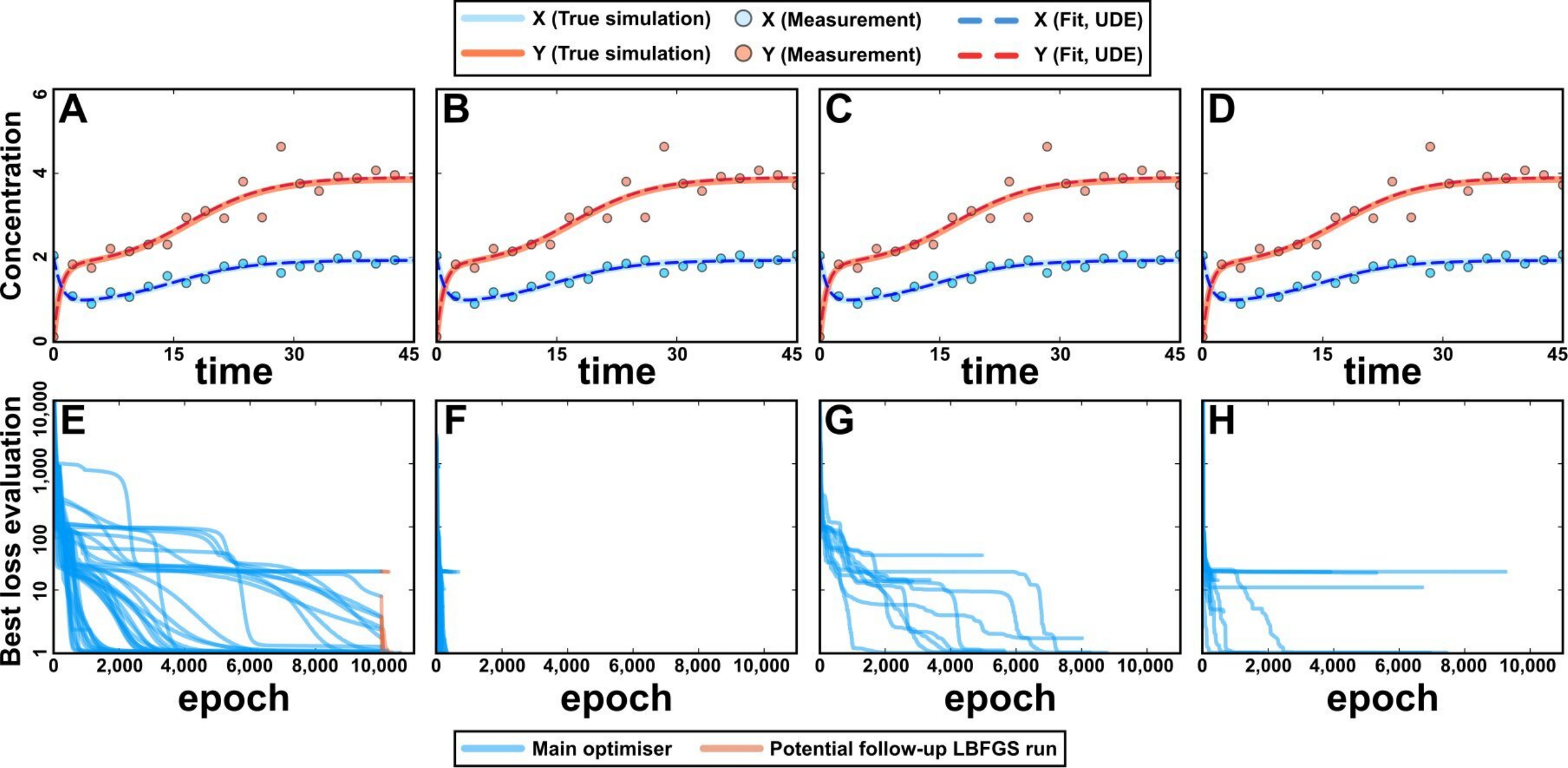}
   \caption{\textbf{Several optimisation schemes perform well for UDE parameter fitting.} Throughout this work, we use Adam followed by an L-BFGS tuning run to fit all models to data (Section~\ref{section:detailed_methodology_UDE_fitting}). Here, we briefly confirm that alternative, simpler, schemes also perform well. We try Adam and L-BFGS (A,E), L-BFGS only (B,F), Nelder-Mead (C,G), and conjugate gradient (D,H). (A--D) Using each optimisation scheme, we fit the extended self-activation loop UDE to synthetic data (using the same settings as in Figure~\ref{figure:extended_sal_pract_ident}B, where both $X$ and $Y$ are measured). All schemes yield a good model-to-data fit. (E--H) Each case in A-D is performed using 40 individual optimisation runs with random initial guesses. Here, for each individual run, we plot the progression of the best loss function evaluation throughout that run. While the rate of convergence is different for each scheme, they all reach the same minimum (all $y$-axes are normalised so that $1$ is equivalent to the best loss function evaluations achieved across all runs across all optimisation schemes).}
   \label{sfigure:alt_optimisation}
\end{figure}


\subsection{Alternative thresholds for ensemble plots} \label{ssection:ensemble_plot_thresholds}

In Section~\ref{section:detailed_methodology_identifiability_analysis}, we describe how, when generating the ensemble plots used to assess practical functional identifiability, we use a threshold based on the likelihood of the data under the ground-truth model to decide which functional forms to plot. As noted there, this threshold is a convenient choice for this article, but it is not available in real-world applications, where the ground truth is unknown. Here, in Supplementary Figure~\ref{sfigure:alt_ens_plot_thres}, we consider an alternative threshold: the one traditionally used in profile likelihood analysis to determine the $95\%$ confidence interval of a parameter. This alternative threshold yields similar results. We also note that plots of the loss function's progression during model fitting can be used to assess the validity of any chosen threshold.


 \begin{figure}[htbp]
    \centering
    \includegraphics[width=\textwidth]{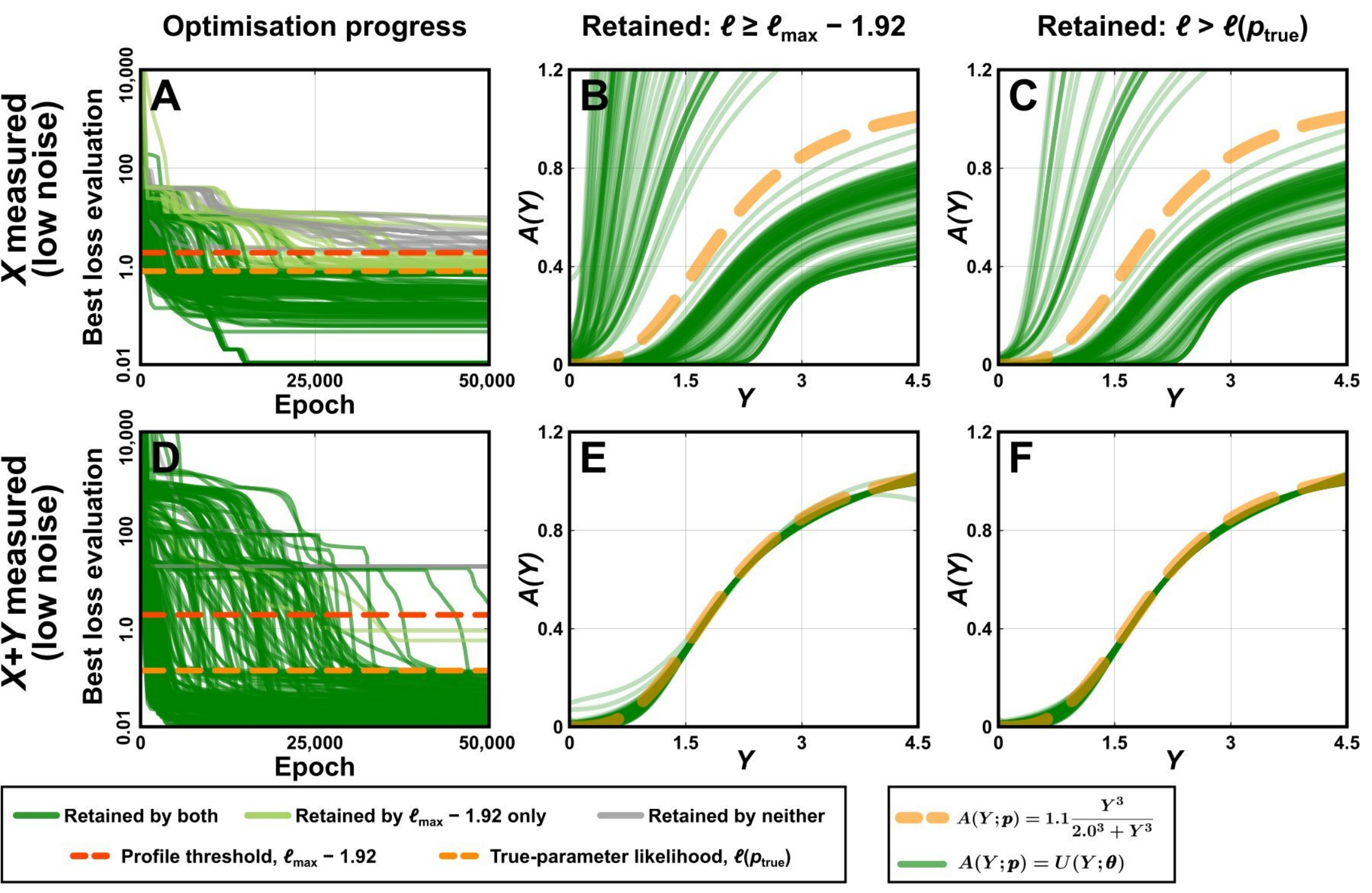}
    \caption{\textbf{Alternative thresholds for ensemble plots are available, and yield consistent results.} In Section~\ref{section:detailed_methodology_identifiability_analysis}, we described how the ensemble plots used to assess practical functional identifiability were generated, and how we used a threshold based on the ground-truth model's likelihood to decide which functions to plot. Here, we demonstrate that alternative thresholds can be used. Specifically, we use $\ell_{\max} - 1.92$, where $\ell_{\max}$ is the maximum log-likelihood. (A) For the fitted extended self-activation loop model (specifically the case used in Figure~\ref{figure:extended_sal_pract_ident}A,D,G), we plot the progression of the best loss function evaluation throughout each run. The $y$-axis shows the negative log-likelihood value of the best fit achieved in that run. A total of $200$ individual runs are shown. The two ensemble thresholds are shown as red and orange dashed horizontal lines. Optimisation runs that are retained by the $\ell_{\max} - 1.92$-based, but not the ground-truth-based, threshold are shown in light green. (B,C) The ensembles retained by the $\ell_{\max} - 1.92$ threshold (B) and by the ground-truth threshold (C). The differences between the two plots are minor, and both clearly indicate practical non-identifiability. (D--F) Repeats of the analysis in A--C, but using the dataset in Figure~\ref{figure:extended_sal_pract_ident}B,E,H. Again, both thresholds consistently indicate practical identifiability.}
    \label{sfigure:alt_ens_plot_thres}
\end{figure}

For both thresholds, the ensemble plots give consistent descriptions of the model's practical functional identifiability. The analysis also shows how loss function progression plots can be used to compare different thresholds, and to choose an appropriate threshold in applications. These tools are not necessarily principled in the way that the $95\%$ confidence interval of profile likelihood analysis is for the parameters of standard differential equation models. However, they provide an effective diagnostic of the degree of practical identifiability of a given UDE.


\subsection{Learning noise parameters from data} \label{ssection:noise_parameter_learning}

In Section~\ref{section:detailed_methodology_UDE_fitting}, we described the likelihood-based loss function used to fit our models to data (Equation~\eqref{eq:likelihood_formula}). It depends on the values of the mechanistic parameters, $\vec{p}$, and of the neural network parameters, $\vec{\theta}$. It also depends on the value of $\sigma$, which sets the scale of the noise (Equation~\eqref{eq:measurement_formula}) and which we assume to be known. It is, however, possible to treat $\sigma$ as unknown, and to learn it by optimising the same loss function with $\sigma$ as an additional free parameter. To simplify the analysis, we omit this step elsewhere in this work. Here, however, we carry out this analysis for the extended self-activation loop model (Section~\ref{section:results_extended_sa}), in Supplementary Figure~\ref{sfigure:learning_noise_params}. The results are similar to those obtained when $\sigma$ is assumed to be known: the accuracy with which $d$ is recovered is not significantly affected, and $\sigma$ itself is estimated close to its true value.


 \begin{figure}[htbp]
    \centering
    \includegraphics[width=\textwidth]{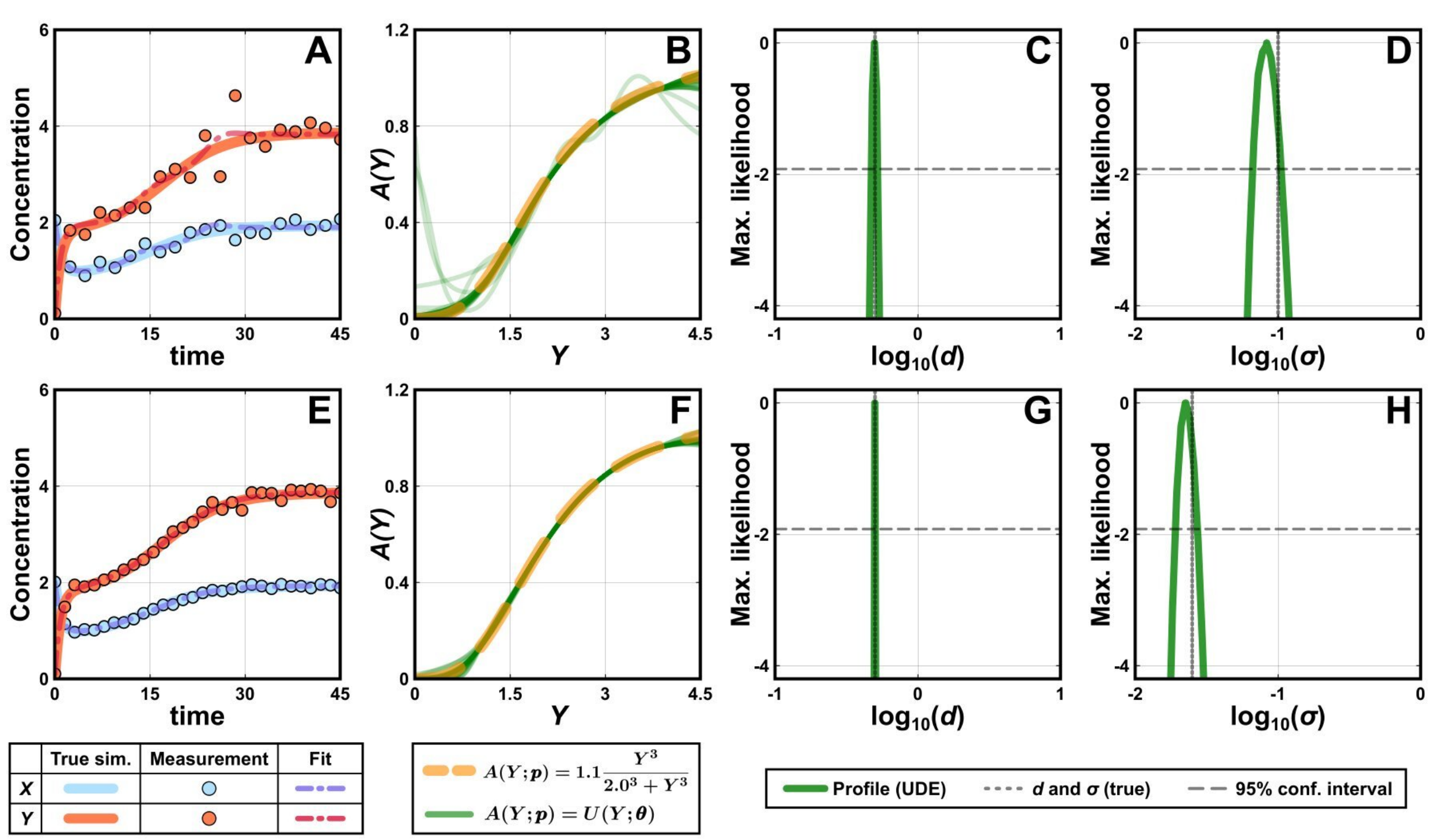}
    \caption{\textbf{A UDE's noise parameters can be learnt in tandem with its mechanistic parameters.} In Supplementary Section~\ref{ssection:noise_parameter_learning}, we describe how the model fitting workflows used in this article can also be used to learn parameters of the measurement noise process. Specifically, we here assume that the shape of the noise is known (a normal distribution that scales with the magnitude of the measured variable), but that the scaling factor $\sigma$ is unknown. (A--D) Repeat of the fit in Figure~\ref{figure:extended_sal_pract_ident}B,E,H, but where we also learn $\sigma$. (A) The ground-truth simulation, data, and fitted simulation (corresponding to Figure~\ref{figure:extended_sal_pract_ident}B). (B) The ensemble of fitted functions $A(Y)$ (corresponding to Figure~\ref{figure:extended_sal_pract_ident}H). There is some additional variability in the learnt functional forms, but practical functional identifiability remains. (C) Profile likelihood for the mechanistic parameter $d$ (corresponding to Figure~\ref{figure:extended_sal_pract_ident}E). The profile is not affected by also learning $\sigma$. (D) Profile likelihood for the noise parameter $\sigma$. The estimate of $\sigma$ is close to its true value, and the profile likelihood is narrow, indicating that $\sigma$ is practically identifiable. (E--H) Repeat of the analysis in A--D, but using a different dataset with more measurements ($30$ rather than $20$ time points) and a smaller noise parameter ($\sigma = 0.025$).}
    \label{sfigure:learning_noise_params}
 \end{figure}

\pagebreak

\section{Sample code for chemical reaction network/UDE implementation and training}\label{ssection:code_example}

Here, we provide a simple demo script showing how the Lux.jl, Catalyst.jl, and ModelingToolkitNeuralNets.jl packages (all implemented in the Julia programming language) can be combined to generate a chemical reaction network/UDE for the extended self-activation loop (Section~\ref{section:results_extended_sa}).

\enlargethispage{0.5cm}


\begin{juliacode}[numbers=left]
# Import packages.
using Catalyst, ModelingToolkitNeuralNets, Lux

# Create a symbolic NN representation using a selected architecture.
nn_arch = Lux.Chain(
    Lux.Dense(1 => 3, Lux.softplus, use_bias = false),
    Lux.Dense(3 => 3, Lux.softplus, use_bias = false),
    Lux.Dense(3 => 1, Lux.softplus, use_bias = false)
)
@SymbolicNeuralNetwork U, theta = nn_arch
A(x) = U(x, theta)[1]

# Creates a chemical reaction network/UDE. 
# The dollar sign is used to insert U into the chemical reaction network block.
extended_sa = @reaction_network begin
    $A(Y), 0 --> X
    X, 0 --> Y
    d, (X,Y) --> 0
end
\end{juliacode}
This model can now be trained using standard workflows. Alternatively, one can use PEtab.jl, which implements specific UDE fitting frameworks~\cite{PEtabBioinformatics2025}. For example, assuming the data are in ``data.csv'', that $X$ only is measured, that observational noise is normally distributed with standard deviation $X/10$. Then, we can run $10{,}000$ iterations of Adam with $20$ independent starts using:
\begin{juliacode}[numbers=left]
# Import packages.
using CSV, DataFrames, Optimisers, PEtab

# Create a likelihood-based parameter fitting problem using PEtab.jl
ps_est = [PEtabMLParameter(:theta), PEtabParameter(:d)]
observables = [PEtabObservable(:obs_X, :X, extended_sa.X/10)]
measurements = CSV.read("data.csv", DataFrame)
petab_model = PEtabModel(extended_sa, observables, measurements, ps_est)
petab_prob = PEtabODEProblem(petab_model)

# Fits the model to the data.
options = OptimisersOptions(iterations = 10000)
fit_result = calibrate_multistart(petab_prob, Optimisers.Adam(1e-3), 20; options)
\end{juliacode}

\pagebreak

\section{Goodwin oscillator model description}\label{ssection:gwo} 

The Goodwin oscillator model provides an example of a circuit able to generate oscillatory behaviours~\cite{goodwin_oscillatory_1965}. We model this using the following chemical reaction network
\begin{equation}
    \ce{ \text{\o} ->C[$R(P)$] $M$},
    \quad
    \ce{ $M$ ->C[$d_m$] \text{\o}},
    \quad
    \ce{ \text{\o} ->C[$p_e M$] $E$},
    \quad
    \ce{ $E$ ->C[$d_e$] \text{\o}},
    \quad
    \ce{ \text{\o} ->C[$p_p E$] $P$},
    \quad
    \ce{ $P$ ->C[$d_p$] \text{\o}},
\end{equation}
which yields the following ODE model
\begin{align}
    \frac{\mathrm{d}M}{\mathrm{d}t} &= R(P) - d_m M,\\
    \frac{\mathrm{d}E}{\mathrm{d}t} &= p_e M  - d_e E, \\
    \frac{\mathrm{d}P}{\mathrm{d}t} &= p_p E  - d_p P.
\end{align}
We assume that the repression function $R$ takes the following form:
\begin{equation*}
    R(P;\vec{p}) = vK^n/(P^n + K^n),
\end{equation*}
with $\vec{p} = (v,K,n)$. 

Throughout all analyses, we consider various model cases, including: (1) where $R(P;\vec{p})$ and $\vec{p}$ are known; (2) where $R(P;\vec{p})$ is known but $\vec{p}$ is unknown; and (3) where $R(P;\vec{p})$ is approximated using the universal function approximator $U(P;\vec{\theta})$. In all cases, we assume that the only unknown non-functional parameter is $d_p$, with $p_e$, $p_p$, $d_m$, and $d_e$ all being known (ensuring the model has the same number of unknown mechanistic parameters as the model to which it is compared). We also note that, when the full set of mechanistic parameters is unknown, the model exhibits structural non-identifiability.

\pagebreak

\section{Supplementary tables}\label{ssection:sup_tables}


\begin{table}[H]
    \centering
    \small
    \renewcommand{\arraystretch}{1.3}
    \setlength{\tabcolsep}{4pt}

    \caption{\textbf{Synthetic data description.}
        For each synthetic dataset used throughout this work, the table gives the initial conditions, ground-truth parameter values, final simulation time, number of samples, and noise level. The measured species are
        specified in the corresponding figures.}
    \label{stable:data_settings}
    
    \begin{tabularx}{\textwidth}{
        @{}
        >{\centering\arraybackslash}m{2.6cm}
        >{\centering\arraybackslash}m{2.2cm}
        >{\centering\arraybackslash}X
        >{\centering\arraybackslash}m{1.4cm}
        >{\centering\arraybackslash}m{1.4cm}
        >{\centering\arraybackslash}m{1.4cm}
        @{}
    }
        \toprule
        \textbf{Figure}
        &
        \textbf{Initial condition}
        &
        \textbf{Ground-truth parameters}
        &
        $t_{\mathrm{end}}$
        &
        \textbf{Samples}
        &
        \textbf{Noise}
        \\
        \midrule

        Fig.~\ref{figure:simple_sal_struc_nonident}
        &
        $X=0.18$
        &
        $(v,K,n,d)=(0.6,0.3,3.0,0.5)$
        &
        $15$
        &
        $40$
        &
        $0.05$
        \\

        \midrule

        \makecell{Fig.~\ref{figure:extended_sal_pract_ident}\\A, B, D, E, G, H}
        &
        $(X,Y)=(2.0,0.1)$
        &
        $(v,K,n,d)=(1.1,2.0,3,0.5)$
        &
        $45$
        &
        $20$
        &
        $0.1$
        \\

        \midrule

        \makecell{Fig.~\ref{figure:extended_sal_pract_ident}\\C, F, I}
        &
        $(X,Y)=(2.0,0.1)$
        &
        $(v,K,n,d)=(1.1,2.0,3,0.5)$
        &
        $45$
        &
        $20$
        &
        $0.3$
        \\

        \midrule

        Fig.~\ref{figure:ic_feedforward_2nn}
        &
        $\begin{gathered}
            (X,Y,Z)\\
            =(0.0,0.0,1.0)
        \end{gathered}$
        &
        $\begin{gathered}
            (v_A,K_A,v_R,K_R,n,d)\\
            =(1.0,0.4,5.0,0.8,3,0.5)
        \end{gathered}$
        &
        $8$
        &
        $20$
        &
        $0.25$
        \\

        \midrule

        Fig.~\ref{figure:ln_pathway_nn_constraints}
        &
        $\begin{gathered}
            (X_1,X_2,X_3,X_4)\\
            =(1.0,1.0,1.0,1.0)
        \end{gathered}$
        &
        $\begin{gathered}
            (v_0,v,K,n,d)\\
            =(0.05,1.0,0.4,7.0,1.0)
        \end{gathered}$
        &
        $5$
        &
        $50$
        &
        $0.15$
        \\

        \midrule

        Fig.~\ref{figure:modified_sir_misspecification}
        &
        $\begin{gathered}
            (S,I,R)\\
            =(99.0,1.0,0.0)
        \end{gathered}$
        &
        $\begin{gathered}
            (\beta,\gamma,\alpha_1,\alpha_2)\\
            =(0.003,0.005,-0.05,0.5)
        \end{gathered}$
        &
        $1000$
        &
        $50$
        &
        $0.2$
        \\

        \bottomrule
    \end{tabularx}
    
\end{table}


\small
\renewcommand{\arraystretch}{1.25}
\setlength{\tabcolsep}{3pt}

\captionof{table}{\textbf{Synthetic data description for the Supplementary Figures.} Continuation of Table~\ref{stable:data_settings}, giving the synthetic data settings used in the supplementary figures.}
\label{stable:data_settings_sfigs}
\begin{xltabular}{\textwidth}{
    @{}
    >{\centering\arraybackslash}p{0.17\textwidth}
    >{\centering\arraybackslash}p{0.14\textwidth}
    >{\centering\arraybackslash}X
    >{\centering\arraybackslash}p{0.07\textwidth}
    >{\centering\arraybackslash}p{0.12\textwidth}
    >{\centering\arraybackslash}p{0.12\textwidth}
    @{}
}

\toprule
\textbf{Figure}
&
\textbf{Initial condition}
&
\textbf{Ground-truth parameters}
&
$t_{\mathrm{end}}$
&
\textbf{Samples}
&
\textbf{Noise}
\\
\midrule
\endfirsthead

\toprule
\textbf{Figure}
&
\textbf{Initial condition}
&
\textbf{Ground-truth parameters}
&
$t_{\mathrm{end}}$
&
\textbf{Samples}
&
\textbf{Noise}
\\
\midrule
\endhead

\midrule
\multicolumn{6}{r}{\textit{Continued on next page}}\\
\endfoot

\bottomrule
\endlastfoot

Fig.~\ref{sfigure:alt_optimisation}
&
$(X,Y)=(2.0,0.1)$
&
$(v,K,n,d)=(1.1,2.0,3,0.5)$
&
$45.0$
&
$20$
&
$0.1$
\\
\midrule

\makecell{Fig.~\ref{sfigure:alt_ens_plot_thres}\\
A, B, C}
&
$(X,Y)=(2.0,0.1)$
&
$(v,K,n,d)=(1.1,2.0,3,0.5)$
&
$45.0$
&
$20$
&
$0.1$
\\
\midrule

\makecell{Fig.~\ref{sfigure:alt_ens_plot_thres}\\
D, E, F}
&
$(X,Y)=(2.0,0.1)$
&
$(v,K,n,d)=(1.1,2.0,3,0.5)$
&
$45.0$
&
$20$
&
$0.1$
\\
\midrule

\makecell{Fig.~\ref{sfigure:learning_noise_params}\\
A, B, C, D}
&
$(X,Y)=(2.0,0.1)$
&
$(v,K,n,d)=(1.1,2.0,3,0.5)$
&
$45.0$
&
$20$
&
$0.1$
\\
\midrule

\makecell{Fig.~\ref{sfigure:learning_noise_params}\\
E, F, G, H}
&
$(X,Y)=(2.0,0.1)$
&
$(v,K,n,d)=(1.1,2.0,3,0.5)$
&
$45.0$
&
$30$
&
$0.025$
\\
\midrule

\makecell{Fig.~\ref{sfigure:extended_sal_pred_uncertain}\\
A, B, D, E, G, H}
&
$(X,Y)=(2.0,0.1)$
&
$(v,K,n,d)=(1.1,2.0,3,0.5)$
&
$45.0$
&
$20$
&
$0.1$
\\
\midrule

\makecell{Fig.~\ref{sfigure:extended_sal_pred_uncertain}\\
C, F, I}
&
$(X,Y)=(2.0,0.1)$
&
$(v,K,n,d)=(1.1,2.0,3,0.5)$
&
$45.0$
&
$20$
&
$0.3$
\\
\midrule

Fig.~\ref{sfigure:extended_sal_pract_ident_maxdata}
&
$(X,Y)=(2.0,0.1)$
&
$(v,K,n,d)=(1.1,2.0,3,0.5)$
&
$45.0$
&
$80$
&
$0.01$
\\
\midrule

\makecell{Fig.~\ref{sfigure:extended_sal_pract_ident_rep1}\\
A, B, D, E, G, H}
&
$(X,Y)=(2.0,0.1)$
&
$(v,K,n,d)=(1.1,2.0,3,0.5)$
&
$45.0$
&
$20$
&
$0.1$
\\
\midrule

\makecell{Fig.~\ref{sfigure:extended_sal_pract_ident_rep1}\\
C, F, I}
&
$(X,Y)=(2.0,0.1)$
&
$(v,K,n,d)=(1.1,2.0,3,0.5)$
&
$45.0$
&
$20$
&
$0.3$
\\
\midrule

\makecell{Fig.~\ref{sfigure:extended_sal_pract_ident_rep2}\\
A, B, D, E, G, H}
&
$(X,Y)=(2.0,0.1)$
&
$(v,K,n,d)=(1.1,2.0,3,0.5)$
&
$45.0$
&
$20$
&
$0.1$
\\
\midrule

\makecell{Fig.~\ref{sfigure:extended_sal_pract_ident_rep2}\\
C, F, I}
&
$(X,Y)=(2.0,0.1)$
&
$(v,K,n,d)=(1.1,2.0,3,0.5)$
&
$45.0$
&
$20$
&
$0.3$
\\
\midrule

\makecell{Fig.~\ref{sfigure:extended_sal_pract_ident_alt_ps}\\
A, B, D, E, G, H}
&
$(X,Y)=(2.0,0.1)$
&
$(v,K,n,d)=(2.3,2.5,4,0.5)$
&
$45.0$
&
$20$
&
$0.1$
\\
\midrule

\makecell{Fig.~\ref{sfigure:extended_sal_pract_ident_alt_ps}\\
C, F, I}
&
$(X,Y)=(2.0,0.1)$
&
$(v,K,n,d)=(2.3,2.5,4,0.5)$
&
$45.0$
&
$20$
&
$0.3$
\\
\midrule

\makecell{Fig.~\ref{sfigure:gwo_pract_ident}\\
A, B, D, E, G, H}
&
$(M,E,P)=(1.0,0.1,0.1)$
&
$\begin{gathered}
(v,K,n,d_m,p_e,d_e,p_p,d_p)\\
=(1.0,0.4,4,1.0,0.25,0.2,0.1,0.05)
\end{gathered}$
&
$30.0$
&
$30$
&
$0.05$
\\
\midrule

\makecell{Fig.~\ref{sfigure:gwo_pract_ident}\\
C, F, I}
&
$(M,E,P)=(1.0,0.1,0.1)$
&
$\begin{gathered}
(v,K,n,d_m,p_e,d_e,p_p,d_p)\\
=(1.0,0.4,4,1.0,0.25,0.2,0.1,0.05)
\end{gathered}$
&
$30.0$
&
$30$
&
$0.2$
\\
\midrule

\makecell{Fig.~\ref{sfigure:extended_sal_pract_ident_lowconc_ps}\\
A, B, D, E, G, H}
&
$(X,Y)=(2.0,0.1)$
&
$(v,K,n,d)=(1.6,1.5,4,0.5)$
&
$15.0$
&
$20$
&
$0.1$
\\
\midrule

\makecell{Fig.~\ref{sfigure:extended_sal_pract_ident_lowconc_ps}\\
C, F, I}
&
$(X,Y)=(2.0,0.1)$
&
$(v,K,n,d)=(1.6,1.5,4,0.5)$
&
$15.0$
&
$20$
&
$0.3$
\\
\midrule

\makecell{Fig.~\ref{sfigure:extended_sal_sample_fits}\\
A, E, I}
&
$(X,Y)=(2.0,0.1)$
&
$(v,K,n,d)=(1.1,2.0,4,0.5)$
&
$20.0$
&
$15$
&
$0.40$
\\
\midrule

\makecell{Fig.~\ref{sfigure:extended_sal_sample_fits}\\
B, F, J}
&
$(X,Y)=(2.0,0.1)$
&
$(v,K,n,d)=(1.1,2.0,4,0.5)$
&
$20.0$
&
$30$
&
$0.2$
\\
\midrule

\makecell{Fig.~\ref{sfigure:extended_sal_sample_fits}\\
C, G, K}
&
$(X,Y)=(2.0,0.1)$
&
$(v,K,n,d)=(1.1,2.0,4,0.5)$
&
$20.0$
&
$45$
&
$0.1$
\\
\midrule

\makecell{Fig.~\ref{sfigure:extended_sal_sample_fits}\\
D, H, L}
&
$(X,Y)=(2.0,0.1)$
&
$(v,K,n,d)=(1.1,2.0,4,0.5)$
&
$20.0$
&
$60$
&
$0.05$
\\
\midrule

\makecell{Fig.~\ref{sfigure:insects_example_fits}\\
A, E, I}
&
$(E,L,M)=(0.0,0.0,1.0)$
&
$\begin{gathered}
(r,d_E,k_{L,M},d_{L,1},d_{L,2},d_{L,3},k_{E,L},d_M)\\
=(0.3,0.05,0.15,0.05,0.1,0.2,0.2,0.1)
\end{gathered}$
&
$80.0$
&
$30$
&
$0.07$
\\
\midrule

\makecell{Fig.~\ref{sfigure:insects_example_fits}\\
B, F, J}
&
$(E,L,M)=(0.0,0.0,1.0)$
&
$\begin{gathered}
(r,d_E,k_{L,M},d_{L,1},d_{L,2},d_{L,3},k_{E,L},d_M)\\
=(0.3,0.05,0.15,0.05,0.1,0.2,0.2,0.1)
\end{gathered}$
&
$80.0$
&
$53$
&
$0.05$
\\
\midrule

\makecell{Fig.~\ref{sfigure:insects_example_fits}\\
C, G, K}
&
$(E,L,M)=(0.0,0.0,1.0)$
&
$\begin{gathered}
(r,d_E,k_{L,M},d_{L,1},d_{L,2},d_{L,3},k_{E,L},d_M)\\
=(0.3,0.05,0.15,0.05,0.1,0.2,0.2,0.1)
\end{gathered}$
&
$80.0$
&
$77$
&
$0.03$
\\
\midrule

\makecell{Fig.~\ref{sfigure:insects_example_fits}\\
D, H, L}
&
$(E,L,M)=(0.0,0.0,1.0)$
&
$\begin{gathered}
(r,d_E,k_{L,M},d_{L,1},d_{L,2},d_{L,3},k_{E,L},d_M)\\
=(0.3,0.05,0.15,0.05,0.1,0.2,0.2,0.1)
\end{gathered}$
&
$80.0$
&
$100$
&
$0.01$
\\
\midrule

\makecell{Fig.~\ref{sfigure:ic_feedforward_measurement_type_identifiability}\\
A, B}
&
$(X,Y,Z)=(0.0,0.0,0.0)$
&
$\begin{gathered}
(v_A,K_A,v_R,K_R,n,d)\\
=(1.0,0.6,5.0,0.8,5,0.5)
\end{gathered}$
&
$8.0$
&
$30$
&
$0.15$
\\
\midrule

\makecell{Fig.~\ref{sfigure:ic_feedforward_measurement_type_identifiability}\\
C.1}
&
$(X,Y,Z)=(0.0,0.0,0.0)$
&
$\begin{gathered}
(v_A,K_A,v_R,K_R,n,d)\\
=(1.0,0.6,5.0,0.8,5,0.5)
\end{gathered}$
&
$8.0$
&
$30$
&
$0.15$
\\
\midrule

\makecell{Fig.~\ref{sfigure:ic_feedforward_measurement_type_identifiability}\\
C.2}
&
$(X,Y,Z)=(1.0,0.0,0.0)$
&
$\begin{gathered}
(v_A,K_A,v_R,K_R,n,d)\\
=(1.0,0.6,5.0,0.8,5,0.5)
\end{gathered}$
&
$8.0$
&
$30$
&
$0.15$
\\
\midrule

\makecell{Fig.~\ref{sfigure:ic_feedforward_measurement_type_identifiability}\\
C.3}
&
$(X,Y,Z)=(3.0,3.0,0.0)$
&
$\begin{gathered}
(v_A,K_A,v_R,K_R,n,d)\\
=(1.0,0.6,5.0,0.8,5,0.5)
\end{gathered}$
&
$8.0$
&
$30$
&
$0.15$
\\
\midrule

\makecell{Fig.~\ref{sfigure:ic_feedforward_measurement_type_identifiability}\\
C.4}
&
$(X,Y,Z)=(0.0,3.0,0.0)$
&
$\begin{gathered}
(v_A,K_A,v_R,K_R,n,d)\\
=(1.0,0.6,5.0,0.8,5,0.5)
\end{gathered}$
&
$8.0$
&
$30$
&
$0.15$
\\
\midrule

Fig.~\ref{sfigure:ic_feedforward_2nn_repeat}
&
$(X,Y,Z)=(0.0,0.0,1.0)$
&
$\begin{gathered}
(v_A,K_A,v_R,K_R,n,d)\\
=(1.0,0.4,5.0,0.8,3,0.5)
\end{gathered}$
&
$8.0$
&
$20$
&
$0.25$
\\
\midrule

Fig.~\ref{sfigure:ln_pathway_nn_constraints_repeat2}
&
$\begin{gathered}
(X_1,X_2,X_3,X_4)\\
=(1.0,1.0,1.0,1.0)
\end{gathered}$
&
$(v_0,v,K,n,d)=(0.05,1.0,0.4,7.0,1.0)$
&
$5.0$
&
$50$
&
$0.15$
\\
\midrule

Fig.~\ref{sfigure:ln_pathway_nn_constraints_repeat3}
&
$\begin{gathered}
(X_1,X_2,X_3,X_4)\\
=(1.0,1.0,1.0,1.0)
\end{gathered}$
&
$(v_0,v,K,n,d)=(0.05,1.0,0.4,7.0,1.0)$
&
$5.0$
&
$50$
&
$0.15$
\\
\midrule

Fig.~\ref{sfigure:gwo_nn_constraints_repeat1}
&
$(M,E,P)=(1.0,0.1,0.1)$
&
$\begin{gathered}
(v,K,n,d_m,p_e,d_e,p_p,d_p)\\
=(1.0,0.4,4.0,1.0,0.25,0.2,0.1,0.05)
\end{gathered}$
&
$30.0$
&
$20$
&
$0.25$
\\
\midrule

Fig.~\ref{sfigure:gwo_nn_constraints_repeat2}
&
$(M,E,P)=(1.0,0.1,0.1)$
&
$\begin{gathered}
(v,K,n,d_m,p_e,d_e,p_p,d_p)\\
=(1.0,0.4,4.0,1.0,0.25,0.2,0.1,0.05)
\end{gathered}$
&
$30.0$
&
$20$
&
$0.25$
\\
\midrule

Fig.~\ref{sfigure:modified_sir_misspecification_repeat}
&
$(S,I,R)=(99.0,1.0,0.0)$
&
$\begin{gathered}
(\beta,\gamma,\alpha_1,\alpha_2)\\
=(0.003,0.005,-0.05,0.5)
\end{gathered}$
&
$1000.0$
&
$50$
&
$0.2$
\\
\midrule

Fig.~\ref{sfigure:ic_feedforward_misspecification}
&
$(X,Y,Z)=(0.0,0.0,0.0)$
&
$\begin{gathered}
(v_A,K_A,v_R,K_R,n,d)\\
=(1.0,0.6,5.0,0.8,5,0.5)
\end{gathered}$
&
$10.0$
&
$25$
&
$0.15$
\\

\end{xltabular}

\normalsize


\begin{table}[htbp]
    \centering
    \small
    \renewcommand{\arraystretch}{1.45}
    \setlength{\tabcolsep}{5pt}

    \caption{\textbf{Synthetic data scan description.} For each scan over data quality, the table gives the initial conditions, ground-truth parameter values, final simulation time, and the numbers of samples and noise parameters corresponding to the lowest and highest data qualities.}
    \label{stable:data_variation_scan_settings}

    \begin{tabularx}{\textwidth}{
        @{}
        >{\centering\arraybackslash}m{2.3cm}
        >{\centering\arraybackslash}m{2.1cm}
        >{\raggedright\arraybackslash}X
        >{\centering\arraybackslash}m{1.0cm}
        >{\centering\arraybackslash}m{1.3cm}
        >{\centering\arraybackslash}m{1.3cm}
        @{}
    }
        \toprule
        \textbf{Figure}
        &
        \textbf{Initial condition}
        &
        \textbf{Ground-truth parameters}
        &
        $t_{\mathrm{end}}$
        &
        $(n_1,n_{10})$
        &
        $(\sigma_1,\sigma_{10})$
        \\
        \midrule

        Fig.~\ref{figure:mosquitoes_data_quality_scan}
        &
        $\begin{gathered}
            (E,L,M)\\[-2pt]
            =(0,0,1)
        \end{gathered}$
        &
        $\begin{gathered}
            r=0.3,\quad d_E=0.05,\quad k_{L,M}=0.15,\\
            d_{L,1}=0.05,\quad d_{L,2}=0.1,\quad d_{L,3}=0.2,\\
            k_{E,L}=0.2,\quad d_M=0.1
        \end{gathered}$
        &
        $80$
        &
        $(30,100)$
        &
        $(0.07,0.01)$
        \\

        \addlinespace[4pt]

        Fig.~\ref{sfigure:extended_sal_datavar}
        &
        $\begin{gathered}
            (X,Y)\\[-2pt]
            =(2.0,0.1)
        \end{gathered}$
        &
        $\begin{gathered}
            v=1.1,\quad K=2.0,\\
            n=4.0,\quad d=0.5
        \end{gathered}$
        &
        $20$
        &
        $(15,60)$
        &
        $(0.40,0.05)$
        \\

        \bottomrule
    \end{tabularx}

\end{table}


\begin{table}[htbp]
    \centering
    \small
    \renewcommand{\arraystretch}{1.45}
    \setlength{\tabcolsep}{5pt}

    \caption{\textbf{Ensemble plot mean $L^2$ distances.} For each ensemble plot in the main text, the mean $L^2$ distance between functional forms randomly sampled from the fitted ensemble, computed as described in Section~\ref{section:detailed_methodology_identifiability_analysis}.}
    \label{stable:figure_ensemble_mean_l2_distances}

    \begin{tabular}{
        @{}
        >{\centering\arraybackslash}m{2.3cm}
        >{\raggedright\arraybackslash}m{7.2cm}
        >{\centering\arraybackslash}m{3.0cm}
        @{}
    }
        \toprule
        \textbf{Figure} & \textbf{Model} & \textbf{Mean $L^2$ distance} \\
        \midrule

        Fig.~\ref{figure:simple_sal_struc_nonident}C & Parameterised & $0.0253$ \\
         & UDE & $3.98$ \\

        \midrule

        Fig.~\ref{figure:extended_sal_pract_ident}G & Parameterised & $1.01$ \\
         & UDE & $8.70$ \\

        \midrule

        Fig.~\ref{figure:extended_sal_pract_ident}H & Parameterised & $5.79\times10^{-4}$ \\
         & UDE & $0.129$ \\

        \midrule

        Fig.~\ref{figure:extended_sal_pract_ident}I & Parameterised & $5.97\times10^{-4}$ \\
         & UDE & $0.161$ \\

        \midrule

        Fig.~\ref{figure:ic_feedforward_2nn}C & $R$ Parameterised, $A$ Parameterised & $1.72\times10^{-3}$ \\
         & $R$ UDE, $A$ Parameterised & $3.28$ \\
         & $R$ Parameterised, $A$ UDE & $1.07$ \\
         & $R$ UDE, $A$ UDE & $1.07$ \\

        \midrule

        Fig.~\ref{figure:ic_feedforward_2nn}D & $R$ Parameterised, $A$ Parameterised & $8.74\times10^{-4}$ \\
         & $R$ UDE, $A$ Parameterised & $1.06$ \\
         & $R$ Parameterised, $A$ UDE & $1.53$ \\
         & $R$ UDE, $A$ UDE & $1.26$ \\

        \midrule

        Fig.~\ref{figure:ln_pathway_nn_constraints}F & Parameterised & $10.4$ \\
         & UDE (monotonic, bounded) & $0.744$ \\
         & UDE (monotonic) & $7.41$ \\
         & UDE & $3.66\times10^{13}$ \\

        \bottomrule
    \end{tabular}

\end{table}


\pagebreak

\section{Supplementary figures}\label{ssection:sup_figs}


\begin{figure}[H]
    \centering
    \includegraphics[width=\textwidth]{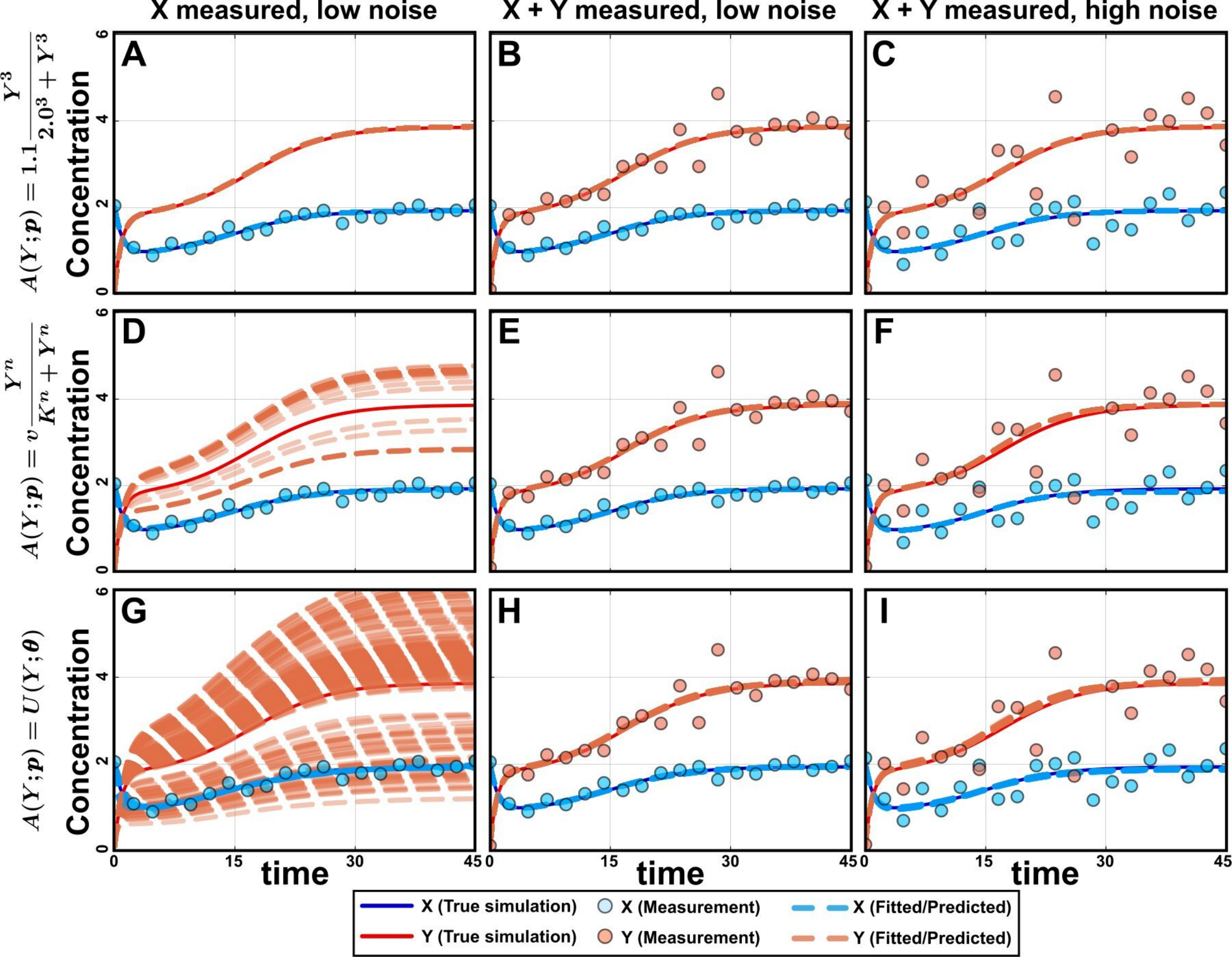}
    \caption{\textbf{Model non-identifiability indicates predictive uncertainty.} In Figure~\ref{figure:extended_sal_pract_ident}A--C, we compare the simulations of the best-fitting models to the ground-truth simulations and measured data. Here, we instead show multiple fitted simulations. Each plot shows the ground-truth simulations (solid lines), measured data (dots), and the ensemble of fitted simulations (dashed lines). We plot the simulations for all fits achieving the same loss threshold used to determine the functions to plot in the ensemble plots (i.e.~a likelihood at least equal to that of the ground-truth case). (A,D,G) Predictions for the low-noise case with only $X$ measured (Figure~\ref{figure:extended_sal_pract_ident}A). (B,E,H) Predictions for the low-noise case with both species measured (Figure~\ref{figure:extended_sal_pract_ident}B). (C,F,I) Predictions for the high-noise case with both species measured (Figure~\ref{figure:extended_sal_pract_ident}C). Separate plots show the fits for when both $A(Y;\vec{p}) = vY^n/(Y^n+K^n)$ and $\vec{p}$ are known (A--C), $A(Y;\vec{p}) = vY^n/(Y^n+K^n)$ is known but $\vec{p}$ is unknown (D--F), and when $A(Y;\vec{p})$ is approximated as $U(Y;\vec{\theta})$ (G--I). When both species are measured, the fitted simulations, unsurprisingly, align with the ground-truth ones. Next, when only $X$ is measured, the only case where the fitted predictions align with the ground-truth $Y$ values is the case which achieves parametric and functional identifiability (which is when $A(Y;\vec{p})$ is fully known). In the other two cases, a range of different $Y$ trajectories is possible, yielding predictive uncertainty.}
    \label{sfigure:extended_sal_pred_uncertain}
\end{figure}


\begin{figure}
    \centering
    \includegraphics[width=\textwidth]{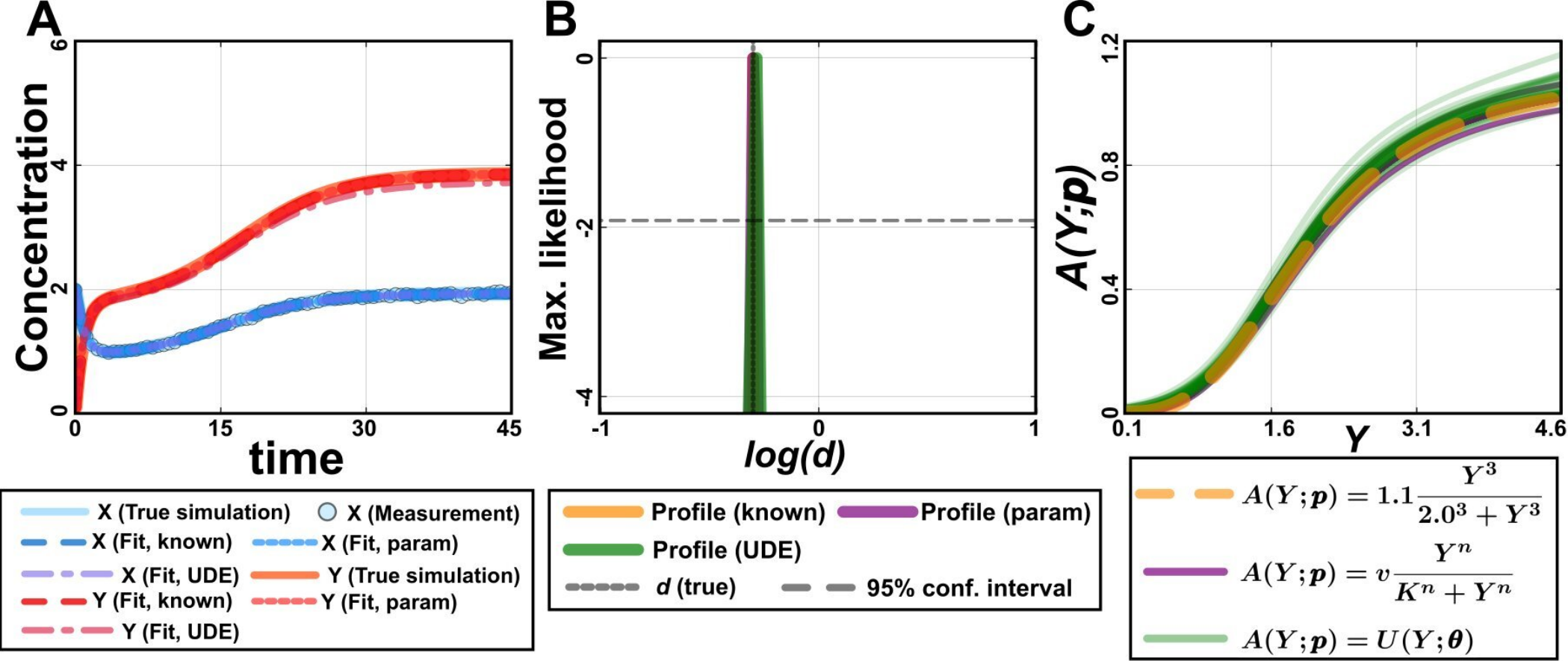}
    \caption{\textbf{The extended self-activation loop model can achieve identifiability when $X$ only is measured, given sufficiently high quality data.} In Figure~\ref{figure:extended_sal_pract_ident}, we noted that the models where $A$ was fitted either as a Hill function with unknown parameters, or using a UDE, were non-identifiable in the case where only $X$ was measured. Here, we confirm that practical identifiability can be recovered by increasing the data quality. We do so by generating a new dataset where only $X$ is measured, but with significantly more samples and less noise (A). We confirm that, for each of the three cases describing our different levels of knowledge of $A$, we can fit the model to the data and achieve both parametric (B) and functional (C) identifiability.}
    \label{sfigure:extended_sal_pract_ident_maxdata}
\end{figure}


\begin{figure}
    \centering
    \includegraphics[width=\textwidth]{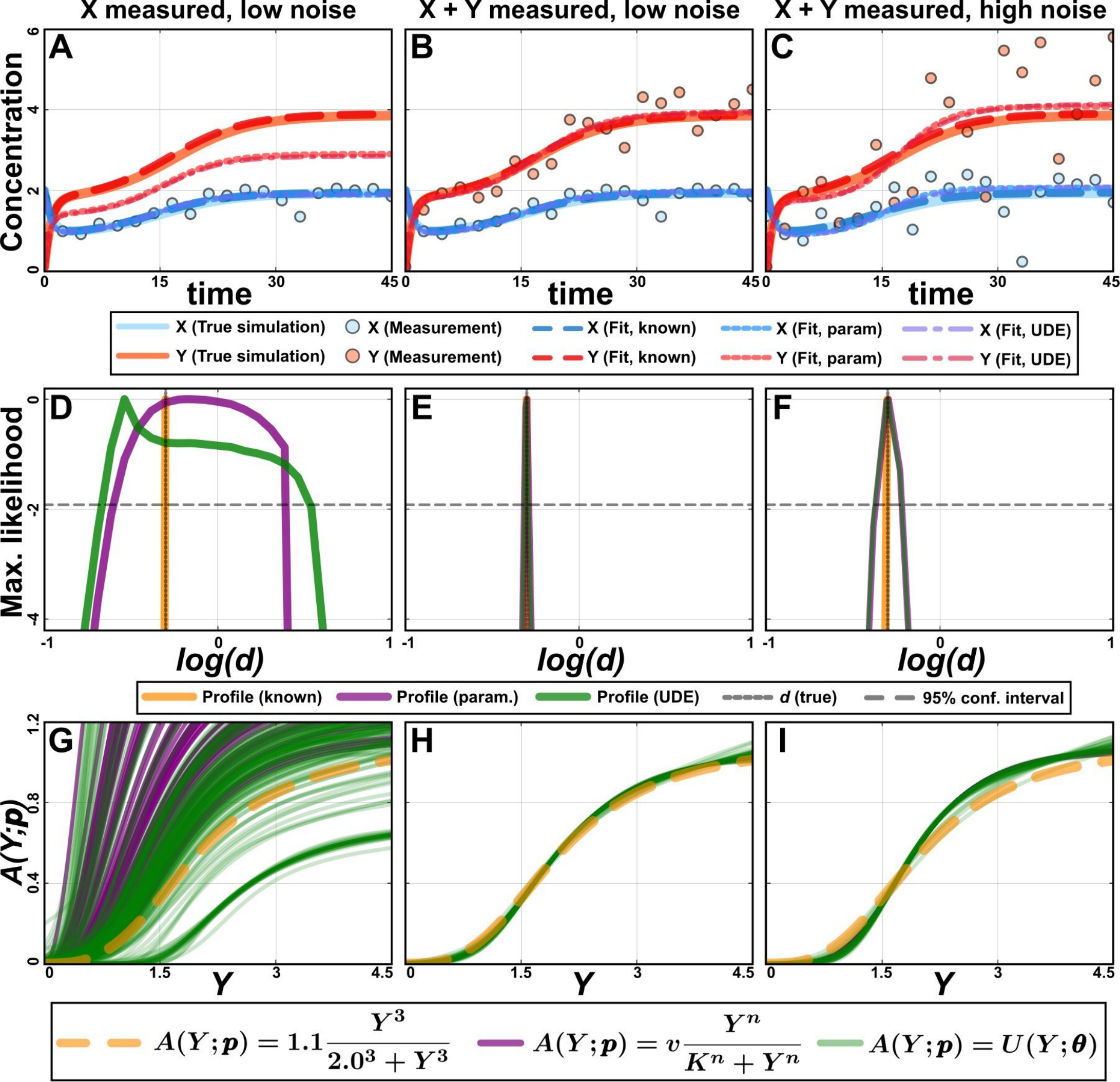}
    \caption{\textbf{Practical identifiability analysis of the extended self-activation loop: Repeat 1.} This figure is identical to Figure~\ref{figure:extended_sal_pract_ident} of the main text, however, here we use a different set of (equally distributed) measurements. This figure exhibits similar results and conclusions as Figure~\ref{figure:extended_sal_pract_ident}.}
    \label{sfigure:extended_sal_pract_ident_rep1}
\end{figure}

\begin{figure}
    \centering
    \includegraphics[width=\textwidth]{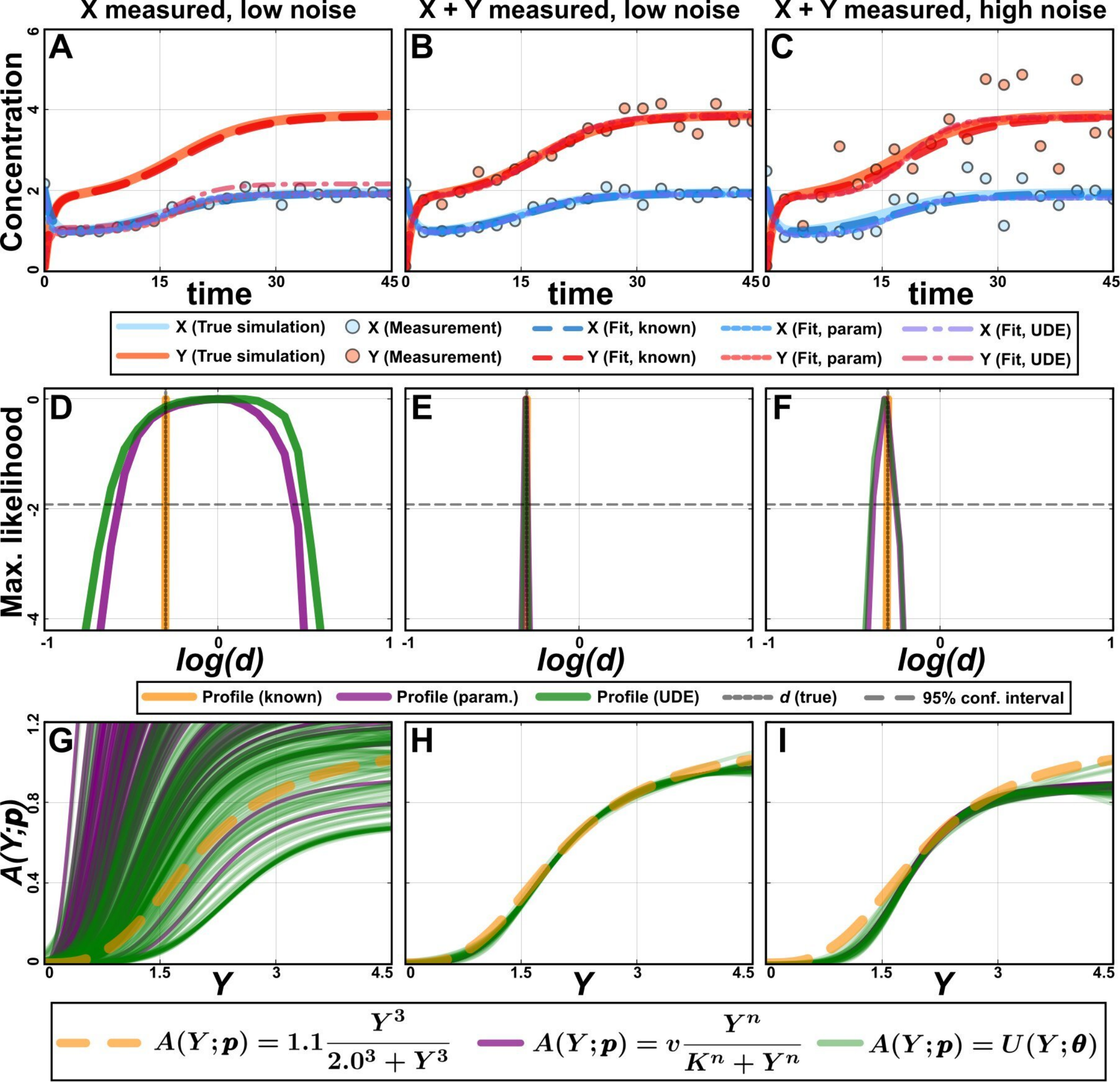}
    \caption{\textbf{Practical identifiability analysis of the extended self-activation loop: Repeat 2.} This figure is identical to Figure~\ref{figure:extended_sal_pract_ident} of the main text, however, here we use a different set of (equally distributed) measurements. This figure exhibits similar results and conclusions as Figure~\ref{figure:extended_sal_pract_ident}.}
    \label{sfigure:extended_sal_pract_ident_rep2}
\end{figure}

\begin{figure}
    \centering
    \includegraphics[width=\textwidth]{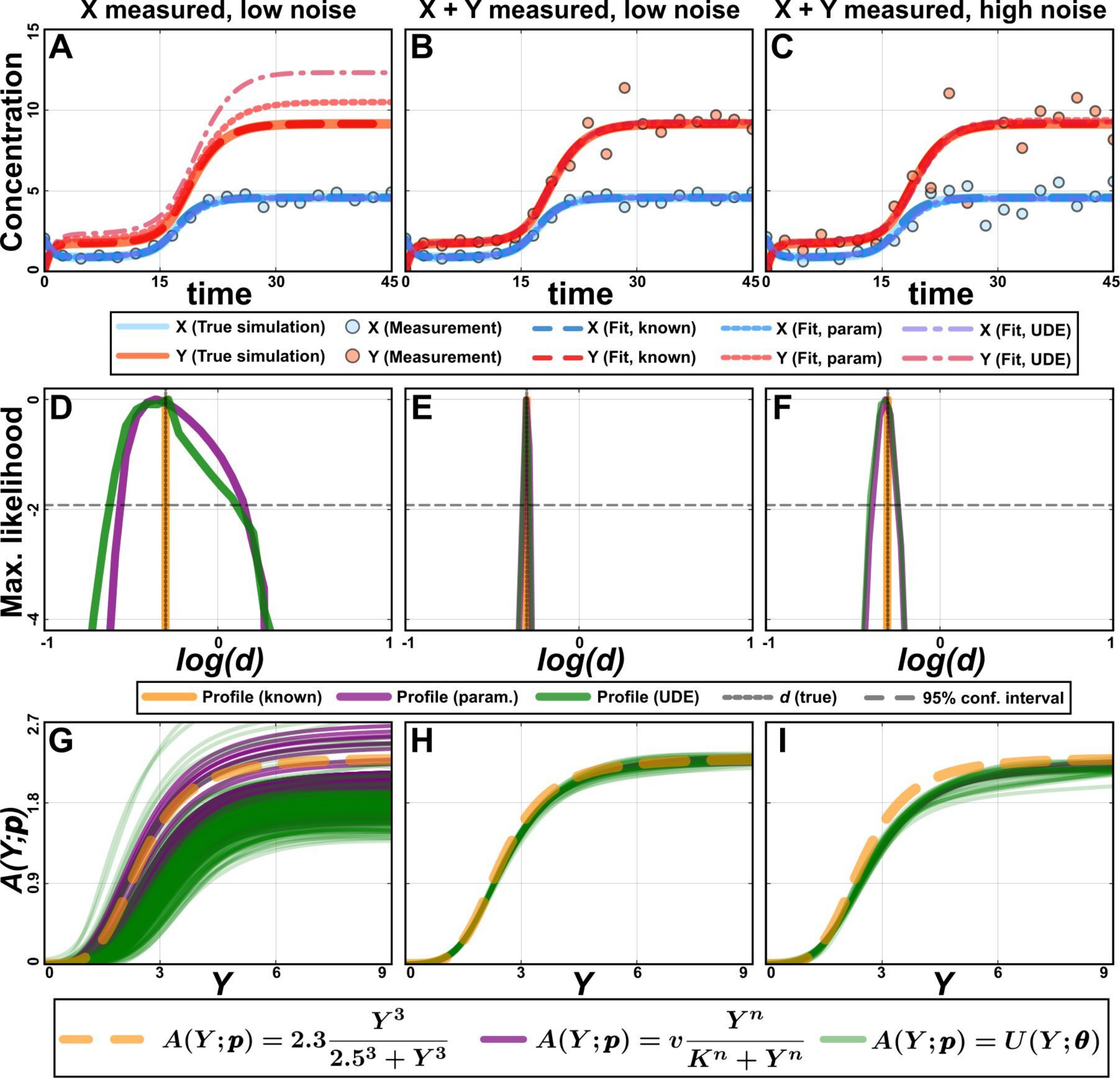}
    \caption{\textbf{Practical identifiability analysis of the extended self-activation loop: Repeat using a different parameter set.} This figure is identical to Figure~\ref{figure:extended_sal_pract_ident} of the main text, however, the data were generated from a different ground-truth parameter set (described in Supplementary Table~\ref{stable:data_settings_sfigs}). This figure exhibits similar results and conclusions as Figure~\ref{figure:extended_sal_pract_ident}.}
    \label{sfigure:extended_sal_pract_ident_alt_ps}
\end{figure}


\begin{figure}
    \centering
    \includegraphics[width=\textwidth]{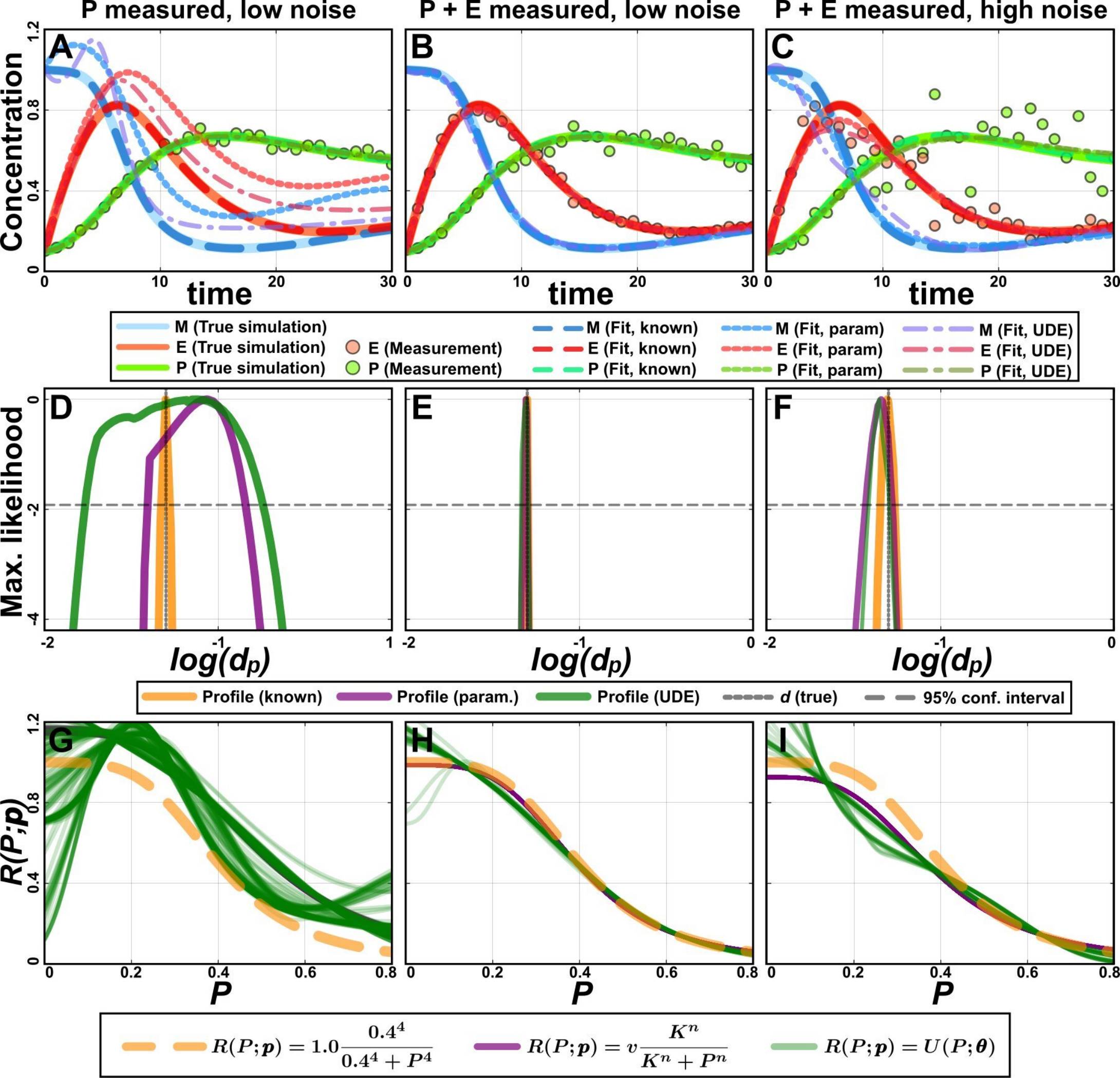}
    \caption{\textbf{Practical identifiability analysis: Repeat using the Goodwin oscillator model.} This figure is a repeat of Figure~\ref{figure:extended_sal_pract_ident}, but using the Goodwin oscillator model (described in Supplementary Section~\ref{ssection:gwo}). In A, only $P$ is measured, while in B and C, both $E$ and $P$ are measured. The model exhibits similar identifiability patterns as the extended self-activation loop model (Figure~\ref{figure:extended_sal_pract_ident}). Notable differences include that we can here confirm accurate predictions for a non-measured species ($M$) in the datasets where two species are measured (B and C). Furthermore, the case where the form of the unknown function, but not its parameters, is known, exhibits functional identifiability in G.}
    \label{sfigure:gwo_pract_ident}
\end{figure}


\begin{figure}
    \centering
    \includegraphics[width=\textwidth]{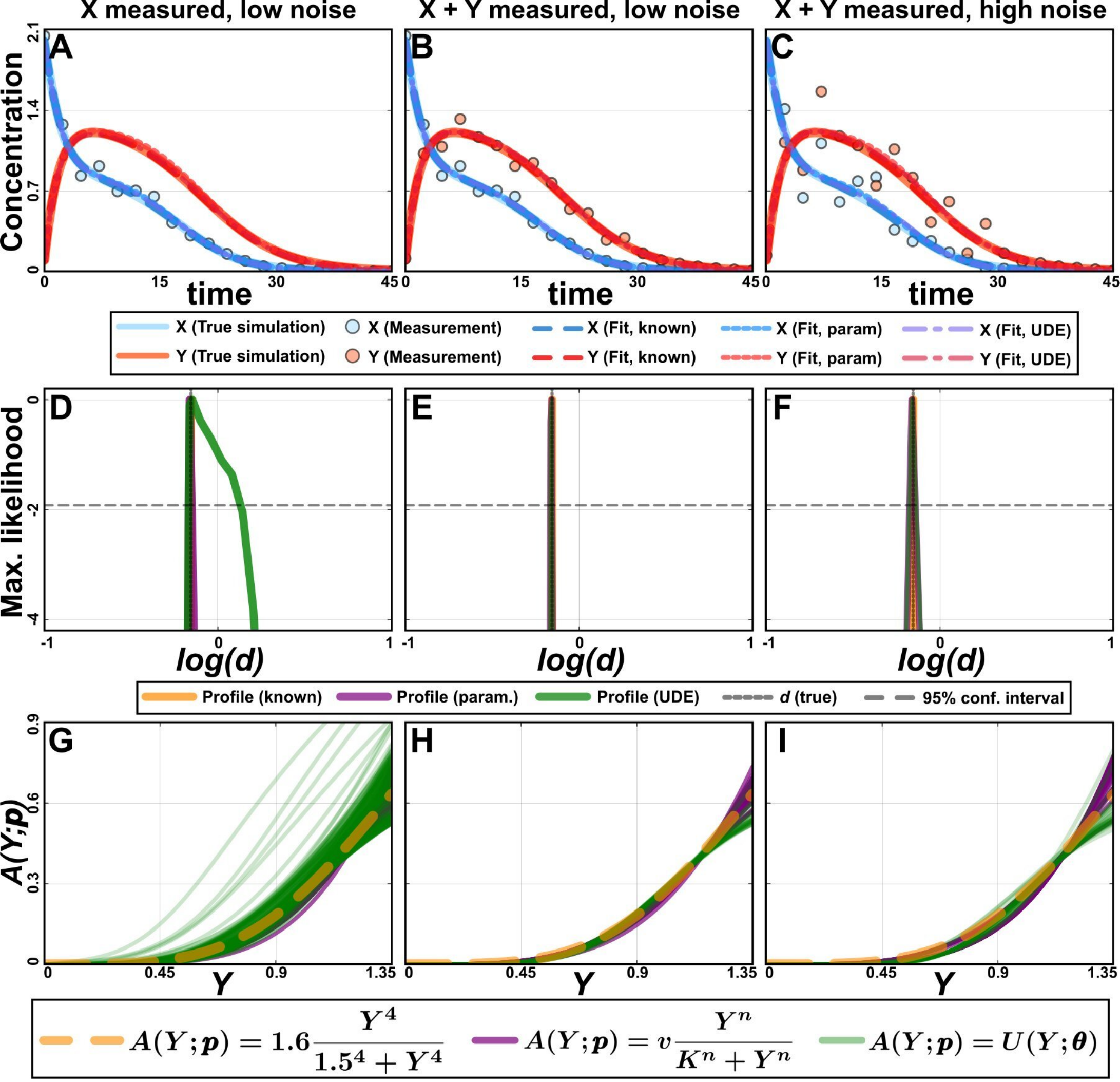}
    \caption{\textbf{For certain parameter sets, the identifiability of the extended self-activation loop changes.} This figure is identical to Figure~\ref{figure:extended_sal_pract_ident} of the main text, however, the data were generated from a different ground-truth parameter set and using a different initial condition (described in Supplementary Table~\ref{stable:data_settings_sfigs}). Notably, for this dataset and unlike the previous ones, the trained model only ``sees'' $Y$ values that are lower than the inflexion point of the Hill function (here at $K=1.5$). This means that we can only learn the true functions for this range of $Y$ values, i.e.~not the full sigmoidal shape. Furthermore, unlike previous cases, we here see a noteworthy difference between the parameterised Hill and UDE models. The former exhibits practical parametric and functional identifiability, while the latter exhibits some amount of practical non-identifiability (although not as much as in previous figures where only $X$ was measured).}
    \label{sfigure:extended_sal_pract_ident_lowconc_ps}
\end{figure}


\begin{figure}
    \centering
    \includegraphics[width=\textwidth]{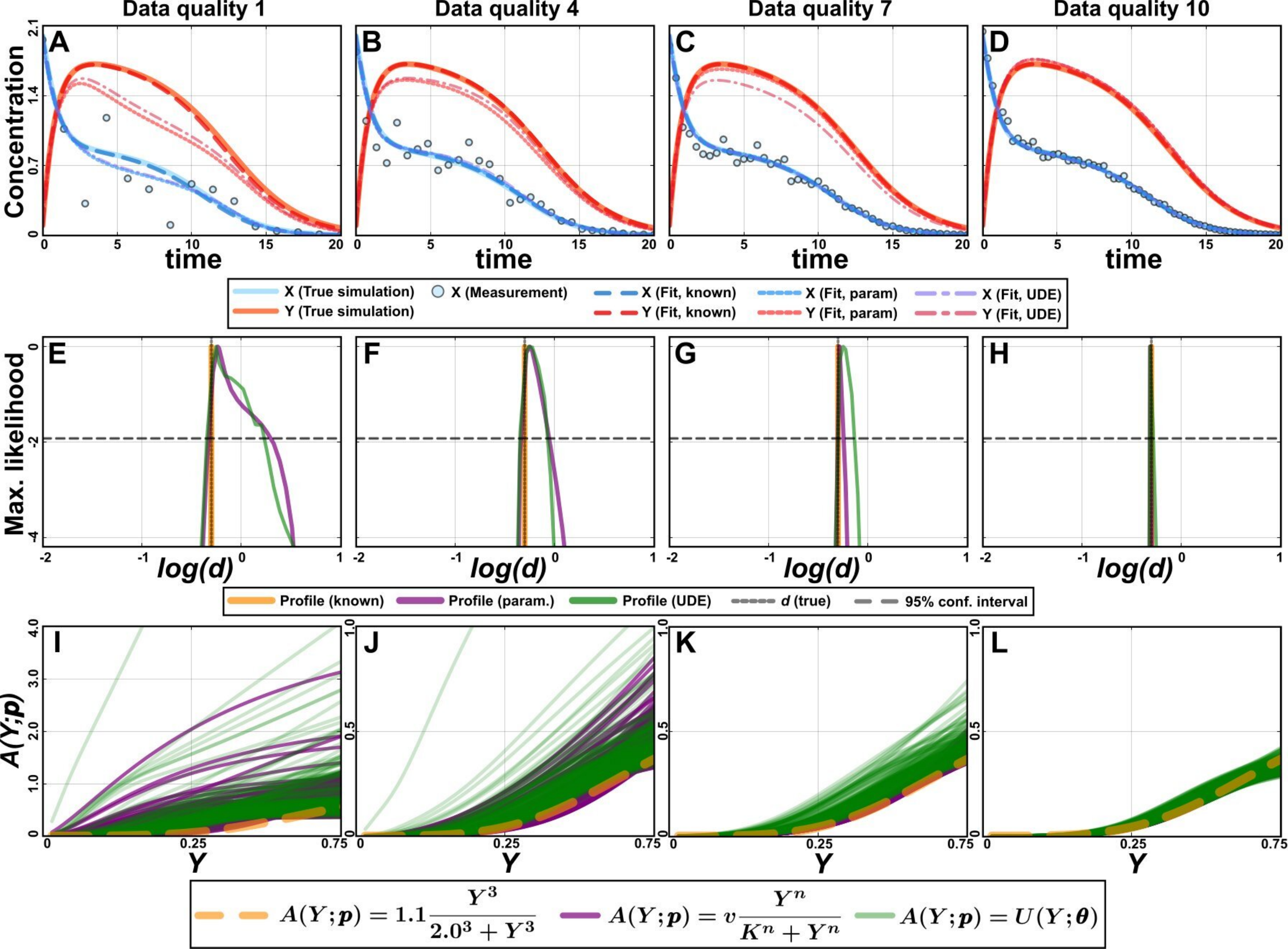}
    \caption{\textbf{Extended self-activation loop data quality scan example fitted simulations, profile likelihoods, and ensemble plots.} Supplementary Figure~\ref{sfigure:extended_sal_datavar} shows the relationship between practical parametric identifiability, practical functional identifiability, and predictive power. Here, we show example fitted simulations (A--D), likelihood profiles (E--H), and function ensemble plots (I--L). Note that the $y$-axis in I is on a different scale as compared to those in J-L. The leftmost column (A,E,I) shows examples from the lowest data quality ($x$-axis value of 1 in Supplementary Figure~\ref{sfigure:extended_sal_datavar}). The following columns show examples from data qualities 4, 7, and 10, respectively.}
    \label{sfigure:extended_sal_sample_fits}
\end{figure}


\begin{figure}[htbp]
    \centering
    \includegraphics[width=\textwidth]{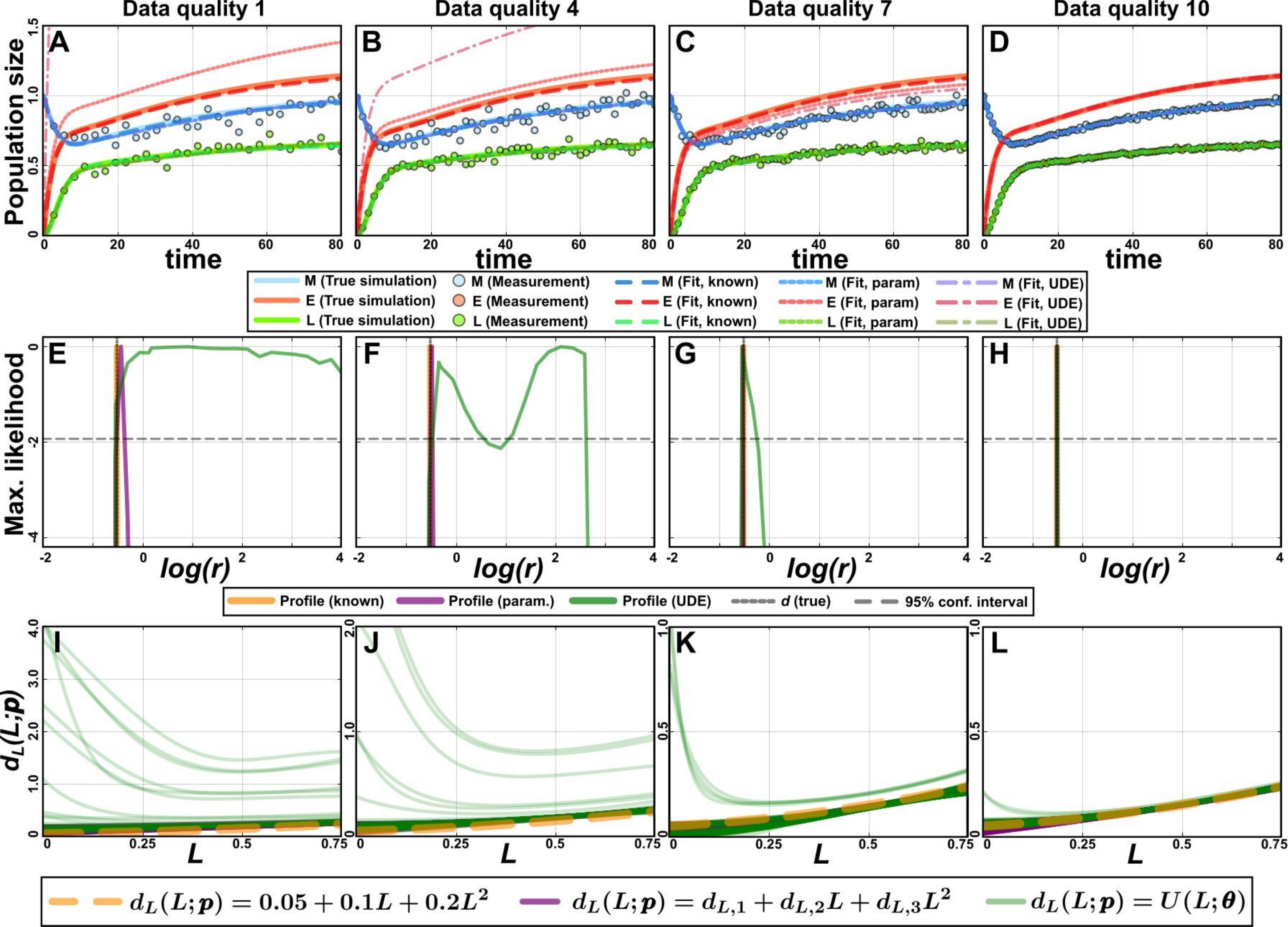}
    \caption{\textbf{Insect population data quality scan: example fitted simulations, profile likelihoods, and ensemble plots.} Figure~\ref{figure:mosquitoes_data_quality_scan} shows the relationship between parametric identifiability, functional identifiability, and predictive power. Here, we show example fitted simulations (A--D), likelihood profiles (E--H), and function ensemble plots (I--L). Note that the $y$-axis is changed across plots I to L. The leftmost column (A,E,I) shows examples from the lowest data quality ($x$-axis value of 1 in Figure~\ref{figure:mosquitoes_data_quality_scan}). The following columns show examples from data qualities 4, 7, and 10, respectively.}
    \label{sfigure:insects_example_fits}
\end{figure}


\begin{figure}
    \centering
    \includegraphics[width=\textwidth]{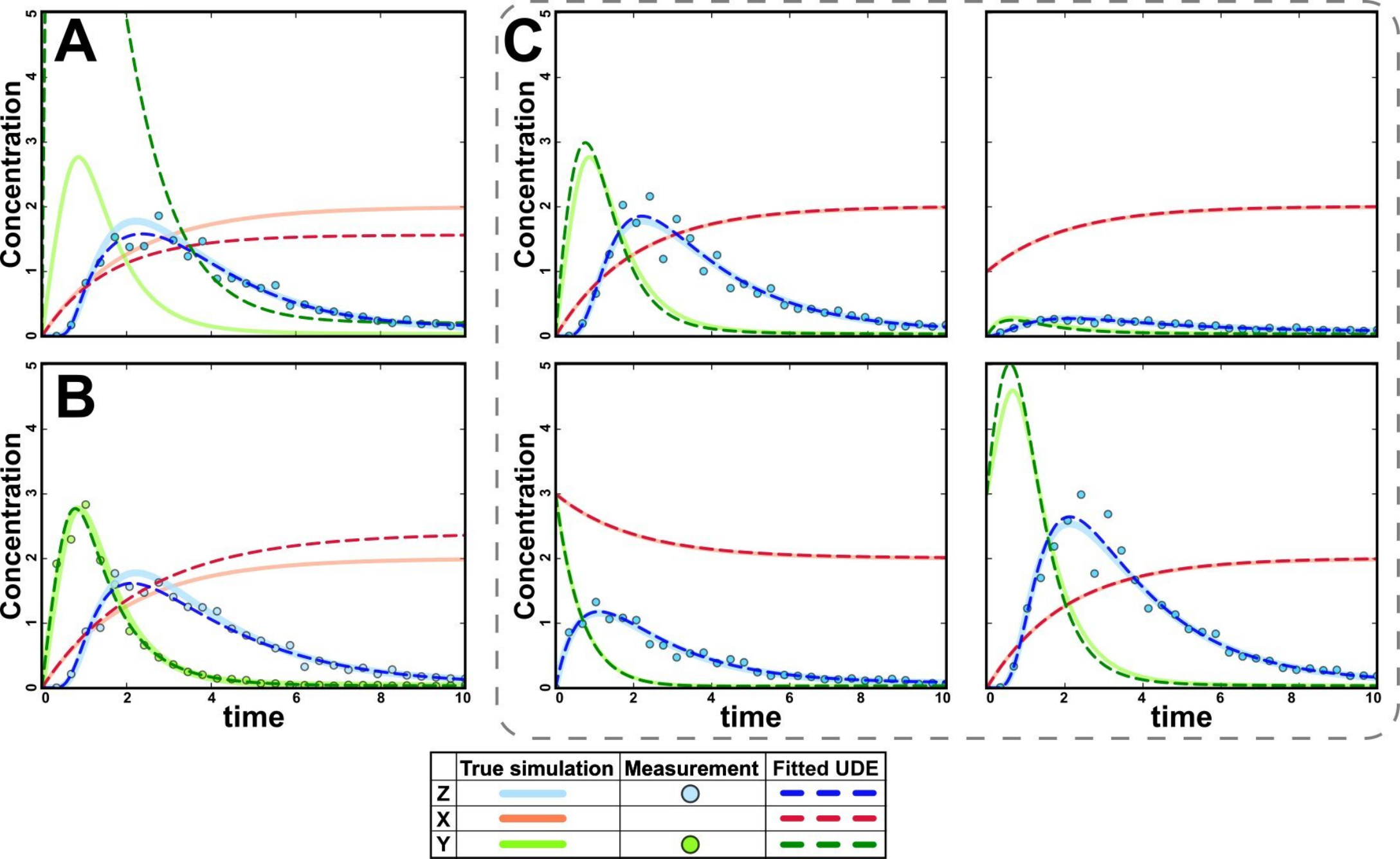}
    \caption{\textbf{The incoherent feedforward loop UDE can be successfully fitted to data from a single measured species if multiple experimental conditions are measured.} For different experimental settings, we fit a UDE of the incoherent feedforward loop (where both $R$ and $A$ are modelled using neural networks). (A) When fitting to measurements of $Z$ only, the prediction for $Y$ is wrong, and the prediction for $X$ is poor. (B) When measuring $Y$ as well, the trajectory aligns with $Y$. However, the prediction for $X$ remains incorrect. (C) While measuring $Z$ only, but using multiple (synthetic) experimental conditions, we can now successfully predict all species. This shows how the number of available experimental conditions is important for achieving a correct fit.}
    \label{sfigure:ic_feedforward_measurement_type_identifiability}
\end{figure}


\begin{figure}
    \centering
    \includegraphics[width=\textwidth]{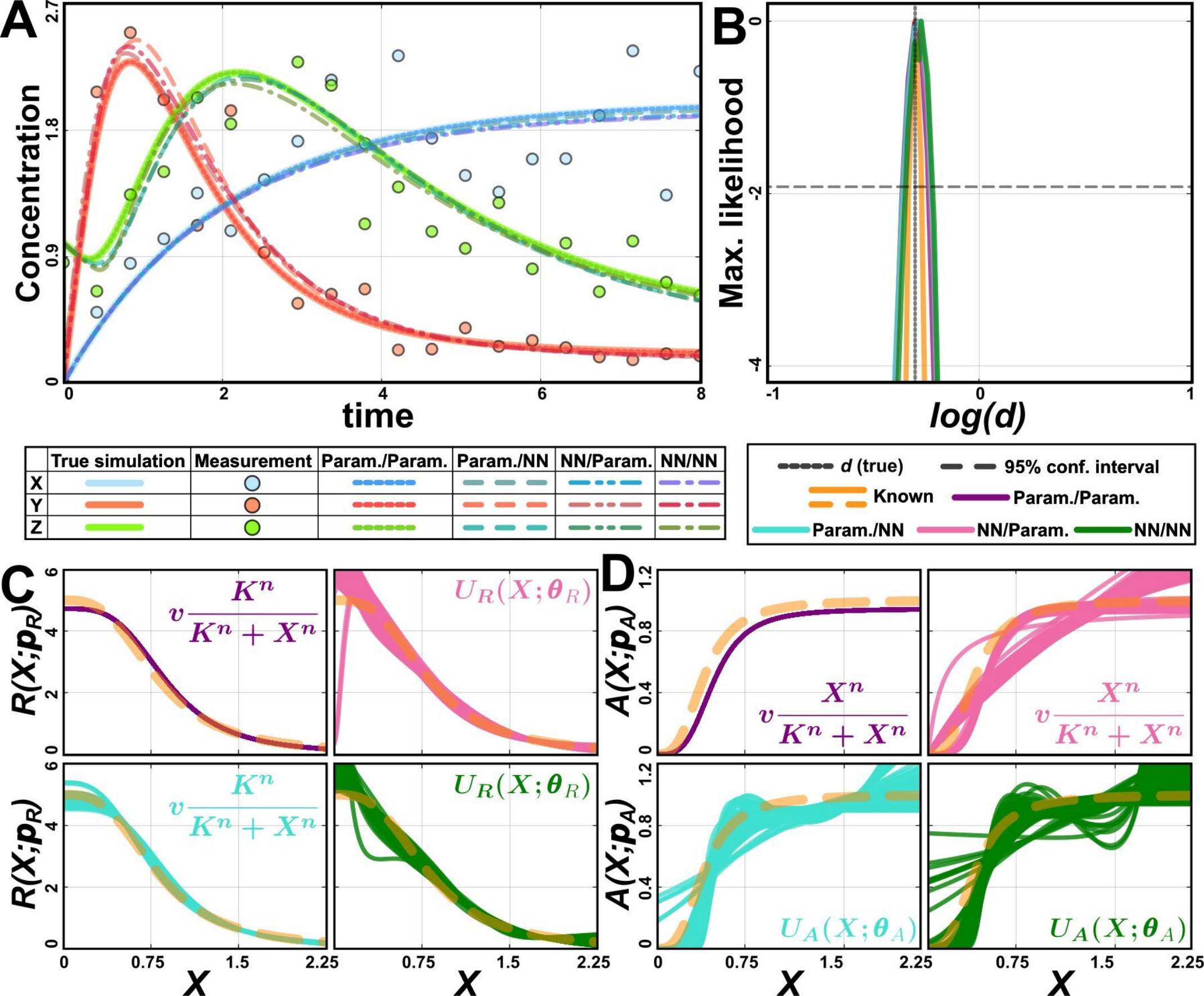}
    \caption{\textbf{The effect of multiple neural network approximation of functions in the incoherent feedforward loop model: Repeat.} This figure is identical to Figure~\ref{figure:ic_feedforward_2nn} of the main text, however, here we use a different set of (equally distributed) measurements. This figure exhibits similar results and conclusions as Figure~\ref{figure:ic_feedforward_2nn}.}
    \label{sfigure:ic_feedforward_2nn_repeat}
\end{figure}


\begin{figure}
    \centering
    \includegraphics[width=\textwidth]{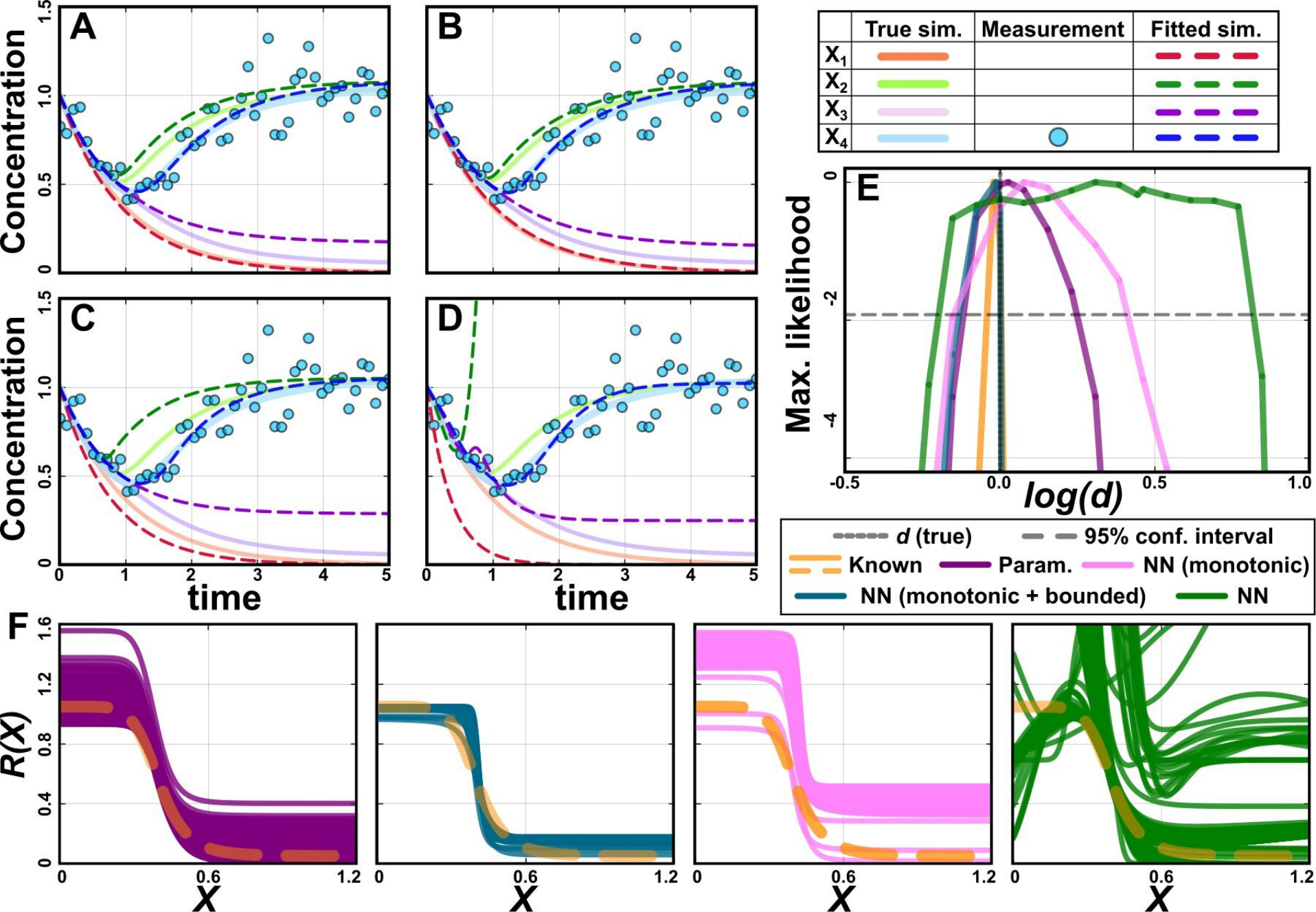}
    \caption{\textbf{The effect of neural network constraints on functional identifiability in the linear pathway model: Repeat 1.} This figure is identical to Figure~\ref{figure:ln_pathway_nn_constraints} of the main text, however, here we use a different set of (equally distributed) measurements. This figure exhibits similar results and conclusions as Figure~\ref{figure:ln_pathway_nn_constraints}.}
    \label{sfigure:ln_pathway_nn_constraints_repeat2}
\end{figure}

\begin{figure}
    \centering
    \includegraphics[width=\textwidth]{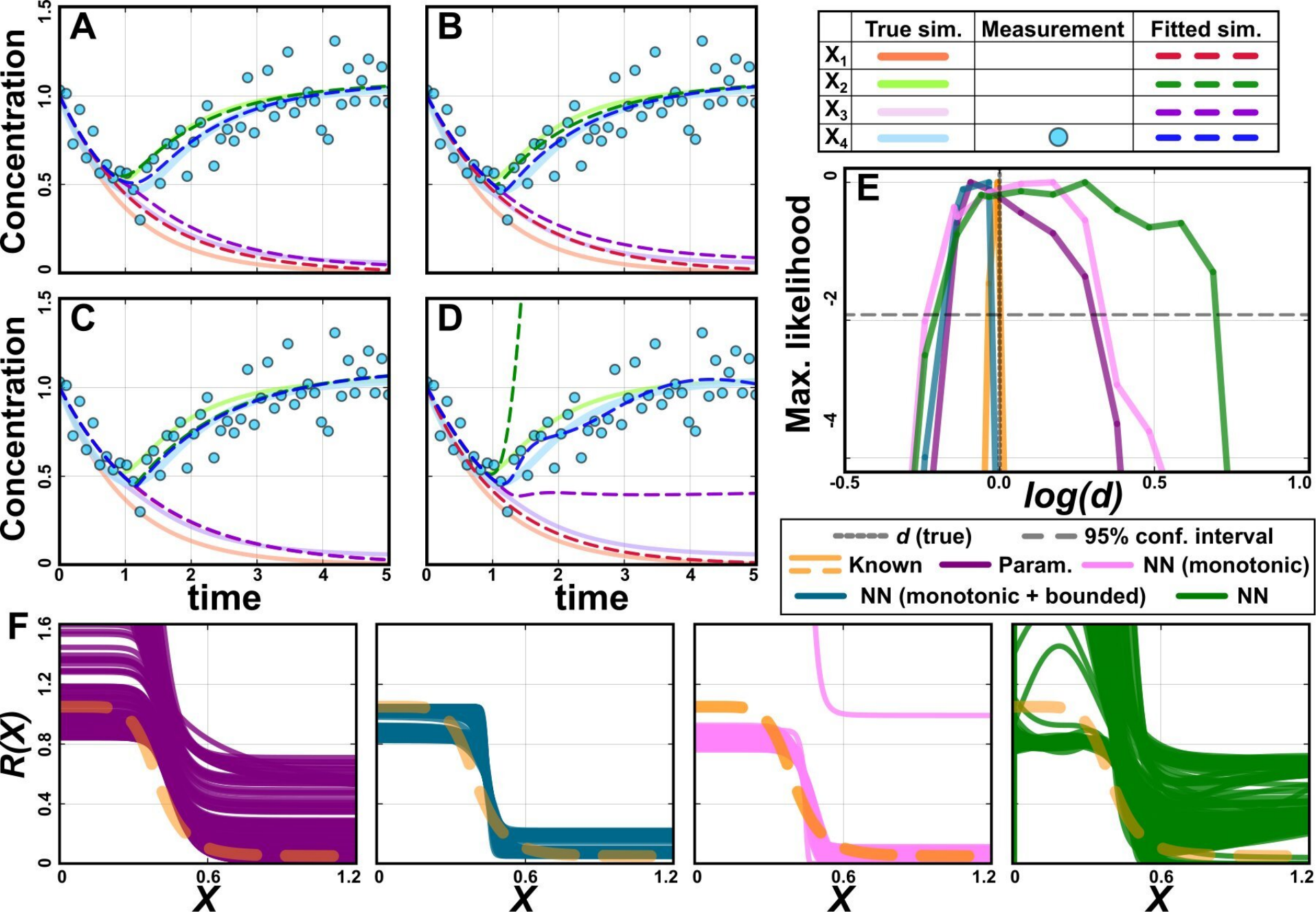}
    \caption{\textbf{The effect of neural network constraints on functional identifiability in the linear pathway model: Repeat 2.} This figure is identical to Figure~\ref{figure:ln_pathway_nn_constraints} of the main text, however, here we use a different set of (equally distributed) measurements. This figure exhibits similar results and conclusions as Figure~\ref{figure:ln_pathway_nn_constraints}.}
    \label{sfigure:ln_pathway_nn_constraints_repeat3}
\end{figure}



\begin{figure}
    \centering
    \includegraphics[width=\textwidth]{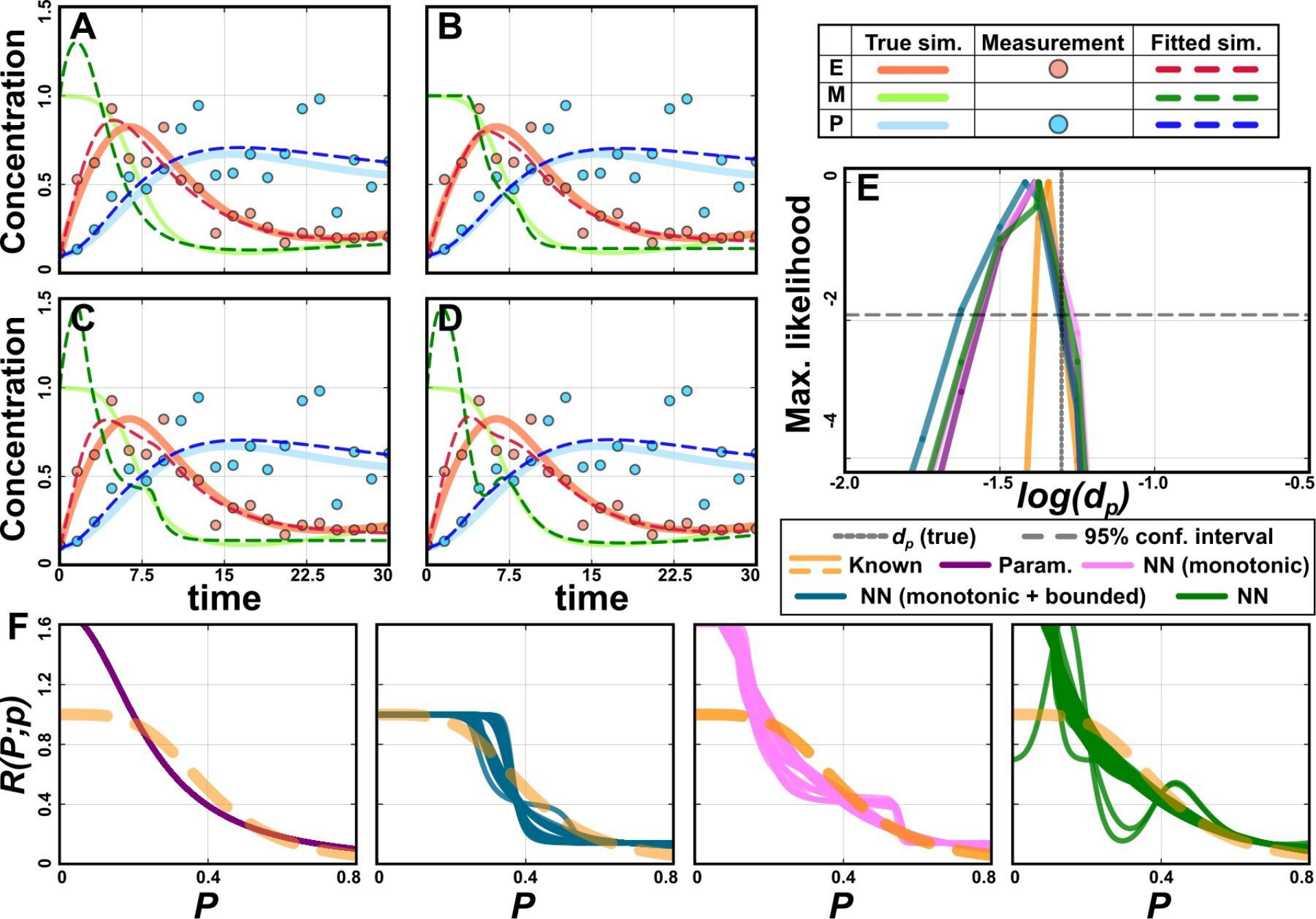}
    \caption{\textbf{The effect of neural network constraints on functional identifiability in the Goodwin oscillator model.}  This figure is a repeat of Figure~\ref{figure:ln_pathway_nn_constraints}, but using the Goodwin oscillator model (Supplementary Section~\ref{ssection:gwo}). Synthetic data are used to fit four model variants: (A) $R(P;\vec{p})$ known but $\vec{p}$ is unknown; (B) $R(P;\vec{p})$ is approximated as $U(P;\vec{\theta})$, where $U$ is constrained to be monotonically decreasing and bounded by $0 \leq U \leq v$; (C) $R(P;\vec{p})$ is approximated as $U(P;\vec{\theta})$, where $U$ is constrained to be monotonically decreasing; and (D) $R(P;\vec{p})$ is approximated as $U(P;\vec{\theta})$. (A--D) Fitted trajectories compared to ground-truth simulations. The monotonic and bounded neural network models produce good predictions for the unobserved species $M$. Next, the fully parameterised, monotonic neural network, and unconstrained neural network produce successively worse predictions. (E) Profile likelihoods for the fitted decay parameter $d_p$, including for the fully specified model as a reference. Here, all models exhibit similar identifiability, which is notably worse than for the fully specified reference. This differs from the linear pathway model (Figure~\ref{figure:ln_pathway_nn_constraints}), where the fully constrained model exhibits much better parameter identifiability as compared to the other models. (F) Ensembles of fitted functions $R$. The monotonic and bounded model yields functional identifiability, the unconstrained model exhibits non-identifiability. The remaining two models exhibit functional identifiability, but fail to converge to the correct values for low $P$. These results are reproduced in Supplementary Figure~\ref{sfigure:gwo_nn_constraints_repeat2} (where, notably, the predictions of the fully parameterised model are more similar to those of the monotonic and bounded neural network model).}    
    \label{sfigure:gwo_nn_constraints_repeat1}
\end{figure}

\begin{figure}
    \centering
    \includegraphics[width=\textwidth]{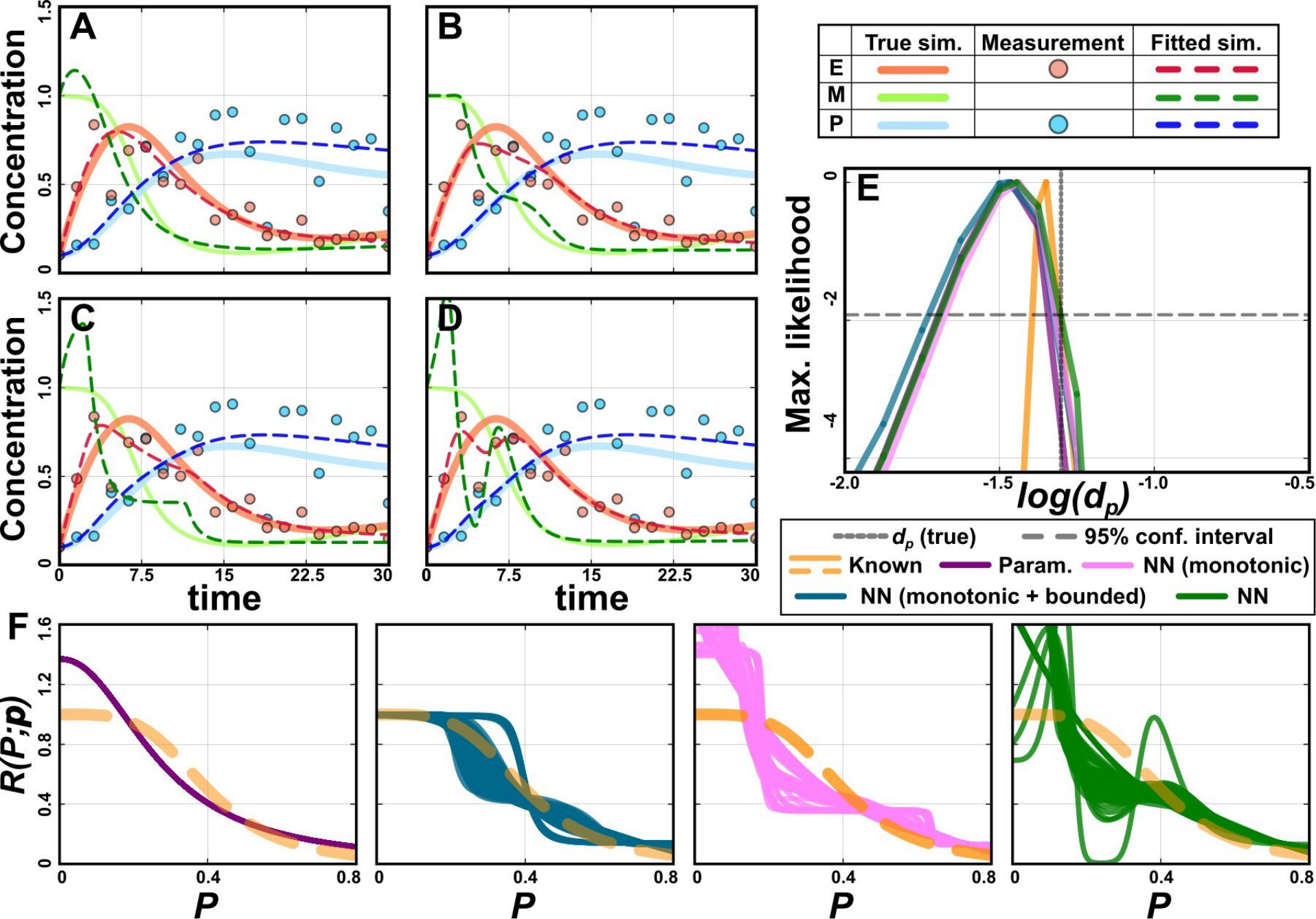}
    \caption{\textbf{The effect of neural network constraints on functional identifiability in the Goodwin oscillator model: Repeat.} This figure is identical to Figure~\ref{sfigure:gwo_nn_constraints_repeat1}, however, here we use a different set of (equally distributed) measurements. This figure exhibits similar results and conclusions as Supplementary Figure~\ref{sfigure:gwo_nn_constraints_repeat1} (with the notable exception that here the predictions of the fully parameterised model are more similar to those of the monotonic and bounded neural network model).}
    \label{sfigure:gwo_nn_constraints_repeat2}
\end{figure}


\begin{figure}
    \centering
    \includegraphics[width=\textwidth]{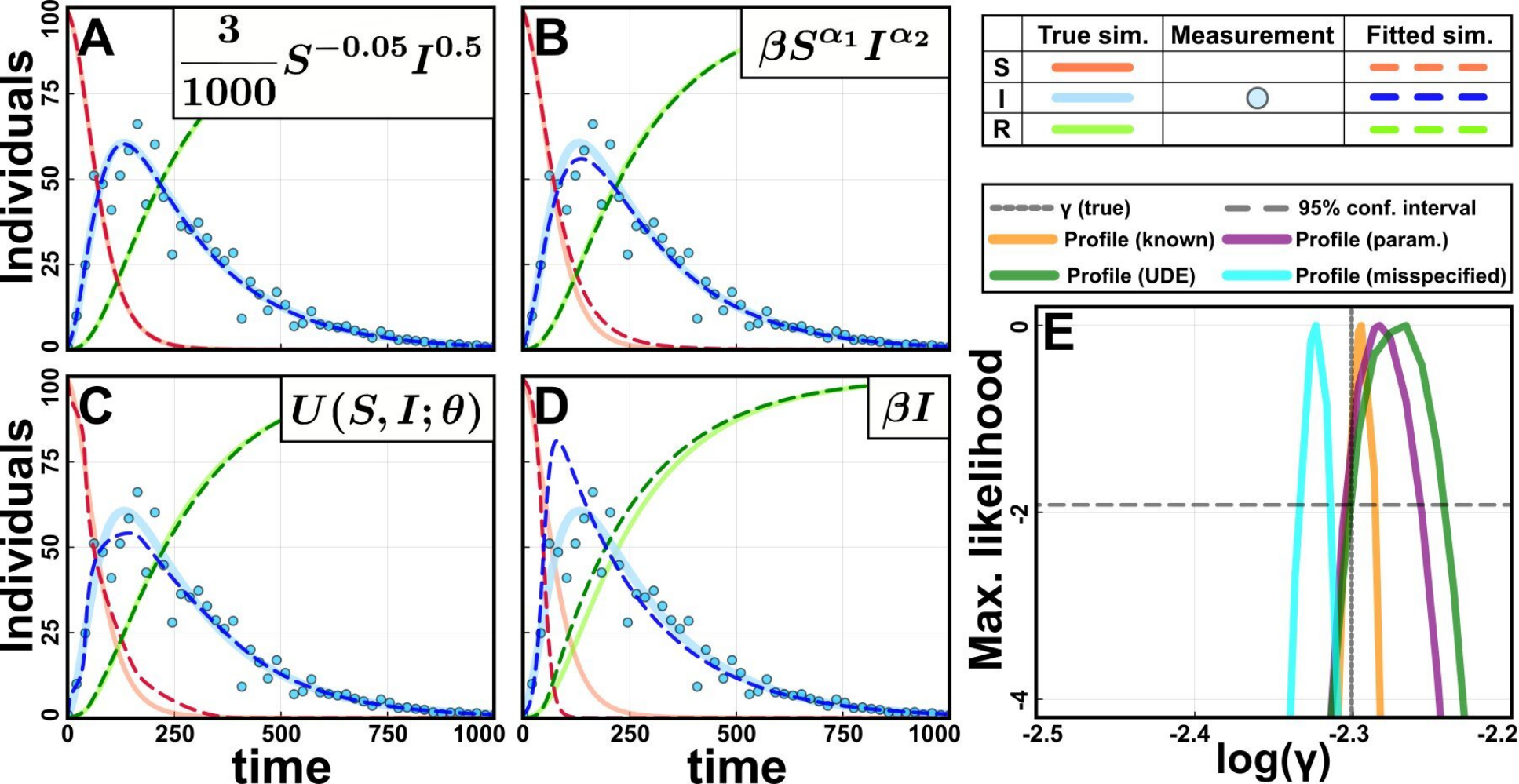}
    \caption{\textbf{A flexible infection-rate representation mitigates functional-form misspecification with little loss of practical parameter identifiability: Repeat.} This figure is identical to Figure~\ref{figure:modified_sir_misspecification} of the main text, however, here we use a different set of (equally distributed) measurements.}
    \label{sfigure:modified_sir_misspecification_repeat}
\end{figure}

\clearpage


\begin{figure}
    \centering
    \includegraphics[width=\textwidth]{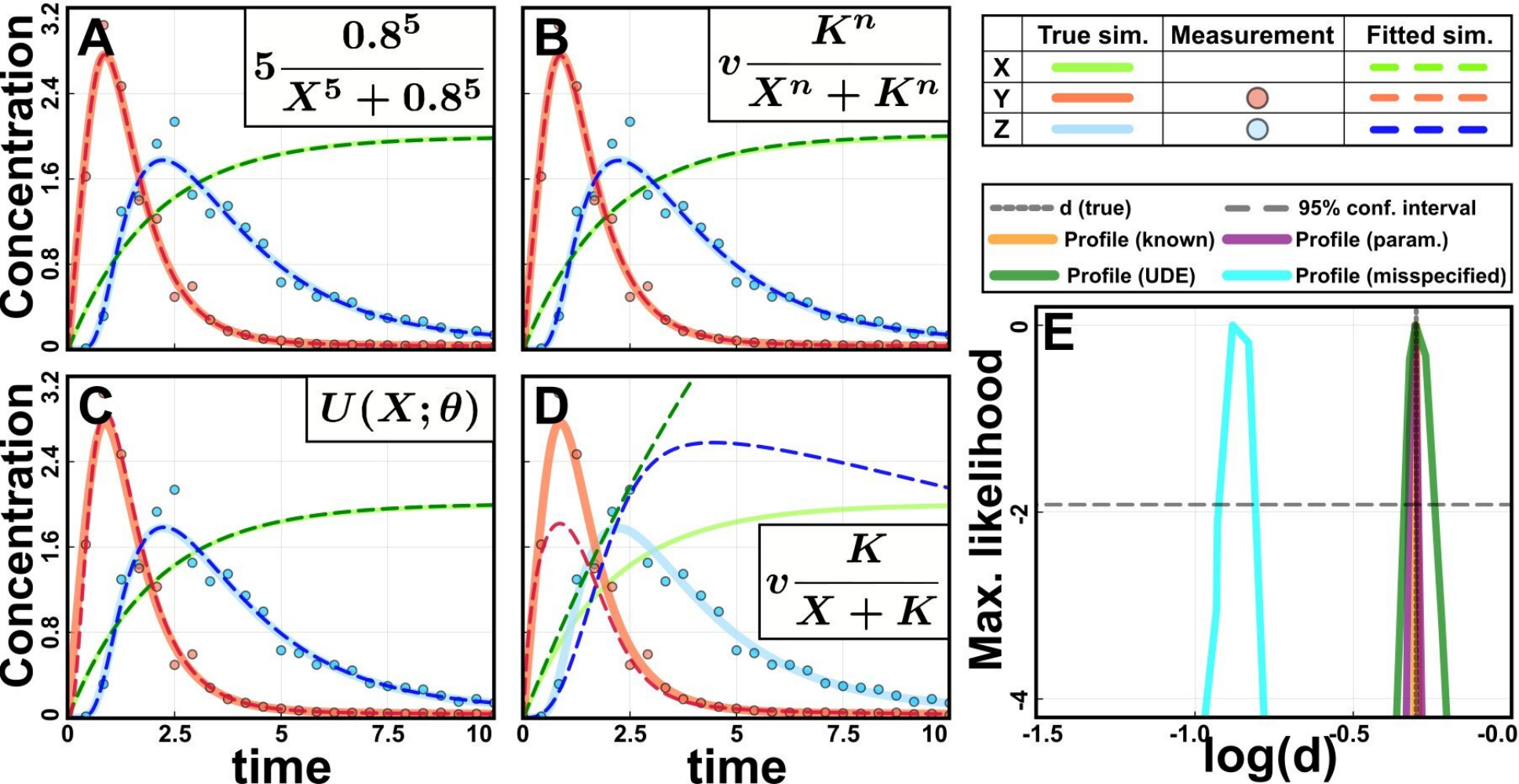}
    \caption{\textbf{The effect of function misspecification: Repeat using the incoherent feedforward loop model.}  This figure is a repeat of Figure~\ref{figure:modified_sir_misspecification}, but using the incoherent feedforward loop model (Section~\ref{section:results_ic_feedforward}). We use four versions of the model. In all versions, we assume that the activation function, $A(X;\vec{p}_A)$, but not its parameters ($\vec{p}_A$) is known. Next, we assume that the repressive function, $R(X;\vec{p}_R)$, is either: fully known (A), known while $\vec{p}_R$ are unknown (B), approximated as $U(X;\vec{\theta})$ (C), or misspecified as a repressive Michaelis--Menten function, $R(X;\vec{p}_R) = v\frac{K}{X + K}$, where $v$ and $K$ are here unknown parameters fitted to the data (D). (A--D) The first three models yield good fits to the data. The fourth model, however, fits the data poorly. (E) Likelihood profiles for the four models. The misspecified model predicts a somewhat wide interval of potential values that are off from the true value by almost an order of magnitude. The remaining three models yield good predictions of the true parameter value.}
    \label{sfigure:ic_feedforward_misspecification}
\end{figure}


\end{document}